\documentclass[final,3p]{elsarticle}
\usepackage{color}
\usepackage{amssymb}
\usepackage{amsmath}
\usepackage{booktabs}
\usepackage{multirow}
\usepackage{tabularx}
\usepackage{tabu}
\usepackage{stmaryrd}
\usepackage{epstopdf}
\usepackage{subfigure} 
\usepackage{graphicx}
\usepackage{caption}
\usepackage{amsthm} 
\usepackage{amssymb} 
\usepackage{mathrsfs}

\newcommand{\mbf}{\mathbf}

\numberwithin{equation}{section}
\newtheorem{thm}{Theorem}[section]

\newtheorem{rem}[thm]{Remark}
\newtheorem{ex}[thm]{Example}
\allowdisplaybreaks[4]

\biboptions{}
\journal{*}
\begin{document}
	
\begin{frontmatter}
\title{ Three kinds of novel multi-symplectic methods for stochastic Hamiltonian partial differential equations }

\author[ad1,ad2]{Jialin Hong}
\ead{hjl@lsec.cc.ac.cn} 

\author[ad1]{Baohui Hou\corref{cof1}}
\ead{houbaohui@lsec.cc.ac.cn} 

\author[ad1]{Qiang Li}
\ead{liqiang@amss.ac.cn} 

\author[ad1]{Liying Sun}
\ead{liyingsun@lsec.cc.ac.cn}

\cortext[cof1]{Corresponding author}
\address[ad1]{Institute of Computational Mathematics and Scientific/Engineering Computing, Academy of Mathematics and Systems Science, Chinese Academy of Sciences, Beijing 100190, China}
\address[ad2]{School of Mathematical Sciences, University of Chinese Academy of Sciences, Beijing 100049, China}
\begin{abstract}	
	
Stochastic Hamiltonian partial differential equations, which possess the multi-symplectic conservation law, are an important and fairly large class of systems. 
The multi-symplectic methods inheriting the geometric features of stochastic Hamiltonian partial differential equations provide  numerical approximations with better numerical stability, and are of vital significance  for obtaining correct
numerical results. 
 In this paper, we propose three novel multi-symplectic methods for stochastic Hamiltonian partial differential equations based on  the local radial basis function collocation method, the splitting technique, and the partitioned Runge--Kutta method. 
Concrete numerical methods are presented for nonlinear stochastic wave equations, stochastic nonlinear Schr\"odinger equations, stochastic Korteweg-de Vries  equations and stochastic Maxwell equations.  
We take stochastic wave equations as examples to perform numerical experiments, which indicate the validity of the proposed methods.
\end{abstract}		

\begin{keyword}
	stochastic Hamiltonian partial differential equations \sep multi-symplecticity  \sep  
	local radial basis function collocation method \sep splitting technique \sep partitioned Runge--Kutta method
\end{keyword}		

\end{frontmatter}


\section{Introduction}

A common way to describe the physical and engineering phenomena in the area of fluid dynamics, nonlinear optics, and quantum field theory (see e.g., \cite{Cao, Jiang, Roach, Song} and references therein) is by means of 
stochastic Hamiltonian partial differential equations (PDEs). 
Stochastic Hamiltonian PDEs,   which include stochastic wave equations,  stochastic Schr\"odinger equations, stochastic Korteweg-de Vries (KdV) equations, stochastic Maxwell equations, etc., are proposed in \cite{Chen, Hong1,Jiang}, and they have a prominent characteristic, that is,  multi-symplectic conservation law. 
The multi-symplecticity is the concatenation of differential 2-forms in both space and time,  decomposes
neatly the different facets of the governing equation, and characterizes the geometric invariants of the solution manifold.
Theoretical results concerning such multi-symplecticity reformulation can be found in \cite{Hong1,Hong2,Jiang} and references therein. 

 When designing a numerical method, a basic principle is that it should inherit the intrinsic properties of the original system as much as possible.  
Numerical methods that are incorporated more physical and geometric properties, especially multi-symplectic methods admitting the discrete multi-symplectic conservation law, have remarkable numerical superiority to conventional numerical methods.   
Recently, multi-symplectic methods possessing good
performance in preserving local conservation laws of the original system have been developed in the field of stochastic geometric integration of stochastic Hamiltonian PDEs (see e.g., \cite{Chen,Cui, Hong1,Hong2,Hong3} and references therein). 
For instance, \cite{Cui, Hong3} present multi-symplectic methods for stochastic nonlinear Schr\"odinger equations by making use of the central finite difference method in spatial direction combined with the midpoint method in temporal direction.   \cite{Hong2} proposes a 
multi-symplectic energy-conserving method, based on the wavelet collocation method in space and the symplectic method in time, for a three-dimensional stochastic Maxwell equation with multiplicative noise. 
 \cite{Chen, Hong1} investigate multi-symplectic methods for the stochastic Maxwell equation with additive noise via the implicit midpoint method and the leapfrog method.    
To the best of our knowledge, there is few work on the study of constructions of multi-symplectic methods for the general stochastic Hamiltonian PDEs.  
The first attempt to show the multi-symplectic method for the general 1-dimensional stochastic Hamiltonian PDEs is given in \cite{ZhangJ}, which takes advantage of symplectic Runge--Kutta methods with two Butcher tableaux. 
Our results in this paper not only present multi-symplectic partitioned Runge--Kutta methods with more Butcher tableaux, which increase diversity and flexibility of numerical methods, but also propose another two multi-symplectic methods via the local radial basis function (LRBF) collocation method and the splitting technique for the general stochastic Hamiltonian PDEs. 

Inspired by the fact that the LRBF collocation method has been successfully utilized  to numerically solve deterministic Hamiltonian PDEs, we apply the LRBF collocation method in space and midpoint method in time to derive the first kind of multi-symplectic method, that is, the messless LRBF collocation midpoint method shown in Section 3. 
The method is made on the overlapping sub-domains, which significantly reduces the size of the collocation matrix at the cost of solving lots of small matrices, and thus leads to cost efficiency.  
Moreover, it performs stably, can deal with complicated irregular domains and moving boundary, and possesses a long-time tracking capability. 
The second strategy of constructing the multi-symplectic method is utilizing the splitting technique allowing one to deal with sequentially a deterministic Hamiltonian PDE and a stochastic system, which are simpler than the original equation. 
For the numerical study of deterministic Hamiltonian PDEs, a lot of reliable and efficient numerical methods preserving the multi-symplecticity have been given (see \cite{Cohen, Hong4, McLachlan}). 
In Section 4, we first adopt the multi-symplectic Runge--Kutta method, that is, a derivative-free numerical method, to approximate deterministic Hamiltonian PDEs. 
Then combining the symplectic Euler method applied to the stochastic system, we arrive at the second kind of multi-symplectic method, that is, the splitting multi-symplectic Runge--Kutta method.  
We would like to mention that the splitting method does not need to handle the interaction between the nonlinear potential and the driving stochastic process. 
 In Section 5, we propose the third kind of multi-symplectic method for four stochastic Hamiltonian PDEs, i.e., stochastic wave equation,  stochastic nonlinear Schr\"odinger equation, stochastic KdV equation, and stochastic Maxwell equation, by employing the partitioned Runge--Kutta method in both temporal and spatial directions. 
The resulting method maybe explicit for some stochastic Hamiltonian PDEs. For instance, the method based on the symplectic Euler method in both space and time is explicit for the stochastic wave equation.  

The paper is organized as follows. In Section 2, we introduce the multi-symplectic  conservation law and the definition of multi-symplectic method for stochastic Hamiltonian PDEs. 
Section 3 presents the first kind of multi-symplectic method, which is constructed by the meshless LRBF collocation method and the midpoint  method.  
Section 4 is devoted to the second kind of multi-symplectic Runge--Kutta method based on the splitting technique and symplectic Runge--Kutta method.  
In Section 5, we apply the partitioned Runge--Kutta method to deriving the third kind of multi-symplectic method. 
We take stochastic wave equations as examples to perform numerical experiments, which indicate the validity of the proposed methods.
Finally, we give a conclusion in Section 6.

\section{Stochastic Hamiltonian PDEs}

Stochastic Hamiltonian PDEs, as natural extensions of stochastic 
Hamiltonian ordinary differential equations,  play important roles in the fields of fluid dynamics, nonlinear optics, plasma physics, communications and medical science and so forth.
They  are due to \cite{Jiang} and given by
\begin{equation} \label{eq2.1}
M d z + K z_{x} d t = \nabla S_{1}(z) d t+\nabla S_{2}(z) \circ dW(t),
\end{equation}
where $x\in\mathcal O,$ $M$ and $K$ are skew-symmetric matrices,  $S_{1}$ and $S_{2}$ are smooth functions of the variable $z$, and `$\circ$' stands for the Stratonovich product.
Moreover, $\{W(t)\}_{t\geq 0}$ is an $\mathbb L^2(\mathcal O,\mathbb R)$-valued $Q$-Wiener process with respect to a normal filtration $\{\mathcal{F}_t\}_{t\geq 0}$ on a filtered probability space $(\Omega, \mathcal{F},\{\mathcal{F}_t\}_{t\geq 0}, \mathbb{P})$ and has the expansion form
$$W(t) = \sum_{k=1}^{\infty} \sqrt{q_k}e_k\beta_k(t),$$ 
where $\{(q_k, e_k)\}_{k=1}^{\infty}$ is a sequence of eigenpairs of symmetric, positive definite and finite trace operator $Q$ with orthonormal eigenvectors and $\{\beta_k(t)\}_{k=1}^{\infty}$ is a sequence of real-valued mutually independent standard Brownian motions. 
Stochastic PDEs that can be rewritten as \eqref{eq2.1}, include and are not limited to nonlinear stochastic wave equation, stochastic nonlinear Schr\"odinger equation, stochastic  KdV equation, etc. 
More precisely, 

\begin{enumerate}
\item by introducing  $v = u_t$ and $w = u_x$,  we reformulate the nonlinear stochastic wave equation with homogenous Dirichlet boundary condition
	\begin{equation*} 
	du_{t} - u_{xx} dt +f(u)dt = g(u)\circ dW(t)
	\end{equation*}  into
	\begin{equation}\label{eq2.4}
	\left\{\begin{aligned}
	&d u = v dt , \\
	&u_ x =  w,\\
	&d v -  w_ x dt = -f(u)dt + g(u) \circ dW(t),
	\end{aligned}\right.
	\end{equation}
where $f: \mathbb L^2(\mathcal O ,\mathbb R) \rightarrow  \mathbb L^2(\mathcal O,\mathbb R)$ and $ g: \mathbb L^2(\mathcal O,\mathbb R) \rightarrow \mathscr{L}_2(\mathbb L^2(\mathcal O,\mathbb R),Q^{\frac{1}{2}}(\mathbb L^2(\mathcal O,\mathbb R)))$  satisfy the global Lipschitz continuous  condition 
with $\mathcal O=\left[x_{L}, x_{R}\right],$ $x_{L}, x_{R}\in\mathbb R,$ and  $\mathscr{L}_2(\mathbb L^2(\mathcal O,\mathbb R),Q^{\frac{1}{2}}(\mathbb L^2(\mathcal O,\mathbb R)))$ being the separable Hilbert space of Hilbert--Schmidt operators. 
Denoting $z=(u, p, v, w)^{\top},$  then \eqref{eq2.4} can be transformed into the multi-symplectic formulation 
	\begin{equation*}
	M d z+K z_{x} d t=\nabla S_{1}(z) d t+\nabla S_{2}(z) \circ dW(t)
	\end{equation*}
	with
	\begin{equation*}\begin{aligned}
	&M=\left(\begin{array}{cccc}
	0 & 0 & 1& 0 \\
	0 & 0 & 0& 0 \\
	-1 & 0 & 0& 0\\
	0 & 0 & 0& 0\\
	\end{array}\right), \quad
	K=\left(\begin{array}{cccc}
	0 & 0 & 0& -1 \\
	0 & 0 & 0 & 0\\
	0 & 0 & 0& 0\\
	1 & 0 & 0& 0\\
	\end{array}\right),\\
	& S_{1}(z)=\frac{1}{2}\left(w^{2}-v^{2}\right)-\tilde f(u), \quad S_{2}(z)=\tilde g(u),
	\end{aligned}\end{equation*}
	where $f=\tilde f_u$ and $g=\tilde g_u$ (see e.g., \cite{Song}).

	\item consider the stochastic nonlinear Schr\"odinger equation under the homogenous Dirichlet boundary condition
	\begin{equation*}
	\mathbf{i} du +  u_ {xx} dt + |u|^2 u dt = u\circ dW(t)
	\end{equation*}
	with $\mathcal O=[0,1]$ and $\mathbf{i}^2=-1$. Setting $u = p + \mathbf{i} q$  and letting $v=p_{x}$ and $w=q_{x}$, we rewrite the above equation as
	\begin{equation}\label{eq2.6}
	\left\{\begin{aligned}
	&d q - v_ {x} dt = \left(p^{2}+q^{2}\right)p dt  - p\circ dW(t),\\
	&d p +  w_ {x} dt = -\left(p^{2}+q^{2}\right)q dt + q\circ dW(t),\\
	&p_ {x} = v,\\
	&q_ {x} = w.
	\end{aligned}\right.
	\end{equation}
	Defining a state variable $z=(p, q, v, w)^{\top}$, \cite{Jiang} presents the associated multi-symplectic form of \eqref{eq2.6} as follows
	\begin{equation*}
	M d z+K z_{x} d t=\nabla S_{1}(z) d t+\nabla S_{2}(z) \circ dW(t)
	\end{equation*} 
	with
	\begin{equation*}\begin{aligned}
	&M=\left(\begin{array}{cccc}
	0 & -1 & 0 & 0\\
	1 & 0 & 0 & 0\\
	0 & 0 & 0 & 0\\
	0 & 0 & 0 & 0 \\
	\end{array}\right), \quad
	K=\left(\begin{array}{cccc}
	0 & 0 & 1 & 0\\
	0 & 0 & 0 & 1\\
	-1 & 0 & 0 & 0\\
	0 & -1 & 0 & 0 \\
	\end{array}\right),\\
	& S_{1}(z)=-\frac{1}{4}\left(p^{2}+q^{2}\right)^2-\frac{1}{2}\left(v^{2}+w^{2}\right), \quad S_{2}(z)=\frac{1}{2}\left(p^{2}+q^{2}\right).
	\end{aligned}\end{equation*}

	\item the stochastic  Korteweg-de Vries equation with the homogenous Dirichlet boundary condition takes the form 
	\begin{equation*}
	du + u u_ {x} dt + \beta u_ {xxx} dt = \lambda dW(t),
	\end{equation*}
	where  $\beta,\lambda>0$  and $\mathcal O=[0,1]$.
	Given new variables $v, \rho, w$ satisfying
	\begin{equation}\label{eq2.10}
	\left\{\begin{aligned}
	&-\frac{1}{2}  d\rho -\beta w_ {x} dt = \frac{1}{2}u^2 dt -v dt,\\
	&\rho_ {x}=u,\\
	&\frac{1}{2} du + v_ {x} dt = \lambda dW(t),\\
	&u_ {x} = w,
	\end{aligned}
	\right.
	\end{equation}
	we have the compact form (see \cite{Bouard}) 
	\begin{equation*}\label{eq2.11}
	M d z+K z_{x} d t=\nabla S_{1}(z) d t+\nabla S_{2}(z) \circ dW(t)
	\end{equation*}
	with  $z=(u, v, \rho, w)^{\top},$
	\begin{equation*}\begin{aligned}
	&M=\left(\begin{array}{cccc}
	0 & 0 & -\frac{1}{2} & 0\\
	0 & 0 & 0 & 0\\
	\frac{1}{2} & 0 & 0 & 0\\
	0 & 0 & 0 & 0 \\
	\end{array}\right), ~
	K=\left(\begin{array}{cccc}
	0 & 0 & 0 & - \beta\\
	0 & 0 &-1 & 0\\
	0 & 1 & 0 & 0\\
	\beta & 0 & 0 & 0 \\
	\end{array}\right),~ S_{1}(z)=\frac{1}{6}u^3-uv+\frac{1}{2}\beta w^2, \quad S_{2}(z)=\lambda \rho.
	\end{aligned}\end{equation*}

		\item take the stochastic  Maxwell equation with multiplicative noise
	\begin{equation}\label{Maxwell's equation}
	\left\{\begin{aligned}
	&d\mathbf E(t,x,y,z)=\nabla \times \mathbf H(t,x,y,z)-\lambda \mathbf H(t,x,y,z)\circ dW(t),\\
	&d\mathbf H(t,x,y,z)= -\nabla \times \mathbf E(t,x,y,z) + \lambda \mathbf E(t,x,y,z) \circ dW(t)
	\end{aligned}\right.
	\end{equation}
	into account,
	where $\lambda\in\mathbb R,$  $\mathcal O\subset \mathbb{R}^{3}$ is a bounded and simply connected domain with smooth boundary $\partial \mathcal O$. 
	We employ the perfectly electric conducting (PEC) boundary condition
	$
	\mathbf{E} \times \mathbf{n}=\mathbf{0}
	$
	on $(0, T] \times \partial \mathcal O$ with $\mathbf{n}$ being the unit outward normal of $\partial \mathcal O$ (see \cite{Hong2}). 
	Denote $\mathbf u=(\mathbf H^\top, \mathbf E^\top)^\top=\left(H_{1}, H_{2}, H_{3}, E_{1}, E_{2}, E_{3}\right)^{\top}$ and $S({\bf u})=\frac{\lambda}{2}\left(\left|E_{1}\right|^{2}+\left|E_{2}\right|^{2}+\left|E_{3}\right|^{2}+\left|H_{1}\right|^{2}+\left|H_{2}\right|^{2}+\left|H_{3}\right|^{2}\right).$ 
	Then \eqref{Maxwell's equation}  can be rewritten as 
	\begin{equation}\label{symMaxwell}
	M d\mathbf u+K_{1}\mathbf u_{x} d t+K_{2}\mathbf u_{y} d t+K_{3}\mathbf u_{z} d t=\nabla S({\bf u}) \circ d W,
	\end{equation}
	where
	$$
	M=\left(\begin{array}{cc}
	0 & -I_{3 \times 3} \\
	I_{3 \times 3} & 0
	\end{array}\right), \quad K_{i}=\left(\begin{array}{cc}
	\mathscr{D}_{i} & 0 \\
	0 & \mathscr{D}_{i}
	\end{array}\right), \quad i=1,2,3
	$$
	with $I_{3 \times 3}$ being a $3 \times 3$ identity matrix,
	$$
	\mathscr{D}_{1}=\left(\begin{array}{ccc}
	0 & 0 & 0 \\
	0 & 0 & -1 \\
	0 & 1 & 0
	\end{array}\right),\quad \mathscr{D}_{2}=\left(\begin{array}{ccc}
	0 & 0 & 1 \\
	0 & 0 & 0 \\
	-1 & 0 & 0
	\end{array}\right),\quad \mathscr{D}_{3}=\left(\begin{array}{ccc}
	0 & -1 & 0 \\
	1 & 0 & 0 \\
	0 & 0 & 0
	\end{array}\right) .
	$$		
\end{enumerate}

Analogues to the symplecticity-preserving property of stochastic Hamiltonian ordinary differential equations, \cite{Jiang} shows that stochastic Hamiltonian PDEs possess the multi-symplectic conservation law. 
In detail, denote two differential 2-forms by $\omega = \mathrm{d} z \wedge M \mathrm{d} z$ and $\kappa = \mathrm{d} z \wedge K \mathrm{d} z,$ where `$\wedge$' represents the wedge product.  
Then the multi-symplecticity, as a local invariant, is given by
\begin{equation} \label{eq2.2}
d \omega(t, x)+\partial_{x} \kappa(t, x)dt=0, \quad a.s.,
\end{equation}
i.e.,
\begin{equation*}
\int_{x_{0}}^{x_{1}} \omega\left(t_{1}, x\right) d x -\int_{x_{0}}^{x_{1}} \omega\left(t_{0}, x\right) d x+\int_{t_{0}}^{t_{1}} \kappa\left(t, x_{1}\right) d t-\int_{t_{0}}^{t_{1}} \kappa\left(t, x_{0}\right) d t =0,  \quad a.s.,
\end{equation*}
where $\left(x_{0}, x_{1}\right) \times\left(t_{0}, t_{1}\right)$ is the local definition domain of $z$. 
From the multi-symplectic conservation law \eqref{eq2.2} it can be found that symplecticity changes locally and synchronously both in temporal and  spatial directions. 
We would like to remark that  the word `local' means that such conservative property does not depend on the specific domain or on boundary conditions of  stochastic PDEs.  
In addition, the multi-symplectic conservation law \eqref{eq2.2} for stochastic Hamiltonian PDEs holds almost surely. 
To simplify the notation, below we shall suppress the notation `a.s.’ unless it is necessary to avoid confusion.
The multi-symplectic conservation law 
\begin{itemize}
\item for nonlinear stochastic wave equation  \eqref{eq2.4}  is
\begin{equation*}
d[\mathrm{d} u \wedge \mathrm{d} v]+\partial_{x}[\mathrm{d} w \wedge \mathrm{d} u]dt=0. 
\end{equation*}

\item for stochastic nonlinear Schr\"odinger equation \eqref{eq2.6} is
\begin{equation*}\label{eq2.8}
d[\mathrm{d} q \wedge \mathrm{d} p] +\partial_{x}[\mathrm{d} p \wedge \mathrm{d} v + \mathrm{d} q \wedge \mathrm{d} w]dt=0. 
\end{equation*}
\item for stochastic KdV equation \eqref{eq2.10} is
\begin{equation*}
d[\mathrm{d} \rho \wedge \mathrm{d} u] +\partial_{x}[2\mathrm{d} \rho \wedge \mathrm{d} v + 2\beta\mathrm{d} w \wedge \mathrm{d} u]dt = 0. 
\end{equation*}

\item for stochastic Maxwell equation  \eqref{Maxwell's equation} is
\begin{align*}
&d[\mathrm{d} \mathbf E \wedge \mathrm{d} \mathbf H] +\partial_{x}[\mathrm{d} H_3 \wedge \mathrm{d} H_2 + \mathrm{d} E_3 \wedge \mathrm{d} E_2]dt\\
&+\partial_{y}[\mathrm{d} H_1 \wedge \mathrm{d} H_3 + \mathrm{d} E_1 \wedge \mathrm{d} E_3]dt
+\partial_{z}[\mathrm{d} H_2 \wedge \mathrm{d} H_1 + \mathrm{d} E_2 \wedge \mathrm{d} E_1]dt=0.
\end{align*} 
\end{itemize}

In order to keep more intrinsic properties of the original system into numerical simulations, there has been growing interest in the geometric integration of stochastic Hamiltonian PDEs, namely in the multi-symplectic method, which can more fully capture behaviors of interesting phenomena. 
For the purpose of numerical approximation, we let $\Delta x$, $\Delta y$ and $\Delta z$ be the mesh sizes along $x, y$ and $z$ directions, respectively, and $\Delta t$ be the time step length. 
The temporal-spatial domain we are interested in the following sections is  $[0, T] \times \mathcal O:=[0, T] \times\left[x_{L}, x_{R}\right] \times\left[y_{L}, y_{R}\right] \times\left[z_{L}, z_{R}\right].$ 
It is partitioned by parallel lines, where $t_{n}=n \Delta t,$ $x_{i}=x_{L}+i \Delta x,$ $y_{j}=y_{L}+j \Delta y$ and $z_{k}=z_{L}+k \Delta z$ for $n=0,1, \dots, N,$ $i=0,1, \dots, I,$  $j=0,1, \dots, J$ and  $k=0,1, \dots, K$. 
Now we take $\mathcal O=\left[x_{L}, x_{R}\right]$ for example and denote the approximation of the $z(x, t)$ at the mesh point $(x_j, t_k)$ by $z_{j, k},$ i.e., $z_{j, k}\approx z(x_j, t_k)$.  
The numerical method for \eqref{eq2.1} and \eqref{eq2.2}, can be written, respectively, as
\begin{align}
\label{eq0.3}
\Delta t M \delta_{t}^{j, k} z_{j, k}+\Delta t K \delta_{x}^{j, k} z_{j, k} &=\Delta t(\nabla_{z} S_1(z))_{j,k}+\Delta W_{j}^k(\nabla_{z} S_2(z))_{j,k}, \\
\label{eq0.4}
\delta_{t}^{j, k} \omega_{j, k}+\delta_{x}^{j, k} \kappa_{j, k} &=0,
\end{align}
 where 
\begin{equation*} 
\omega_{j, k} = {\rm d} z_{j, k}\wedge M {\rm d} z_{j, k},~ \kappa_{j, k} = {\rm d} z_{j, k} \wedge K {\rm d} z_{j, k}, 
\end{equation*}
$\Delta W_{j}^k =  W({x_j, t_{k+1}}) -  W({x_j, t_k})$, and $\delta_{t}^{j, k}, \delta_{x}^{j, k}$ are corresponding discretizations of two partial derivatives $ \partial_{t}$ and $ \partial_{x}$, respectively.  
The numerical method \eqref{eq0.3} is called a multi-symplectic method  for stochastic Hamiltonian PDEs if  it satisfies a discrete version of the multi-symplectic conservation law \eqref{eq0.4}. 
In recent years, many researchers have studied various multi-symplectic methods for stochastic Maxwell equations (see e.g., \cite{Hong1,Hong2,ZhangJ}), stochastic nonlinear Schr\"odinger equations (see e.g., \cite{Cui, Hong3}), etc.

In what follows, we propose three  multi-symplectic methods of stochastic Hamiltonian PDEs. Soon afterwards, applications to nonlinear stochastic wave equation,  stochastic nonlinear Schr\"odinger equation, stochastic  KdV equation and stochastic Maxwell equation are given. 

\section{Meshless LRBF collocation midpoint method}
\label{Sec;MLRBF}
In this section, we present a kind of multi-symplectic methods for stochastic Hamiltonian PDEs  by exploiting the meshless LRBF collocation method in space and the midpoint method in time, respectively.

The global radial basis function collocation method, such as the Kansa’s method in \cite{Kansa1, Kansa2}, becomes a powerful tool for numerically solving  deterministic PDEs, especially deterministic Hamiltonian PDEs (see e.g., \cite{Dehghan, Wu} and references therein), since it does not need to evaluate any integral and has both high-order accuracy and geometric flexibility. 
A key ingredient of the global radial basis function collocation method is the radial basis function $\varphi$, such as the Gaussian radial basis function $\varphi(x) = e^{-c^2x^2}$, the multiquadric radial basis function $\varphi(x) = \sqrt{x^2 + c^2},$ and the inverse multiquadric  radial basis function $\varphi(x) =1/ \sqrt{x^2 + c^2}$, where the shape parameter $c$ is a constant.
However, when applying the global radial basis function collocation method  to solve PDEs, large scaled linear systems are needed to solve,  the corresponding coefficient matrices are ill-conditioned and the results are sensitive to the shape parameter $c$.  
To overcome the above problems arised by using the global radial basis function collocation method, the LRBF collocation method was formulated by \cite{Lee, Sarler}, from which the main idea is the collocation on influence domain, and can drastically reduce the collocation matrix size at the expense of solving many small matrices with the dimension of the number of nodes included in the domain of influence for each node.  
Since the LRBF collocation method, as a type of meshless methods, can be employed to cope with complex geometries, complicated irregular domains including moving boundary and high-dimensional problem, it has been applied for solving many problems in engineering and applied mathematics widely (see \cite{Zhang}). 
To be specific, let $\{\mathbf x_i, f(\mathbf x_i)\}$ be the scattered data with $i= 0, 1,\dots, L, L+1$, and $L\in \mathbb N.$ 
Fix $i\in\{0,1, \dots, L, L+1\}.$
Given $\mathbf{x}_{i},$ there exist $n_i$ neighboring nodes which are nearest from $\mathbf{x}_{i}$ in the influence domain $_{i} \Omega=\left\{_{i}\mathbf{x}_{k}\right\}_{k=1}^{n_{i}}$.  
For $_{i} \mathbf{x} \in{ }_{i} \Omega,$ the function $f$ can be approximated by
\begin{equation*}
f^*(_{i}\mathbf{x})=\sum_{k=1}^{n_{i}} {}_{i}\alpha_{k} \varphi\left(\left\|_{i} \mathbf{x}-{}_{i}\mathbf{x}_{k}\right\|\right),
\end{equation*}
where the coefficient $\{{}_i\alpha_k\}_{k=1}^{n_i}$ in the above equation satisfies the interpolation condition $f^*(_{i} \mathbf{x}) = f(_{i} \mathbf{x})$. 
Taking $_{i} \mathbf{x} = {}_{i}\mathbf{x}_{k}$ for $k=1, \dots, n_{i}$, we obtain
\begin{equation}
\begin{aligned}
\label{eq3.1}
{ }_{i}\mathbf{f} &=\left[\begin{array}{cccc}
\varphi\left(\left\|_{i} \mathbf{x}_{1}-{ }_{i} \mathbf{x}_{1}\right\|\right) & \varphi\left(\left\|_{i} \mathbf{x}_{1}-{ }_{i} \mathbf{x}_{2}\right\|\right) & \cdots & \varphi\left(\left\|_{i} \mathbf{x}_{1}-{ }_{i} \mathbf{x}_{n_{i}}\right\|\right) \\
\varphi\left(\left\|_{i} \mathbf{x}_{2}-{ }_{i} \mathbf{x}_{1}\right\|\right) & \varphi\left(\left\|_{i} \mathbf{x}_{2}-{ }_{i} \mathbf{x}_{2}\right\|\right) & \cdots & \varphi\left(\left\|_{i} \mathbf{x}_{2}-{ }_{i}\mathbf{x}_{n_{i}}\right\|\right) \\
\vdots & \vdots & \vdots & \vdots \\
\varphi\left(\left\|_{i} \mathbf{x}_{n_{i}}-i \mathbf{x}_{1}\right\|\right) & \varphi\left(\left\|_{i} \mathbf{x}_{n_{i}}-{ }_{i} \mathbf{x}_{2}\right\|\right) & \cdots & \varphi\left(\left\|_{i} \mathbf{x}_{n_{i}}-{ }_{i} \mathbf{x}_{n_{i}}\right\|\right)
\end{array}\right]\left[\begin{array}{c}
{ }_{i}\alpha_{1} \\
{ }_{i}\alpha_{2} \\
\vdots \\
{ }_{i}\alpha_{n_{i}}
\end{array}\right] \\
& =:\left({ }_{i} \boldsymbol{\Phi}\right)\left({ }_{i} \boldsymbol \alpha\right)
\end{aligned}
\end{equation}
with 
${ }_{i}\mathbf{f}=[f({ }_{i}\mathbf{x}_{1}), \ldots, f({ }_{i}\mathbf{x}_{n_{i}})]^{\top}$. From \eqref{eq3.1} it follows that  ${ }_{i} \boldsymbol \alpha = ({ }_{i} \boldsymbol{\Phi})^{-1}{ }_{i}\mathbf{f}$ and 
\begin{equation*}
f^*({ }_{i}\mathbf{x})=\left[\varphi\left(\left\|_{i} \mathbf{x}-{ }_{i}\mathbf{x}_{1}\right\|\right), \ldots, \varphi\left(\left\|_{i} \mathbf{x}-{ }_{i} \mathbf{x}_{n_{i}}\right\|\right)\right]\left({ }_{i} \boldsymbol{\Phi}\right)^{-1}{ }_{i}\mathbf{f},
\end{equation*}
whose $l$-order differential approximation reads
\begin{equation}\label{eq3.5}
f^{*(l)}({ }_{i} \mathbf{x})=\left[\varphi^{(l)}\left(\left\|_{i} \mathbf{x}-{ }_{i}\mathbf{x}_{1}\right\|\right), \ldots, \varphi^{(l)}\left(\left\|_{i} \mathbf{x}-{ }_{i}\mathbf{x}_{n_{i}}\right\|\right)\right]\left({ }_{i} \boldsymbol{\Phi}\right)^{-1}{ }_{i}\mathbf{f}.
\end{equation}
Let $n_i=5$ without loss of generality, that is, for each inner node $\mathbf{x}_{i}$, the local influence domain is 
\begin{equation*}
{ }_{i}\Omega=\left\{{ }_{i}\mathbf{x}_{1}, { }_{i}\mathbf{x}_{2}, { }_{i}\mathbf{x}_{3}, { }_{i}\mathbf{x}_{4}, { }_{i}\mathbf{x}_{5}\right\}
\end{equation*}
with $\mathbf{x}_{i}={ }_{i} \mathbf{x}_{3}$ being the center.  
Based on \eqref{eq3.5}, we have the approximation of  $f^{(l)}\left(\mathbf{x}_{i}\right)$ with $l\in\mathbb N_+$ as follows 
\begin{align}
&\quad  f^{*(l)}\left(\mathbf{x}_{i}\right) =f^{*(l)}\left({ }_{i}\mathbf{x}_{3}\right)\nonumber \\
& =[\varphi^{(l)}(\|_{i} \mathbf{x}_{3}-{ }_{i}\mathbf{x}_{1}\|), \varphi^{(l)}(\|_{i} \mathbf{x}_{3}-{ }_{i}\mathbf{x}_{2}\|), \varphi^{(l)}(\|_{i} \mathbf{x}_{3}-{ }_{i}\mathbf{x}_{3}\|),\varphi^{(l)}(\|_{i} \mathbf{x}_{3}-{ }_{i}\mathbf{x}_{4}\|), \varphi^{(l)}(\|_{i} \mathbf{x}_{3}-{ }_{i}\mathbf{x}_{5}\|)]({ }_{i} \boldsymbol{\Phi})^{-1} { }_{i} \mathbf{f} \nonumber\\
& = :\left[_{i} d_{-2}^{(l)}, { }_{i}d_{-1}^{(l)}, { }_{i}d_{0}^{(l)}, { }_{i}d_{1}^{(l)},{ }_{i} d_{2}^{(l)}\right]{ }_{i}\mathbf{f},
\end{align}
which yields
\begin{equation}\begin{aligned}
\mathbf f^{*(l)} &= \left[\begin{array}{c}
f^{*(l)}\left(\mathbf{x}_{0}\right) 
\dots 
f^{*(l)}\left(\mathbf{x}_{i}\right) 
\dots 
f^{*(l)}\left(\mathbf{x}_{L+1}\right)
\end{array}\right]^\top \\
&=\left[\begin{array}{ccccccccccccc}
\cdots & \cdots & \cdots & \cdots & \cdots & \cdots & \cdots & \cdots & \cdots & \cdots & \cdots & \cdots & \cdots \\
\vdots & \vdots & \vdots & \vdots & \vdots & \vdots & \vdots & \vdots & \vdots & \vdots & \vdots & \vdots & \vdots \\
0 & \cdots & { }_{i}d_{-2}^{(l)} & 0 & { }_{i}d_{-1}^{(l)} & 0 & { }_{i}d_{0}^{(l)} & 0 & { }_{i}d_{1}^{(l)} & 0 &{ }_{i}d_{2}^{(l)} & \cdots & 0 \\
\vdots & \vdots & \vdots & \vdots & \vdots & \vdots & \vdots & \vdots & \vdots & \vdots & \vdots & \vdots & \vdots \\
\cdots & \cdots & \cdots & \cdots & \cdots & \cdots & \cdots & \cdots & \cdots & \cdots & \cdots & \cdots & \cdots
\end{array}\right]
\left[\begin{array}{c}
f(\mathbf{x}_{0}) \\
\vdots \\
f(\mathbf{x}_{i}) \\
\vdots \\
f(\mathbf{x}_{L+1})
\end{array}\right] \\
&= : \mathbf D^{(l)} \mathbf f .
\end{aligned}
\label{diffmatrix}
\end{equation}
It can be found that $\mathbf D^{(l)},$ $l\in\mathbb N_+$, is a sparse matrix and there is not any zero between ${ }_{i} d_{k}^{(l)}$ and ${ }_{i} d_{k+1}^{(l)}$ for  $k\in\{-2,-1,0,1\}$, if ${ }_{i} \mathbf{x}_{k}$ and $_{i} \mathbf{x}_{k+1}$ are located next to each other in the full sequence $\left\{_{i}\mathbf{x}_{k}\right\}_{k=1}^{n_{i}}$. 
Especially, under the homogeneous Dirichlet boundary condition, it should be noted that if the node ${ }_{i} \mathbf{x}_{k}$ is out of boundary, the corresponding coefficient ${ }_{i} \alpha_{k}$ is equal to zero. 
Hence, in this case, the form of  $l$-order differential matrix $\mbf D^{(l)}$ for $l\in\mathbb N_+$ becomes
\begin{equation*}
\mbf D^{(l)}=\begin{bmatrix}
{ }_{1}d_{0}^{(l)} & { }_{1}d_{1}^{(l)} & { }_{1}d_{2}^{(l)} & & & & &  &  \\
{ }_{2}d_{-1}^{(l)} & { }_{2}d_{0}^{(l)} & { }_{2}d_{1}^{(l)} & { }_{2}d_{2}^{(l)} & & & & &  \\
{ }_{3}d_{-2}^{(l)} & { }_{3}d_{-1}^{(l)} & { }_{3}d_{0}^{(l)} & { }_{3}d_{1}^{(l)} & { }_{3}d_{2}^{(l)} & & & & \\
&  { }_{4}d_{-2}^{(l)} &  { }_{4}d_{-1}^{(l)} &  { }_{4}d_{0}^{(l)} & { }_{4}d_{1}^{(l)} & { }_{4}d_{2}^{(l)} & & & \\
& & \ddots & \ddots & \ddots & \ddots & \ddots & & \\
& & & { }_{L-3}d_{-2}^{(l)} & { }_{L-3}d_{-1}^{(l)} & { }_{L-3}d_{0}^{(l)} & { }_{L-3}d_{1}^{(l)} & { }_{L-3}d_{2}^{(l)} & \\
& & & & { }_{L-2}d_{-2}^{(l)} & { }_{L-2}d_{-1}^{(l)} & { }_{L-2}d_{0}^{(l)} & { }_{L-2}d_{1}^{(l)} & { }_{L-2}d_{2}^{(l)} \\
& & & & & { }_{L-1}d_{-2}^{(l)} & { }_{L-1}d_{-1}^{(l)} & { }_{L-1}d_{0}^{(l)} & { }_{L-1}d_{1}^{(l)} \\
&  & & & & & { }_{L}d_{-2}^{(l)} & { }_{L}d_{-1}^{(l)} & { }_{L}d_{0}^{(l)}
\end{bmatrix}.
\end{equation*}

Approximating the spatial derivative in \eqref{eq2.1} by $\mathbf D^{(1)}$ of 
the LRBF collocation method leads to a semi-discrete method
\begin{equation}\label{eq3.66}
M d Z_i + K\sum_{k=1}^{n_i} { }_{i}d_k^{(1)} { }_{i}Z_k d t = \nabla S_{1}(Z_i) d t+\nabla S_{2}(Z_i) \circ dW(x_i,t),
\end{equation}
where $i=1,\dots,I-1$, $k=1,\dots,n_i$, $Z_i \approx z(x_i),  {}_{i}Z_k \approx z({ }_{i}x_k)$ and ${}_{i}d_k^{(1)}$ is the element of $\mathbf D^{(1)}.$ 
After making use of the midpoint method  in time, we obtain the meshless LRBF collocation midpoint method of \eqref{eq2.1} as follows
\begin{equation}\label{eq3.6}
M(Z_i^{n+1}- Z_i^{n}) + \Delta tK\sum_{k=1}^{n_i} { }_{i}d_k^{(1)} { }_{i}Z_k^{n+\frac{1}{2}}= \Delta t\nabla S_1(Z_i^{n+\frac{1}{2}}) +\Delta W_i^n \nabla  S_2(Z_i^{n+\frac{1}{2}}),
\end{equation}
where $Z_i^n\approx z(x_i, t_n),  Z_i^{n+\frac{1}{2}} \approx \left(z(x_i, t_n)+ z(x_i, t_{n+1})\right)/2, { }_{i}Z_k^{n+\frac{1}{2}} \approx \left(z( { }_{i}x_k, t_n ) + z( { }_{i}x_k, , t_{n+1})\right)/2$ and  $\Delta W_i^n = W(x_i, t_{n+1}) - W(x_i, t_n) $.

Applying \eqref{eq3.6} to the nonlinear stochastic wave equation with multiplicative noise \eqref{eq2.4}, we derive
\begin{equation}\label{eq3.7}
\left\{\begin{aligned}
&\frac{\mbf U^{n+1} -\mbf U^n}{\Delta t} = \mbf V^{n+\frac{1}{2}}, \\
&\mbf D^{(1)}\mbf U^{n+\frac{1}{2}} = \boldsymbol {\mathcal{W}}^{n+\frac{1}{2}}, \\
&\frac{\mbf V^{n+1} -\mbf V^n}{\Delta t} =  \mbf D^{(1)}\boldsymbol {\mathcal{W}}^{n+\frac{1}{2}} - \mbf F(U^{n+\frac{1}{2}}) + \mbf G(U^{n+\frac{1}{2}})\frac{\Delta \mbf W^n}{\Delta t} , \\
\end{aligned}\right.
\end{equation}
where 
\begin{equation*}
\begin{aligned}
& \mbf U^{n+\frac{1}{2}} = (\mbf U^{n+1} +\mbf U^n)/2,\quad\mbf V^{n+\frac{1}{2}} = (\mbf V^{n+1} +\mbf V^n)/2, \quad \boldsymbol {\mathcal{W}}^{n+\frac{1}{2}} = (\boldsymbol {\mathcal{W}}^{n+1} +\boldsymbol {\mathcal{W}}^n)/2,\\
&\mbf U^n = [U_1^n, \dots, U_{I-1}^n]^\top,\quad~ \mbf V^n = [V_1^n, \dots, V_{I-1}^n]^\top,\quad \boldsymbol {\mathcal{W}}^n = [\mathcal{W}_1^n,  \dots, \mathcal{W}_{I-1}^n]^\top,\\
&\mbf F(U^{n+\frac{1}{2}}) = [f((U_1^{n+1} +U_1^n)/2), \dots, f((U_{I-1}^{n+1} +U_{I-1}^n)/2)]^\top, \\
&\mbf G(U^{n+\frac{1}{2}}) = [g((U_1^{n+1} +U_1^n)/2), \dots, g((U_{I-1}^{n+1} +U_{I-1}^n)/2)]^\top, \\
&\Delta \mbf W^n = [W(x_1,t_{n+1})-W(x_1, t_n), \dots, W(x_{I-1}, t_{n+1})-W(x_{I-1}, t_n)]^\top.
\end{aligned}
\end{equation*}

\begin{thm} \label{thm1}
The  fully-discrete method \eqref{eq3.7} applied to the stochastic wave equation \eqref{eq2.1} with $S_{1}(z)=\frac{1}{2}\left(w^{2}-v^{2}\right)-\tilde f(u)$ and   $ S_{2}(z)=\tilde g(u)$ admits the discrete multi-symplectic conservation law, i.e., 
\begin{equation}\label{eq3.8}
\frac{\omega_i^{n+1} - \omega_i^n}{\Delta t} +\sum_{k=1}^{n_i} { }_{i}d_k^{(1)} { }_{i}\kappa^{n+\frac{1}{2}}_k = 0,\quad 
\end{equation}
where 
\begin{align*}
& \omega_i^n = \frac{1}{2}\mathrm{d} Z_i^n\wedge M\mathrm{d}Z_i^n,  \quad { }_{i}\kappa_k^{n+\frac{1}{2}} = \mathrm{d} Z_i^{n+\frac{1}{2}}
\wedge K\mathrm{d}{ }_{i}Z_k^{n+\frac{1}{2}}, \quad Z_i^n = (U_i^n, V_i^n, \mathcal{W}_i^n)^\top,\\
& _{i}Z_k^{n+\frac{1}{2}} = \left( (_{i}U_k^{n} + _{i}U_k^{n+1} )/2, (_{i}V_k^{n} + _{i}V_k^{n+1} )/2, (_{i}\mathcal{W}_k^{n} + _{i}\mathcal{W}_k^{n+1} )/2 \right)^\top, ~ i =1,\dots, I-1, k =1,\dots, n_i
\end{align*}
with
	\begin{equation*}\begin{aligned}
	&M=\left(\begin{array}{cccc}
	0 & 0 & 1& 0 \\
	0 & 0 & 0& 0 \\
	-1 & 0 & 0& 0\\
	0 & 0 & 0& 0\\
	\end{array}\right), \quad
	K=\left(\begin{array}{cccc}
	0 & 0 & 0& -1 \\
	0 & 0 & 0 & 0\\
	0 & 0 & 0& 0\\
	1 & 0 & 0& 0\\
	\end{array}\right).
	\end{aligned}\end{equation*}
\end{thm}
\noindent{\textbf{Proof.}}
 The system \eqref{eq3.7} can be rewritten in the form of \eqref{eq3.6}, and its discrete variational equation is given by
\begin{equation}\label{eq3.9}
M\frac{\mathrm{d} Z_i^{n+1}- \mathrm{d} Z_i^{n}}{\Delta t} + \sum_{k=1}^{n_i} { }_{i}d_k^{(1)}K \mathrm{d} { }_{i}Z_k^{n+\frac{1}{2}}= \nabla^2 S_1(Z_i^{n+\frac{1}{2}})\mathrm{d} Z_i^{n+\frac{1}{2}}+ \nabla^2S_2(Z_i^{n+\frac{1}{2}})\frac{\Delta W_i^n}{\Delta t}\mathrm{d}  Z_i^{n+\frac{1}{2}}.
\end{equation}
Taking wedge product of the both sides of \eqref{eq3.9} with $\mathrm{d} Z_i^{n+\frac{1}{2}}$ yields 
\begin{equation*}\begin{aligned}
&\frac{\mathrm{d} Z_i^{n+1}+ \mathrm{d} Z_i^{n}}{2} \wedge M\frac{\mathrm{d} Z_i^{n+1}- \mathrm{d} Z_i^{n}}{\Delta t} + \mathrm{d} Z_i^{n+\frac{1}{2}}\wedge  \sum_{k=1}^{n_i} { }_{i}d_k^{(1)}K \mathrm{d} { }_{i}Z_k^{n+\frac{1}{2}}\\
= &\mathrm{d} Z_i^{n+\frac{1}{2}} \wedge \nabla^2 S_1(Z_i^{n+\frac{1}{2}})\mathrm{d} Z_i^{n+\frac{1}{2}}+ \frac{1}{\Delta t} \mathrm{d} Z_i^{n+\frac{1}{2}} \wedge \nabla^2S_2(Z_i^{n+\frac{1}{2}})\Delta W_i^n\mathrm{d} Z_i^{n+\frac{1}{2}}.
\end{aligned}\end{equation*}
By the symmetry of both $\nabla^2 S_1(Z_i^{n+\frac{1}{2}})$ and $\nabla^2 S_2(Z_i^{n+\frac{1}{2}})$, we obtain
 \begin{equation*}
 \frac{1}{\Delta t}\left(\frac{1}{2}\mathrm{d} Z_j^{n+1}\wedge M\mathrm{d}Z_j^{n+1}-\frac{1}{2}\mathrm{d} Z_j^n\wedge M\mathrm{d}Z_j^n\right) +\sum_{k=1}^{n_i} { }_{i}d_k^{(1)} \mathrm{d} Z_i^{n+\frac{1}{2}}\wedge K\mathrm{d}{ }_{i}Z_k^{n+\frac{1}{2}} =0,
\end{equation*}
which is \eqref{eq3.8} by notations $\omega_i^n$ and ${ }_{i}\kappa_k^{n+\frac{1}{2}}$. 
This completes the proof. \hfill$\qed$

\begin{figure}[h]
	\centering
	\subfigure{
		\begin{minipage}{15cm}
			\centering
			\includegraphics[height=4cm,width=4.5cm]{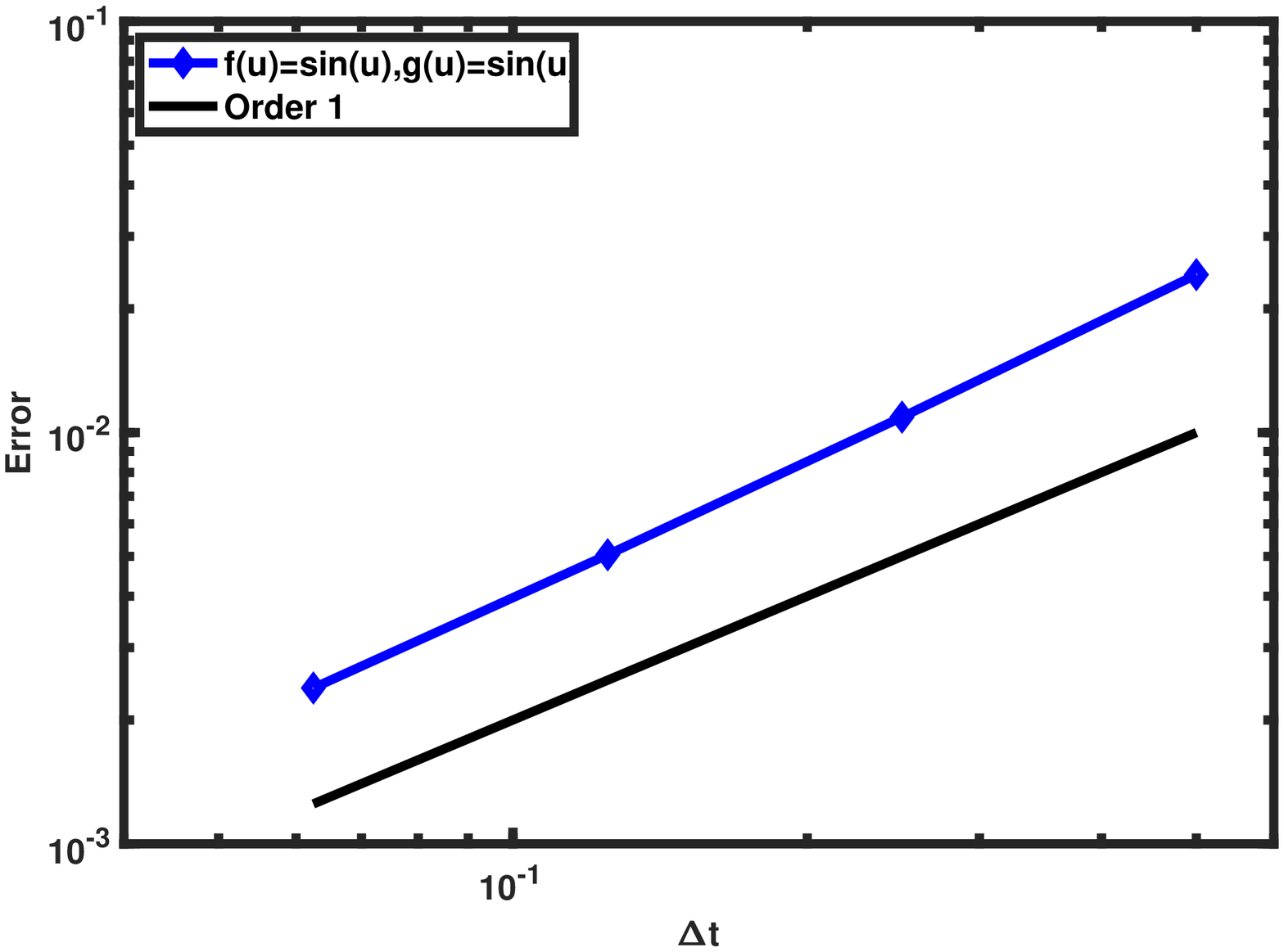}
			\includegraphics[height=4cm,width=4.5cm]{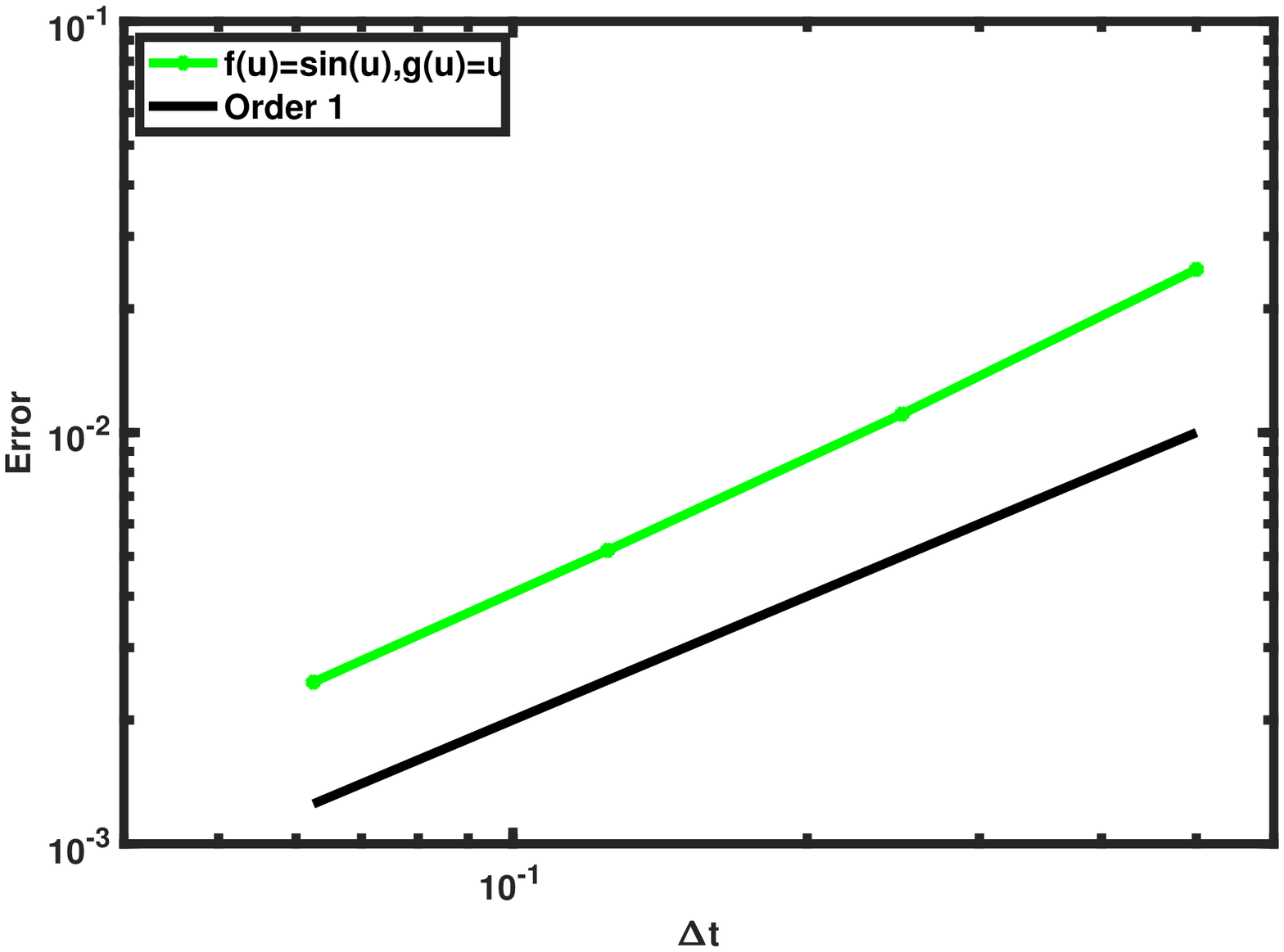}
			\includegraphics[height=4cm,width=4.5cm]{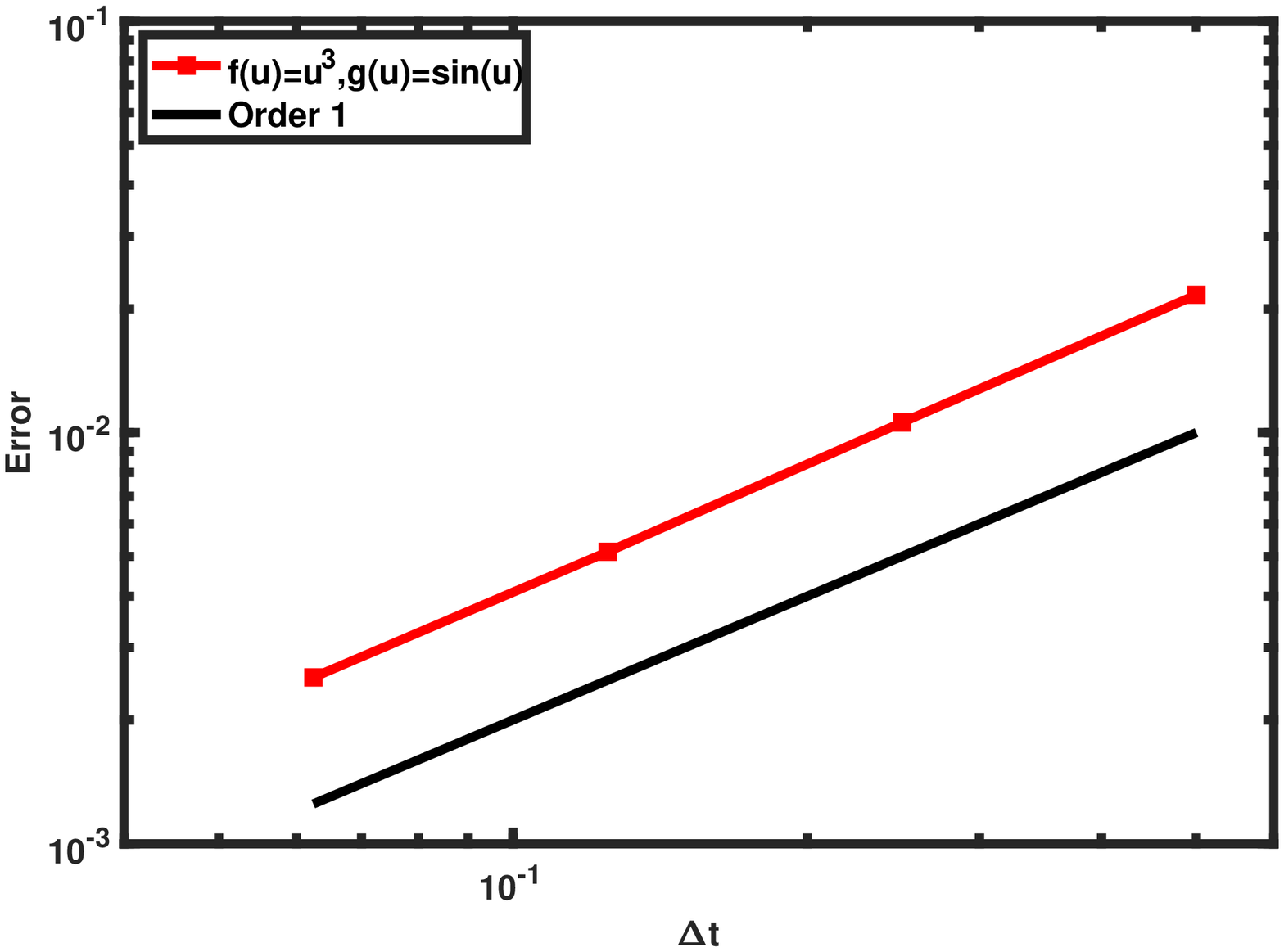}
		\end{minipage}
	}
	\caption{Mean-square convergence order of LRBF collocation midpoint method in temporal direction in the cases of (1) $f(u)= \sin(u), g(u) = \sin(u)$ (2) $f(u)= \sin(u), g(u) = u$ and (3) $f(u)= u^3, g(u) = \sin(u).$} 
	\label{fig1}
\end{figure}

From Theorem \ref{thm1} it is known that the numerical method on the basis of \eqref{eq3.6} for the nonlinear stochastic wave equation possesses the discrete multi-symplectic conservation law. 
Now we perform numerical experiments to illustrate the validity of the proposed method \eqref{eq3.7} for the
1-dimensional stochastic wave equation in different cases: (1)$f(u)=\sin(u), g(u) =\sin(u)$; (2)$f(u)=\sin(u), g(u) =u$; (3)$f(u)=u^3,g(u) =\sin(u).$  
In all the numerical experiments, the expectation is approximated
by taking the average over 1000 realizations. 
Moreover, we choose the orthonormal basis $\left\{e_{k}\right\}_{k\in \mathbb{N}+}$ and the corresponding eigenvalue $\left\{q_{k}\right\}_{k\in \mathbb{N}+}$ of $Q$ as
$e_{k}=\frac{\sqrt{2} }{4}\sin (k \pi x)$ and $q_{k}=\frac{1}{k^6}$, respectively.  And set $x\in(-8,8),$  $u(0) =0,$  $u_t(0) = sech(x),$ and  $u_x(0) = 0.$
The radial basis function is chosen as the inverse multiquadric function $\varphi(x) = 1/\sqrt{1+\|x\|^2}$, i.e., $c=1$. The size of influence domain is taken as $n_i =5$.  Table \ref{tab1} displays strong convergence errors 
against $\Delta t=2^{-s}, s=1, 2, 3, 4,$ on log-log scale  at time $T=1$, which indicates that the meshless LRBF collocation midpoint method offers a good simulation and obtains high precision. 
We regard the numerical approximation obtained by a fine mesh with $\Delta t=2^{-10}, \Delta x = 2^{-5}$ as the exact solution. 
Moreover, compared with the reference line in Fig.  \ref{fig1}, it also can be observed that the mean-square convergence order of the proposed method applied to three cases is 1 in temporal direction.

\begin{table}[htbp]
	\setlength{\abovecaptionskip}{0pt}
	\setlength{\belowcaptionskip}{3pt}
	\centering
       \caption{\label{tab1}Mean-square errors of LRBF collocation midpoint method in time.}
	\begin{tabularx}{0.84\textwidth}{|c|c|c|c|}
		\hline
		& \multicolumn{1}{c|}{$f(u)= \sin(u), g(u) = \sin(u)$}  & \multicolumn{1}{c|}{$f(u)= \sin(u), g(u) = u$}& \multicolumn{1}{c|}{$f(u)= u^3, g(u) = \sin(u)$}\\
		\hline
		$\Delta t$     &$L^2$~error      &$L^2$~error          &$L^2$~error         \\
		\hline
		$2^{-1}$        & 2.4213e-02      &  2.4938e-02         & 2.1630e-02        \\
		$2^{-2}$        & 1.0905e-02      & 1.1091e-02          & 1.0581e-02         \\
		$2^{-3}$        & 5.0449e-03      &  5.1698e-03         & 5.1289e-03         \\
		$2^{-4}$        & 2.3927e-03      &  2.4720e-03         & 2.5391e-03        \\
		\hline
	\end{tabularx}
\end{table}

\begin{rem}
By using \eqref{eq3.6}, we obtain a multi-symplectic method for the stochastic nonlinear Schr\"odinger equation with multiplicative noise \eqref{eq2.6} as follows
\begin{equation}\label{eq3.10}
\left\{\begin{aligned}
&\frac{\mbf P^{n+1} -\mbf P^n}{\Delta t} = - \mbf D^{(1)}\boldsymbol {\mathcal{W}}^{n+\frac{1}{2}}- \left((\mbf P^{n+\frac{1}{2}})^2+(\mbf Q^{n+\frac{1}{2}})^2\right)\mbf Q^{n+\frac{1}{2}} + \mbf Q^{n+\frac{1}{2}}\frac{\Delta \mbf W^n}{\Delta t},\\
&\frac{\mbf Q^{n+1} -\mbf Q^n}{\Delta t} =  \mbf D^{(1)} \mbf V^{n+\frac{1}{2}} + \left((\mbf P^{n+\frac{1}{2}})^2+(\mbf Q^{n+\frac{1}{2}})^2\right)\mbf P^{n+\frac{1}{2}} - \mbf P^{n+\frac{1}{2}}\frac{\Delta \mbf W^n}{\Delta t},\\
&\mbf D^{(1)}\mbf P^{n+\frac{1}{2}} = \mbf V^{n+\frac{1}{2}}, \\
&\mbf D^{(1)}\mbf Q^{n+\frac{1}{2}} = \boldsymbol {\mathcal{W}}^{n+\frac{1}{2}}, 
\end{aligned}\right.
\end{equation}
where $P^{n+\frac 12}_i=\frac{P^{n}_i+P^{n+1}_i}{2},$ $Q^{n+\frac 12}_i=\frac{Q^{n}_i+Q^{n+1}_i}{2},$ $i=1,\dots, I-1,$ and
\begin{equation*}\begin{aligned}
&\left((\mbf P^{n+\frac{1}{2}})^2+(\mbf Q^{n+\frac{1}{2}})^2\right)\mbf Q^{n+\frac{1}{2}} 
= \left[\left((P_1^{n+\frac{1}{2}})^2+(Q_1^{n+\frac{1}{2}})^2\right)Q_1^{n+\frac{1}{2}},  \dots, \left((P_{I-1}^{n+\frac{1}{2}})^2+(Q_{I-1}^{n+\frac{1}{2}})^2\right)Q_{I-1}^{n+\frac{1}{2}}\right]^\top, \\
&\left((\mbf P^{n+\frac{1}{2}})^2+(\mbf Q^{n+\frac{1}{2}})^2\right)\mbf P^{n+\frac{1}{2}} 
= \left[\left((P_1^{n+\frac{1}{2}})^2+(Q_1^{n+\frac{1}{2}})^2\right)P_1^{n+\frac{1}{2}},  \dots, \left((P_{I-1}^{n+\frac{1}{2}})^2+(Q_{I-1}^{n+\frac{1}{2}})^2\right)P_{I-1}^{n+\frac{1}{2}}\right]^\top, \\
&\mbf Q^{n+\frac{1}{2}}\Delta \mbf W^n = [Q_1^{n+\frac{1}{2}}(W(x_1,t_{n+1})-W(x_1, t_n)),  \dots, Q_{I-1}^{n+\frac{1}{2}}(W(x_{I-1},t_{n+1})-W(x_{I-1}, t_n))]^\top,\\
&\mbf P^{n+\frac{1}{2}}\Delta \mbf W^n = [P_1^{n+\frac{1}{2}}(W(x_1,t_{n+1})-W(x_1, t_n)),  \dots, P_{I-1}^{n+\frac{1}{2}}(W(x_{I-1},t_{n+1})-W(x_{I-1}, t_n))]^\top.
\end{aligned}\end{equation*}
Similar to Theorem \ref{thm1}, it can be verified that the fully-discrete method \eqref{eq3.10} possesses the discrete multi-symplectic conservation law 
\begin{equation}
\frac{\omega_i^{n+1} - \omega_i^n}{\Delta t} +\sum_{k=1}^{n_i} { }_{i}d_k^{(1)} { }_{i}\kappa^{n+\frac{1}{2}}_k = 0, 
\end{equation}
where 
\begin{align*}
& \omega_i^n = \frac{1}{2}\mathrm{d} Z_i^n\wedge M\mathrm{d}Z_i^n,  \quad { }_{i}\kappa_k^{n+\frac{1}{2}} = \mathrm{d} Z_i^{n+\frac{1}{2}}\wedge K\mathrm{d}{ }_{i}Z_k^{n+\frac{1}{2}}, \quad Z_i^n = (P_i^n, Q_i^n, V_i^n, \mathcal{W}_i^n)^\top,\\
& _{i}Z_k^{n+\frac{1}{2}} = \left( (_{i}P_k^{n} + _{i}P_k^{n+1} )/2, (_{i}Q_k^{n} + _{i}Q_k^{n+1} )/2, (_{i}V_k^{n} + _{i}V_k^{n+1} )/2, ( _{i}\mathcal{W}_k^{n} + _{i}\mathcal{W}_k^{n+1} )/2 \right)^\top
\end{align*}
with $i =1,\dots, I-1$, $k =1,\dots, n_i,$ and
	\begin{equation*}\begin{aligned}
	&M=\left(\begin{array}{cccc}
	0 & -1 & 0 & 0\\
	1 & 0 & 0 & 0\\
	0 & 0 & 0 & 0\\
	0 & 0 & 0 & 0 \\
	\end{array}\right), \quad
	K=\left(\begin{array}{cccc}
	0 & 0 & 1 & 0\\
	0 & 0 & 0 & 1\\
	-1 & 0 & 0 & 0\\
	0 & -1 & 0 & 0 \\
	\end{array}\right).
	\end{aligned}\end{equation*}
\end{rem}

\begin{rem}
For the stochastic KdV equation \eqref{eq2.10}, making use of \eqref{eq3.6} yields  
\begin{equation}\label{eq3.11}
\left\{\begin{aligned}
&\frac{1}{2}\frac{\mbf U^{n+1} -\mbf U^n}{\Delta t} + \mbf D^{(1)}\mbf V^{n+\frac{1}{2}}=\lambda \frac{\Delta \mbf W^n}{\Delta t},\\
&\frac{1}{2}\frac{ \boldsymbol {\mathcal{P}}^{n+1} - \boldsymbol {\mathcal{P}}^n}{\Delta t} +   \beta\mbf D^{(1)}\boldsymbol {\mathcal{W}}^{n+\frac{1}{2}} = \mbf V^{n+\frac{1}{2}} - \frac{1}{2}(\mbf U^{n+\frac{1}{2}})^2,\\
&\beta \mbf D^{(1)}\mbf U^{n+\frac{1}{2}} = \beta \boldsymbol {\mathcal{W}}^{n+\frac{1}{2}}, \\
&- \mbf D^{(1)} \boldsymbol {\mathcal{P}}^{n+\frac{1}{2}} = - \mbf U^{n+\frac{1}{2}}, 
\end{aligned}\right.
\end{equation}
where $U^{n+\frac 12}_i=\frac{U^{n}_i+U^{n+1}_i}{2}$, $V^{n+\frac 12}_i=\frac{V^{n}_i+V^{n+1}_i}{2}$,  $\mathcal{P}^{n+\frac 12}_i=\frac{\mathcal{P}^{n}_i+\mathcal{P}^{n+1}_i}{2}$,  $\mathcal{W}^{n+\frac 12}_i=\frac{\mathcal{W}^{n}_i+\mathcal{W}^{n+1}_i}{2}$,  ~~$i=1,\dots, I-1$,
\begin{equation*}
(\mbf U^{n+\frac{1}{2}})^2 
= [ (U_1^{n+\frac{1}{2}})^2,  \dots, (U_{I-1}^{n+\frac{1}{2}})^2]^\top,  \Delta \mbf W^n = [W(x_1,t_{n+1})-W(x_1, t_n), \dots, W(x_{I-1}, t_{n+1})-W(x_{I-1}, t_n)]^\top.
\end{equation*}
In fact, the fully-discrete method \eqref{eq3.11} has the discrete multi-symplectic conservation law 
\begin{equation}
\frac{\omega_i^{n+1} - \omega_i^n}{\Delta t} +\sum_{k=1}^{n_i} { }_{i}d_k^{(1)} { }_{i}\kappa^{n+\frac{1}{2}}_k = 0, \qquad i =1,\dots, I-1
\end{equation}
with
\begin{align*}
& \omega_i^n = \frac{1}{2}\mathrm{d} Z_i^n\wedge M\mathrm{d}Z_i^n,  \quad { }_{i}\kappa_k^{n+\frac{1}{2}} = \mathrm{d} Z_i^{n+\frac{1}{2}}\wedge K\mathrm{d}{ }_{i}Z_k^{n+\frac{1}{2}}, \quad Z_i^n = (U_i^n, V_i^n, \mathcal{P}_i^n,\mathcal{W}_i^n)^\top,\\
& _{i}Z_k^{n+\frac{1}{2}} = \left( (_{i}U_k^{n} + _{i}U_k^{n+1} )/2, (_{i}V_k^{n} + _{i}V_k^{n+1} )/2, (_{i}\mathcal{P}_k^{n} + _{i}\mathcal{P}_k^{n+1} )/2, (_{i}\mathcal{W}_k^{n} + _{i}\mathcal{W}_k^{n+1} )/2 \right)^\top, 
\end{align*}
and
	\begin{equation*}\begin{aligned}
	&M=\left(\begin{array}{cccc}
	0 & 0 & -\frac{1}{2} & 0\\
	0 & 0 & 0 & 0\\
	\frac{1}{2} & 0 & 0 & 0\\
	0 & 0 & 0 & 0 \\
	\end{array}\right), \quad
	K=\left(\begin{array}{cccc}
	0 & 0 & 0 & - \beta\\
	0 & 0 &-1 & 0\\
	0 & 1 & 0 & 0\\
	\beta & 0 & 0 & 0 \\
	\end{array}\right).
	\end{aligned}\end{equation*}
\end{rem}

\begin{rem}
For the stochastic Maxwell equation \eqref{Maxwell's equation}, by means  of \eqref{eq3.6}, we obtain  
\begin{equation}\label{eq3.12}
\left\{\begin{aligned}
&\frac{(\mbf E_1)^{n+1} -(\mbf E_1)^n}{\Delta t} = -\mbf D_z^{(1)}(\mbf H_2)^{n+\frac{1}{2}} + \mbf D_y^{(1)}(\mbf H_3)^{n+\frac{1}{2}} 
-\lambda (\mbf H_1)^{n+\frac{1}{2}} \frac{\Delta \mbf W^n}{\Delta t},\\
&\frac{(\mbf E_2)^{n+1} -(\mbf E_2)^n}{\Delta t} = ~\mbf D_z^{(1)}(\mbf H_1)^{n+\frac{1}{2}} - \mbf D_x^{(1)}(\mbf H_3)^{n+\frac{1}{2}} 
- \lambda (\mbf H_2)^{n+\frac{1}{2}} \frac{\Delta \mbf W^n}{\Delta t},\\
&\frac{(\mbf E_3)^{n+1} -(\mbf E_3)^n}{\Delta t} = -\mbf D_y^{(1)}(\mbf H_1)^{n+\frac{1}{2}} + \mbf D_x^{(1)}(\mbf H_2)^{n+\frac{1}{2}} 
- \lambda (\mbf H_3)^{n+\frac{1}{2}} \frac{\Delta \mbf W^n}{\Delta t},\\
&\frac{(\mbf H_1)^{n+1} -(\mbf H_1)^n}{\Delta t} = ~\mbf D_z^{(1)}(\mbf E_2)^{n+\frac{1}{2}} - \mbf D_y^{(1)}(\mbf E_3)^{n+\frac{1}{2}} 
+ \lambda (\mbf E_1)^{n+\frac{1}{2}} \frac{\Delta \mbf W^n}{\Delta t},\\
&\frac{(\mbf H_2)^{n+1} -(\mbf H_2)^n}{\Delta t} = -\mbf D_z^{(1)}(\mbf E_1)^{n+\frac{1}{2}} + \mbf D_x^{(1)}(\mbf E_3)^{n+\frac{1}{2}} 
+ \lambda (\mbf E_2)^{n+\frac{1}{2}} \frac{\Delta \mbf W^n}{\Delta t},\\
&\frac{(\mbf H_3)^{n+1} -(\mbf H_3)^n}{\Delta t} = ~\mbf D_y^{(1)}(\mbf E_1)^{n+\frac{1}{2}} - \mbf D_x^{(1)}(\mbf E_2)^{n+\frac{1}{2}} 
+ \lambda (\mbf E_3)^{n+\frac{1}{2}} \frac{\Delta \mbf W^n}{\Delta t},
\end{aligned}\right.
\end{equation}
where 
\begin{equation*}
\begin{aligned}
& (\mbf E_j)^{n+\frac{1}{2}} = ((\mbf E_j)^{n+1} +(\mbf E_j)^n)/2,\quad(\mbf H_j)^{n+\frac{1}{2}} = ((\mbf H_j)^{n+1} +(\mbf H_j)^n)/2, \\
&(\mbf E_j)^n = [(E_j)_1^n, \dots, (E_j)_{I-1}^n]^\top,\quad~ (\mbf H_j)^n = [(H_j)_1^n, \dots, (H_j)_{I-1}^n]^\top,\quad  j = 1, 2,3,\\
&\Delta \mbf W^n = [W(x_1,t_{n+1})-W(x_1, t_n), \dots, W(x_{I-1}, t_{n+1})-W(x_{I-1}, t_n)]^\top.
\end{aligned}
\end{equation*}
In the  three-dimensional case,  $\mathbf D_{x}^{(1)},$ $\mathbf D_{y}^{(1)}$ and  $\mathbf D_{z}^{(1)}$ are $1$-order differential approximations of partial derivatives $\partial x$, $\partial y$  and  $\partial z$ of LRBF collocation method in \eqref{diffmatrix}, respectively, and
the corresponding elements in above three matrices are denoted by ${ }_{i} d_{x, k}^{(1)}, { }_{i} d_{y, k}^{(1)}, { }_{i} d_{z, k}^{(1)}$ for $i\in\{1,\dots, I-1\}$ and $k\in\{1,\dots, n_i\}$.
The fully-discrete method \eqref{eq3.12} satisfies
\begin{equation}
\frac{\omega_i^{n+1} - \omega_i^n}{\Delta t} 
+\sum_{k=1}^{n_i} { }_{i}d_{x,k}^{(1)} { }_{i}\kappa^{n+\frac{1}{2}}_{1,k}
+\sum_{k=1}^{n_i} { }_{i}d_{y,k}^{(1)} { }_{i}\kappa^{n+\frac{1}{2}}_{2,k}
+\sum_{k=1}^{n_i} { }_{i}d_{z,k}^{(1)} { }_{i}\kappa^{n+\frac{1}{2}}_{3,k}
 = 0, \qquad i =1,\dots, I-1,
\end{equation}
where
\begin{align*}
& \omega_i^n = \frac{1}{2}\mathrm{d} Z_i^n\wedge M\mathrm{d}Z_i^n,  \quad { }_{i}\kappa_{j,k}^{n+\frac{1}{2}} = \mathrm{d} { }Z_i^{n+\frac{1}{2}}\wedge K_j\mathrm{d}{ }_{i}Z_k^{n+\frac{1}{2}}, \quad Z_i^n = ((H_1)_i^n, (H_2)_i^n, (H_3)_i^n, (E_1)_i^n, (E_2)_i^n, (E_3)_i^n )^\top,\\
& _{i}Z_k^{n+\frac{1}{2}} = \left( (_{i}(H_1)_k^{n} + _{i}(H_1)_k^{n+1} )/2, (_{i}(H_2)_k^{n} + _{i}(H_2)_k^{n+1} )/2, (_{i}(H_3)_k^{n} + _{i}(H_3)_k^{n+1} )/2, (_{i}(E_1)_k^{n} + _{i}(E_1)_k^{n+1} )/2,\right.\\
&\left.~\qquad\qquad(_{i}(E_2)_k^{n} + _{i}(E_2)_k^{n+1} )/2, (_{i}(E_3)_k^{n} + _{i}(E_3)_k^{n+1} )/2\right)^\top, 
\end{align*}
and
	$$
	M=\left(\begin{array}{cc}
	0 & -I_{3 \times 3} \\
	I_{3 \times 3} & 0
	\end{array}\right), \quad K_{j}=\left(\begin{array}{cc}
	\mathscr{D}_{j} & 0 \\
	0 & \mathscr{D}_{j}
	\end{array}\right),\quad j=1,2,3
	$$
	with $I_{3 \times 3}$ being a $3 \times 3$ identity matrix,
	$$
	\mathscr{D}_{1}=\left(\begin{array}{ccc}
	0 & 0 & 0 \\
	0 & 0 & -1 \\
	0 & 1 & 0
	\end{array}\right),\quad \mathscr{D}_{2}=\left(\begin{array}{ccc}
	0 & 0 & 1 \\
	0 & 0 & 0 \\
	-1 & 0 & 0
	\end{array}\right),\quad \mathscr{D}_{3}=\left(\begin{array}{ccc}
	0 & -1 & 0 \\
	1 & 0 & 0 \\
	0 & 0 & 0
	\end{array}\right) .
	$$		
\end{rem}

\section{Splitting multi-symplectic Runge--Kutta method}

In this section, we propose the second kind of multi-symplectic methods  for \eqref{eq2.1} via the splitting technique, which avoids the interaction between the nonlinear drift coefficient and the driving process. 
This splitting technique allows us to handle a deterministic Hamiltonian PDE directly, and thus some existing deterministic multi-symplectic method can be exploited.
Motivated by the fact that multi-symplectic Runge--Kutta  methods are a class of efficient derivative-free numerical methods, we concentrate on the splitting multi-symplectic Runge--Kutta method for stochastic Hamiltonian PDEs.

Now we begin our study with the multi-symplectic Runge--Kutta method for deterministic Hamiltonian PDEs
\begin{equation} 
\label{deqhpde}
M d z + K z_{x} d t = \nabla S_{1}(z) d t. 
\end{equation}  
Applying $s$-stage and $r$-stage symplectic Runge--Kutta methods, i.e., $(c, A, b)$ and $(\tilde c, \tilde A, \tilde b)$ as follows
\begin{equation}
\label{butchert1}
\begin{array}{c|ccc}
c_{1} &  a_{11} &\dots &a_{1s}\\
\vdots &   \vdots &&\vdots \\
c_{s} &a_{s1} &\dots &a_{ss}\\
\hline
& b_{1}& \dots&  b_{s}
\end{array},\qquad
\begin{array}{c|ccc}
\tilde c_{1} &  \tilde a_{11} &\dots & \tilde a_{1r}\\
\vdots &   \vdots &&\vdots \\
\tilde c_{r} &\tilde a_{r1} &\dots &\tilde a_{rr}\\
\hline
& \tilde b_{1}& \dots&  \tilde b_{r}
\end{array}, 
\end{equation}
where $s,r\geq 1,$ to \eqref{deqhpde} in space and time, respectively, the
resulting fully-discrete method is as follows:
\begin{equation} 
\label{dMSRK}
\begin{aligned}
&Z_{m}^{k}=z_{i}^{k}+{\Delta x} \sum_{n=1}^{s} a_{m n} \delta_{x}^{n, k} Z_{n}^{k}, \quad \forall~ i=0,1, \ldots, s,\\
&z_{i+1}^{k}=z_{i}^{k}+{\Delta x} \sum_{m=1}^{s} b_{m} \delta_{x}^{m, k} Z_{m}^{k}, \quad \forall~ i=0,1, \ldots, s, \\
&Z_{m}^{k}=z_{m}^{p}+\Delta t \sum_{j=1}^{r} \tilde{a}_{k j} \delta_{t}^{m, j} Z_{m}^{j}, \quad \forall~ p=0,1, \ldots, r,\\
&z_{m}^{p+1}=z_{m}^{p}+\Delta t \sum_{k=1}^{r} \tilde{b}_{k} \delta_{t}^{m, k} Z_{m}^{k}, \quad \forall~ p=0,1, \ldots, r, \\
&M \delta_{t}^{m, k} Z_{m}^{k}+K \delta_{x}^{m, k} Z_{m}^{k}=\nabla_{z} S_{1}\left(Z_{m}^{k}\right),
\end{aligned}
\end{equation}
where $\delta_{t}^{m, k}$ and $\delta_{x}^{m, k}$ are discretizations of two partial derivatives $\partial_{t}$ and $\partial_{x}$, respectively, and
\begin{align}
\label{symplecticcondition}
b_{m} b_{n}-b_{m} a_{m n}-b_{n} a_{n m}=0 \text { and }  \tilde{b}_{k} \tilde{b}_{j}-\tilde{b}_{k} \tilde{a}_{k j}-\tilde{b}_{j} \tilde{a}_{j k}=0 
\end{align}
for all  $m, n=1, \ldots, s$ and $k, j=1, \ldots, r$.
It can be verified that the above stochastic numerical method admits the discrete multi-symplectic conservation law
\begin{align*}
\frac{\omega^{p+1}-\omega^{p}}{\Delta t}+\frac{\kappa_{i+1}-\kappa_{i}}{h}=0, 
\end{align*}
where
$\omega^{p}=\frac{1}{2} \sum\limits_{m=1}^{s} b_{m} \mathrm{d} z_{m}^{p} \wedge M \mathrm{d} z_{m}^{p},$ and $\kappa_{i}=\frac{1}{2} \sum\limits_{k=1}^{r} \tilde b_{k} \mathrm{d} z_{i}^{k} \wedge K \mathrm{d} z_{i}^{k}$ for $p=0,1, \ldots,r$ and $i=0,1,\ldots, s.$  
Applying the  splitting technique to \eqref{eq2.1} in temporal direction, and then we obtain a deterministic Hamiltonian PDE with random input and a stochastic system on $t \in [t_m, t_{m+1}]$ as follows 
\begin{equation} 
\left\{
\begin{aligned}
&M d\bar z + K \bar z_{x} d t = \nabla S_{1}(\bar z) d t,\\
& \bar {z}(t_m) = z(t_{m}),
\end{aligned}
\right. \qquad{\rm and}\qquad
\left\{
\begin{aligned}
& Kz_{x}=0,\\
& M d z = \nabla S_{2}(z) \circ dW(t),\\
& z(t_m) = \overline{z}(t_{m+1}).
\end{aligned}
\right.
\end{equation}
By choosing symplectic methods for the stochastic system and combining \eqref{dMSRK}, 
we obtain the splitting multi-symplectic Runge--Kutta method satisfying the discrete multi-symplectic conservation law. 
Now we construct the splitting multi-symplectic Runge--Kutta method for the nonlinear stochastic wave equation,  stochastic nonlinear Schr\"odinger equation, stochastic KdV equation and stochastic Maxwell equation, one after the other.

We first focus on the nonlinear stochastic wave equation  \eqref{eq2.4} and propose the associated splitting multi-symplectic Runge--Kutta method. 
In detail, we decompose  \eqref{eq2.4} on $[t_0,t_1]$ into a deterministic Hamiltonian PDE with random input
\begin{equation}
\label{eq6.0}
\left\{
\begin{aligned}
&\overline{u}_t = \overline{v}, \\
&\overline{u}_x = \overline{w},\\
&\overline{v}_t  - \overline{w}_x= -f(\overline{u}),\\
&\overline{u} (t_0) = u(t_{0}), ~\overline{v}(t_0) = v(t_{0}),
\end{aligned} \right.
\end{equation}
and a stochastic system
\begin{equation}
\label{eq6.01}
\left\{
\begin{aligned}
& u_x = 0,~w_x =0,\\
& u_t =0,\\
& dv = g(u)  \circ dW(t),\\
&u(t_0)=\overline{u}(t_{1}), ~v(t_0)=\overline{v}(t_{1}).
\end{aligned} \right.
\end{equation}
By making use of $s$-stage and $r$-stage symplectic Runge--Kutta methods  \eqref{butchert1} with $s,r\geq 1$ to approximate \eqref{eq6.0},  together with the application of the symplectic Euler method to the stochastic system \eqref{eq6.01}, we obtain the following fully-discrete method 
\begin{subequations}
\begin{align}
\label{eq6.1}
&U_{i}^{m} = u_{0}^m+ {\Delta x}\sum_{j=1}^{s}a_{ij}\mathcal{W}_{j}^{m},\quad \mathcal{W}_{i}^{m} = w_{0}^{m} + {\Delta x}\sum_{j=1}^{s}a_{ij}\delta_x\mathcal{W}_{j}^{m},\\
\label{eq6.3}
&\overline{u}_{1}^{m}=u_{0}^{m} + {\Delta x}\sum_{i=1}^{s}b_{i}\mathcal{W}_{i}^{m},\quad 
\overline{w}_{1}^{m} = w_{0}^{m} + {\Delta x}\sum_{i=1}^{s}b_{i}\delta_x\mathcal{W}_{i}^{m},\\
\label{eq6.5}
& U_{i}^{m} = u_{i}^{0} +  \Delta t \sum_{n=1}^{r}\tilde{a}_{nm}V_{i}^{n},\quad  
V_{i}^{m} = v_{i}^{0} +  \Delta t \sum_{n=1}^{r}\tilde{a}_{nm}\left(\delta_x\mathcal{W}_{i}^{n} -f(U_{i}^{n})\right),\\
\label{eq6.7}
& \overline{u}_{i}^{1} = u_{i}^{0} +  \Delta t \sum_{m=1}^{r}\tilde{b}_{m}V_{i}^{m},\quad 
\overline{v}_{i}^{1} = v_{i}^{0} +  \Delta t \sum_{m=1}^{r}\tilde{b}_{m}\left(\delta_x\mathcal{W}_{i}^{m} -f(U_{i}^{m})\right),\\
\label{eq6.9}
&u_{1}^{m} = \overline{u}_{1}^{m},\quad w_{1}^{m} = \overline{w}_{1}^{m}, \\
\label{eq6.10}  
&u_{i}^{1} =  \overline{u}_{i}^{1},\quad v_{i}^{1} = \overline{v}_{i}^{1} + g(\overline{u}_{i}^{1})\Delta W_{i}^{1}, 
\end{align}
\end{subequations}
where $i =1,\dots, s$,  $m=1,\dots, r$, $\delta_x$ is the discretization of the partial
derivative $\partial_x$, and  $U_{i}^{m}\approx u(c_i\Delta x, \tilde c_m\Delta t )$, $u_{i}^{0}\approx u(c_i\Delta x,0)$,$u_{i}^{1}\approx u(c_i\Delta x,\Delta t)$, $\overline{u}_{i}^{1}\approx \overline{u}(c_i\Delta x,\Delta t)$, $u_{0}^{m}\approx u(0,\tilde c_m\Delta t)$, $u_{1}^{m}\approx u(\Delta x,\tilde c_m\Delta t)$, $\overline{u}_{1}^{m}\approx \overline{u}(\Delta x,\tilde c_m\Delta t)$,  etc., with $c_i =\sum_{j=1}^{s}a_{ij}$, $\tilde c_m =\sum_{n=1}^{r}\tilde{a}_{mn}$. 

\begin{thm} 
\label{thm4}
Assume that the symplectic condition \eqref{symplecticcondition} or equivalently,
\begin{align*}
B A+A^{\top} B-b b^{\top}=0,\qquad \tilde{B} \tilde{A}+\tilde{A}^{\top} \tilde{B}-\tilde{b} \tilde{b}^{\top}=0,
\end{align*}
where $B=\operatorname{diag}(b)$ and $\tilde{B}=\operatorname{diag}(\tilde{b})$, holds.
Then the fully-discrete method \eqref{eq6.1}-\eqref{eq6.10} admits the discrete multi-symplectic conservation law 
\begin{equation*}
\sum_{i=1}^{s}\frac{b_{i}}{\Delta t}\left(\mathrm{d}u_{i}^{1} \wedge \mathrm{d}v_{i}^{1} - \mathrm{d}u_{i}^{0} \wedge \mathrm{d}v_{i}^{0}\right)
- \sum_{m=1}^{r} \frac{\tilde{b}_m}{{\Delta x}}\left(\mathrm{d} u_{1}^{m} \wedge \mathrm{d}w_{1}^{m}-\mathrm{d} u_{0}^{m} \wedge \mathrm{d} w_{0}^{m}\right)=0,
\end{equation*}
where $s,r\in\mathbb N_+.$
\end{thm}
\noindent{\textbf{Proof.}}  
By utilizing  \eqref{eq6.9}-\eqref{eq6.10}, we obtain
\begin{align*}
&\sum_{i=1}^{s}\frac{b_{i}}{\Delta t}\left(\mathrm{d}u_{i}^{1} \wedge \mathrm{d}v_{i}^{1} - \mathrm{d}u_{i}^{0} \wedge \mathrm{d}v_{i}^{0}\right)
- \sum_{m=1}^{r} \frac{\tilde{b}_m}{{\Delta x}}\left(\mathrm{d} u_{1}^{m} \wedge \mathrm{d}w_{1}^{m}-\mathrm{d} u_{0}^{m} \wedge \mathrm{d} w_{0}^{m}\right)\\
=&\frac{1}{\Delta t}\sum_{i=1}^{s} b_{i}\left(\mathrm{d}\overline{u}_{i}^{1} \wedge \mathrm{d}\overline{v}_{i}^{1} - \mathrm{d}u_{i}^{0} \wedge \mathrm{d}v_{i}^{0}\right)
- \frac{1}{{\Delta x}}\sum_{m=1}^{r} \tilde{b}_m\left(\mathrm{d} \overline{u}_{1}^{m} \wedge \mathrm{d} \overline{w}_{1}^{m}-\mathrm{d} u_{0}^{m} \wedge \mathrm{d} w_{0}^{m}\right).
\end{align*}
For fixed $i \in\{1,\dots, s\}$ and $m\in\{1,\dots,r\}$,  taking advantage of \eqref{eq6.7} leads to 
\begin{equation}
\begin{aligned}\label{eq6.17}
\mathrm{d}\overline{u}_{i}^{1} \wedge \mathrm{d}\overline{v}_{i}^{1} = &\mathrm{d}u_{i}^{0} \wedge \mathrm{d}v_{i}^{0} + \Delta t \sum_{m=1}^{r}\tilde{b}_{m}\mathrm{d}u_{i}^{0} \wedge \mathrm{d}\left(\delta_x\mathcal{W}_{i}^{m} -f(U_{i}^{m})\right)\\
& +\Delta t \sum_{m=1}^{r}\tilde{b}_{m} \mathrm{d}V_{i}^{m} \wedge \mathrm{d}v_{i}^{0}
+\Delta t^2\sum_{m, l=1}^{r}\tilde{b}_{m}\tilde{b}_{l} \mathrm{d}V_{i}^{m} \wedge \mathrm{d}\left(\delta_x\mathcal{W}_{i}^{l} -f(U_{i}^{l})\right).
\end{aligned}\end{equation}
Applying 
$\mathrm{d}U_{i}^{m} = \mathrm{d} u_{i}^{0} +  \Delta t \sum\limits_{n=1}^{r}\tilde{a}_{nm}\mathrm{d}V_{i}^{n}$ and $\mathrm{d}V_{i}^{m} = \mathrm{d}v_{i}^{0} +  \Delta t \sum\limits_{n=1}^{r}\tilde{a}_{nm}\mathrm{d}\left(\delta_xW_{i}^{n} -f(U_{i}^{n})\right)$
to \eqref{eq6.17}, we get
\begin{equation}\begin{aligned}
&\mathrm{d}\overline{u}_{i}^{1} \wedge \mathrm{d}\overline{v}_{i}^{1}\\ 
= &\mathrm{d}u_{i}^{0} \wedge \mathrm{d}v_{i}^{0} + \Delta t \sum_{l=1}^{r}\tilde{b}_{l}\mathrm{d} U_{i}^{l} \wedge \mathrm{d}\left(\delta_x\mathcal{W}_{i}^{l} -f(U_{i}^{l})\right)- \Delta t^2 \sum_{m,l=1}^{r}\tilde{b}_{l}\tilde{a}_{ml}\mathrm{d}V_{i}^{m} \wedge \mathrm{d}\left(\delta_x\mathcal{W}_{i}^{l} -f(U_{i}^{l})\right)\\
& -\Delta t^2 \sum_{m,l=1}^{r}\tilde{b}_{m}\tilde{a}_{lm} \mathrm{d}V_{i}^{m} \wedge \mathrm{d}\left(\delta_x\mathcal{W}_{i}^{l} -f(U_{i}^{l})\right)+\Delta t^2\sum_{m,l=1}^{s}\tilde{b}_{m}\tilde{b}_{l} \mathrm{d}V_{i}^{m} \wedge \mathrm{d}\left(\delta_x\mathcal{W}_{i}^{l} -f(U_{i}^{l})\right).
\end{aligned}\end{equation}
Based on \eqref{symplecticcondition}, we obtain
\begin{equation}\label{eq6.18}
\mathrm{d}\overline{u}_{i}^{1} \wedge \mathrm{d}\overline{v}_{i}^{1} = \mathrm{d}u_{i}^{0} \wedge \mathrm{d}v_{i}^{0} +  \Delta t \sum_{l=1}^{r}\tilde{b}_{l}\mathrm{d} U_{i}^{l} \wedge \mathrm{d}\left(\delta_x\mathcal{W}_{i}^{l}\right).
\end{equation}
Similarly, from \eqref{eq6.3} it follows that
\begin{equation*}
\begin{aligned}
\mathrm{d}\overline{u}_{1}^{m} \wedge \mathrm{d}\overline{w}_{1}^{m}  = &\mathrm{d}u_{0}^{m} \wedge \mathrm{d}w_{0}^{m}
+ {\Delta x} \sum_{i=1}^{s}b_{i}\mathrm{d} u_{0}^{m} \wedge \mathrm{d} \left(\delta_x\mathcal{W}_{i}^{m}\right)
\\
&+{\Delta x} \sum_{i=1}^{s}b_{i} \mathrm{d}\mathcal{W}_{i}^{m} \wedge \mathrm{d}w_{0}^{m}
+{\Delta x}^2\sum_{i, k=1}^{s}b_{i}b_{k} \mathrm{d}\mathcal{W}_{i}^{m} \wedge \mathrm{d}\left(\delta_x\mathcal{W}_{k}^{m}\right).
\end{aligned}
\end{equation*}
By means of \eqref{eq6.1}, we derive
\begin{align*}
&\mathrm{d}U_{i}^{m} = \mathrm{d}u_{0}^{m}  + {\Delta x}\sum_{j=1}^{s}a_{i j}\mathrm{d}\mathcal{W}_{j}^{m},\quad
\mathrm{d}\mathcal{W}_{i}^{m} = \mathrm{d}w_{0}^{m}  + {\Delta x}\sum_{j=1}^{s}a_{ij}\mathrm{d}(\delta_x\mathcal{W}_{j}^{m}),
\end{align*}
which yields
\begin{align}
\mathrm{d}\overline{u}_{1}^{m} \wedge \mathrm{d}\overline{w}_{1}^{m}
= &\mathrm{d}u_{0}^{m} \wedge \mathrm{d}w_{0}^{m}
+ {\Delta x} \sum_{i=1}^{s}b_{i}\mathrm{d}U_{i}^{m}\wedge \mathrm{d} \left(\delta_x\mathcal{W}_{i}^{m}\right)- {\Delta x}^2 \sum_{i,k=1}^{s}b_{k}a_{ki}\mathrm{d}\mathcal{W}_{i}^{m}\wedge \mathrm{d} \left(\delta_x\mathcal{W}_{k}^{m}\right)\nonumber\\
& -{\Delta x}^2 \sum_{i,k=1}^{s}b_{i}a_{ik} \mathrm{d}\mathcal{W}_{i}^{m} \wedge \mathrm{d}(\delta_x\mathcal{W}_{k}^{m})+{\Delta x}^2\sum_{i,k=1}^{s}b_{i}b_{k} \mathrm{d}\mathcal{W}_{i}^{m} \wedge \mathrm{d}\left(\delta_x\mathcal{W}_{k}^{m}\right)\nonumber\\
= &\mathrm{d}u_{0}^{m} \wedge \mathrm{d}w_{0}^{m}
+ {\Delta x} \sum_{i=1}^{s}b_{i}\mathrm{d}U_{i}^{m}\wedge \mathrm{d} \left(\delta_x\mathcal{W}_{i}^{m}\right).\label{eq6.20}
\end{align}
Combining \eqref{eq6.18} and \eqref{eq6.20}, we deduce
\begin{equation*}
\begin{aligned}\label{eq6.21}
&\frac{1}{\Delta t}\sum_{i=1}^{s} b_{i}\left(\mathrm{d}\overline{u}_{i}^{1} \wedge \mathrm{d}\overline{v}_{i}^{1} - \mathrm{d}u_{i}^{0} \wedge \mathrm{d}v_{i}^{0}\right)
- \frac{1}{{\Delta x}}\sum_{m=1}^{r} \tilde{b}_m\left(\mathrm{d} \overline{u}_{1}^{m} \wedge \mathrm{d} \overline{w}_{1}^{m}-\mathrm{d} u_{0}^{m} \wedge \mathrm{d} w_{0}^{m}\right)\\
=&\sum_{i=1}^{s} \sum_{l=1}^{r} b_{i}\tilde{b}_{l}\mathrm{d} U_{i}^{l} \wedge \mathrm{d}\left(\delta_x\mathcal{W}_{i}^{l}\right)
-\sum_{m=1}^{r} \sum_{i=1}^{s}\tilde{b}_{m}b_i\mathrm{d}U_{i}^{m}\wedge \mathrm{d} \left(\delta_x\mathcal{W}_{i}^{m}\right)=0,
\end{aligned}
\end{equation*}
which completes the proof.\hfill$\qed$

\begin{figure}[h]
	\centering
	\subfigure{
		\begin{minipage}{15cm}
\centering
\includegraphics[height=4cm,width=4.5cm]{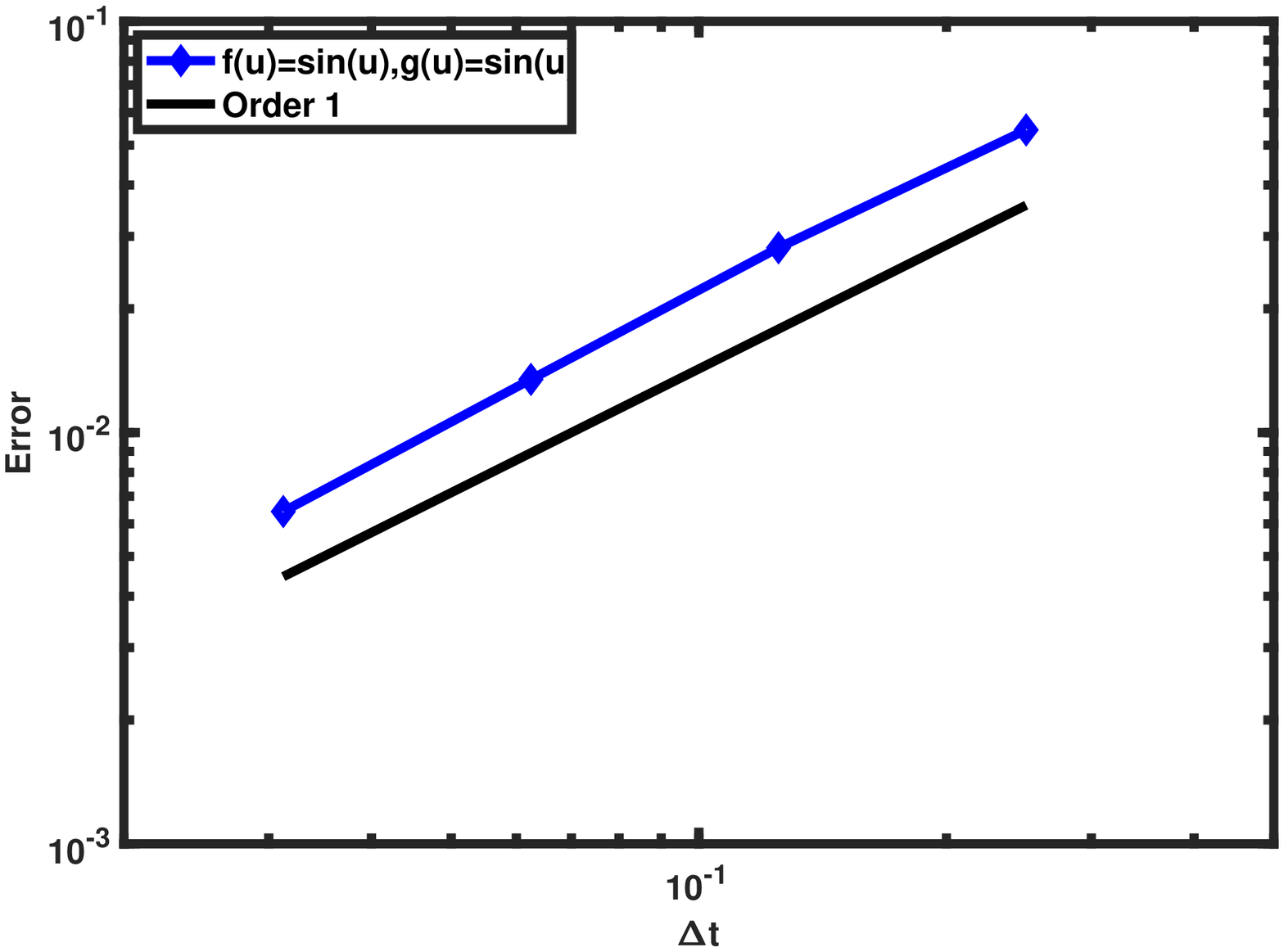}
\includegraphics[height=4cm,width=4.5cm]{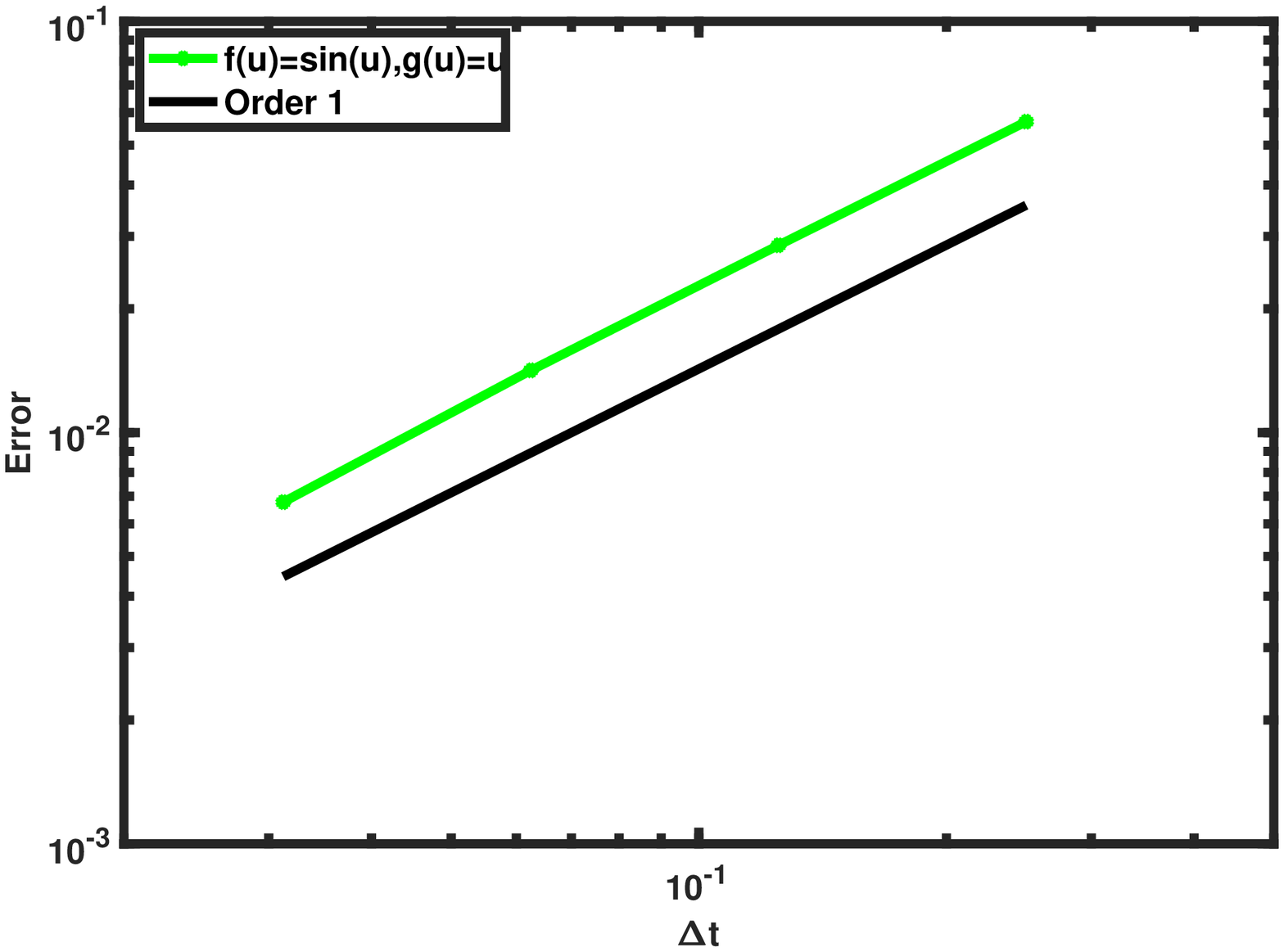}
\includegraphics[height=4cm,width=4.5cm]{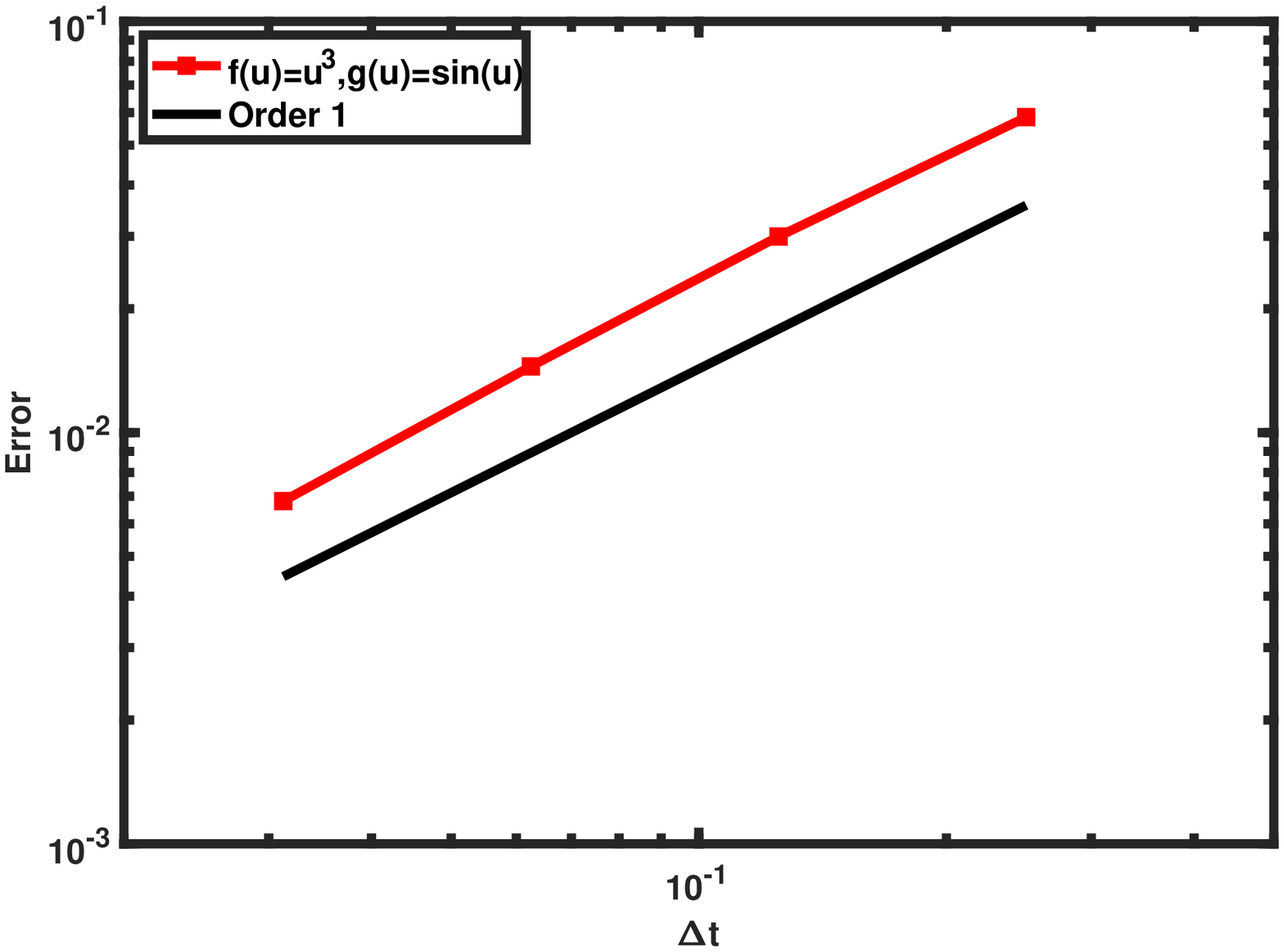}
\end{minipage}
	}
	\caption{Mean-square convergence order of \eqref{smspk11} in temporal direction in the cases of (1) $f(u)= \sin(u), g(u) = \sin(u)$ (2) $f(u)= \sin(u), g(u) = u$ and (3) $f(u)= u^3, g(u) = \sin(u).$} 
	\label{fig2}
\end{figure}

\begin{ex}
	\label{ex2.1}
	If $s=r=1,$ based on symplectic Runge--Kutta methods
	\begin{align*}
	\begin{array}{c|c}
	\frac 12 &  \frac 12\\
	\hline
	&1
	\end{array},\quad 
	\begin{array}{c|c}
	\frac 12 &  \frac 12\\
	\hline
	&1
	\end{array},
	\end{align*}
	we get a numerical method for the nonlinear stochastic wave equation as follows
	\begin{align}
	\label{smspk11}
	&U_{1}^1 = u_{0}^1+ {\Delta x}\frac 12\mathcal{W}_{1}^1,\quad \mathcal{W}_{1}^1 = w_{0}^1 + {\Delta x}\frac 12\delta_x\mathcal{W}_{1}^1,\nonumber\\
	&\overline{u}_{1}^1=u_{0}^1 + {\Delta x}\mathcal{W}_{1}^1,\quad 
	\overline{w}_{1}^1 = w_{0}^1 + {\Delta x}\delta_x\mathcal{W}_{1}^1,\nonumber\\
	& U_{1}^1 = u_{1}^{0} +  \Delta t \frac 12V_{1}^{1},\quad  
	V_{1}^1 = v_{1}^{0} +  \Delta t \frac 12\left(\delta_x\mathcal{W}_{1}^{1} -f(U_{1}^{1})\right),\\
	& \overline{u}_{1}^{1} = u_{1}^{0} +  \Delta t V_{1}^{1},\quad 
	\overline{v}_{1}^{1} = v_{1}^{0} +  \Delta t \left(\delta_x\mathcal{W}_{1}^1 -f(U_{1}^1)\right),\nonumber\\
	&u_{1}^1 = \overline{u}_{1}^1,\quad w_{1}^1 = \overline{w}_{1}^1,\quad v_{1}^{1} = \overline{v}_{1}^{1} + g(\overline{u}_{1}^{1})\Delta W_{1}^{1}. \nonumber
	\end{align}
Similar to the numerical experiments in Section \ref{Sec;MLRBF}, we apply the above multi-symplectic method to approximating the
1-dimensional stochastic wave equation in three cases: {\rm (1)}$f(u)=\sin(u), g(u) =\sin(u)$; {\rm (2)}$f(u)=\sin(u), g(u) =u$; {\rm (3)}$f(u)=u^3,g(u) =\sin(u)$. 
Here, we take $ x \in (-\pi, \pi),$  set $u(0) =0, u_t(0,x) = \sin(x), u_x(0) = 0,$ and let the orthonormal basis $\left\{e_{k}\right\}_{k\in \mathbb{N}+}$ and the corresponding eigenvalue $\left\{q_{k}\right\}_{k\in \mathbb{N}+}$ of $Q$ be
$e_{k}=\frac{1}{\sqrt{\pi}}\sin (k x)$ and $q_{k}=\frac{1}{k^6},$ respectively. 
Table \ref{tab2} shows the mean-square error against $\Delta t=2^{-s}, s=2, 3, 4, 5$ on log-log scale at time $T = 1.$  
We regard the numerical approximation obtained by a fine mesh with $\Delta t=2^{-8}, \Delta x = 2^{-7}\pi$ as the exact solution. 
It can be found from Fig.  \ref{fig2} that the proposed numerical method has accuracy of mean-square order $1$ in temporal direction. 	
\end{ex}

\begin{table}[htbp]
	\setlength{\abovecaptionskip}{0pt}
	\setlength{\belowcaptionskip}{3pt}
	\centering
       \caption{\label{tab2}Mean-square errors of \eqref{smspk11} in time.}
	\begin{tabularx}{0.84\textwidth}{|c|c|c|c|}
		\hline
		& \multicolumn{1}{c|}{$f(u)= \sin(u), g(u) = \sin(u)$}  & \multicolumn{1}{c|}{$f(u)= \sin(u), g(u) = u$}& \multicolumn{1}{c|}{$f(u)= u^3, g(u) = \sin(u)$}\\
		\hline
		$\Delta t$    &$L^2$~error               &$L^2$~error                &$L^2$~error     \\
		\hline
		$2^{-2}$        & 5.4462e-02              & 5.7036e-02               & 5.8488e-02      \\
		$2^{-3}$        & 2.8150e-02              & 2.8546e-02               & 2.9977e-02         \\
		$2^{-4}$        & 1.3469e-02              & 1.4167e-02               & 1.4489e-02          \\
		$2^{-5}$        & 6.4268e-03              & 6.7747e-03               & 6.8146e-03         \\
		\hline
	\end{tabularx}
\end{table}

\begin{figure}[h]
	\centering
	\subfigure{
		\begin{minipage}{12cm}
			\centering
			\includegraphics[height=4cm,width=4.5cm]{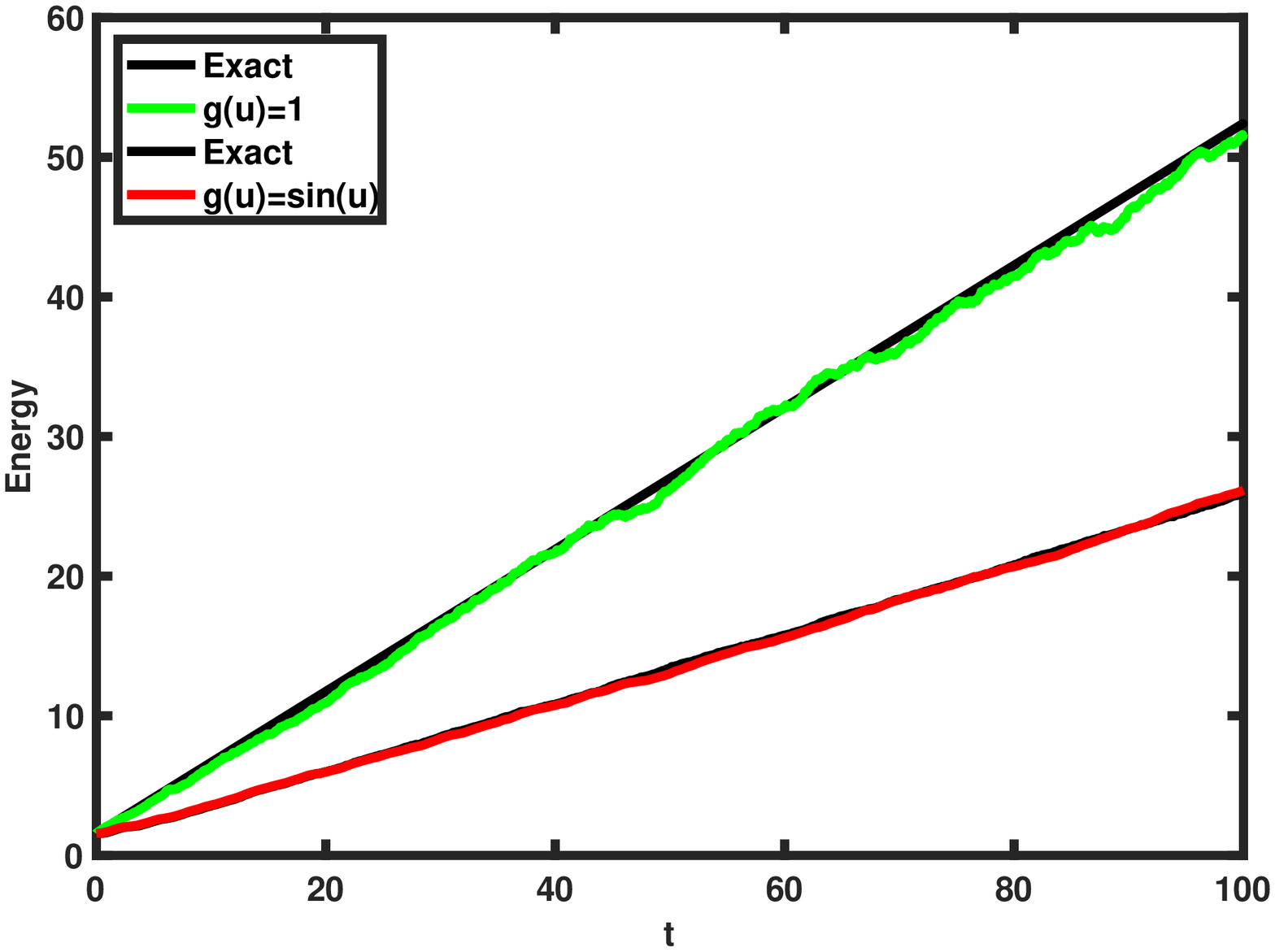}
			\includegraphics[height=4cm,width=4.5cm]{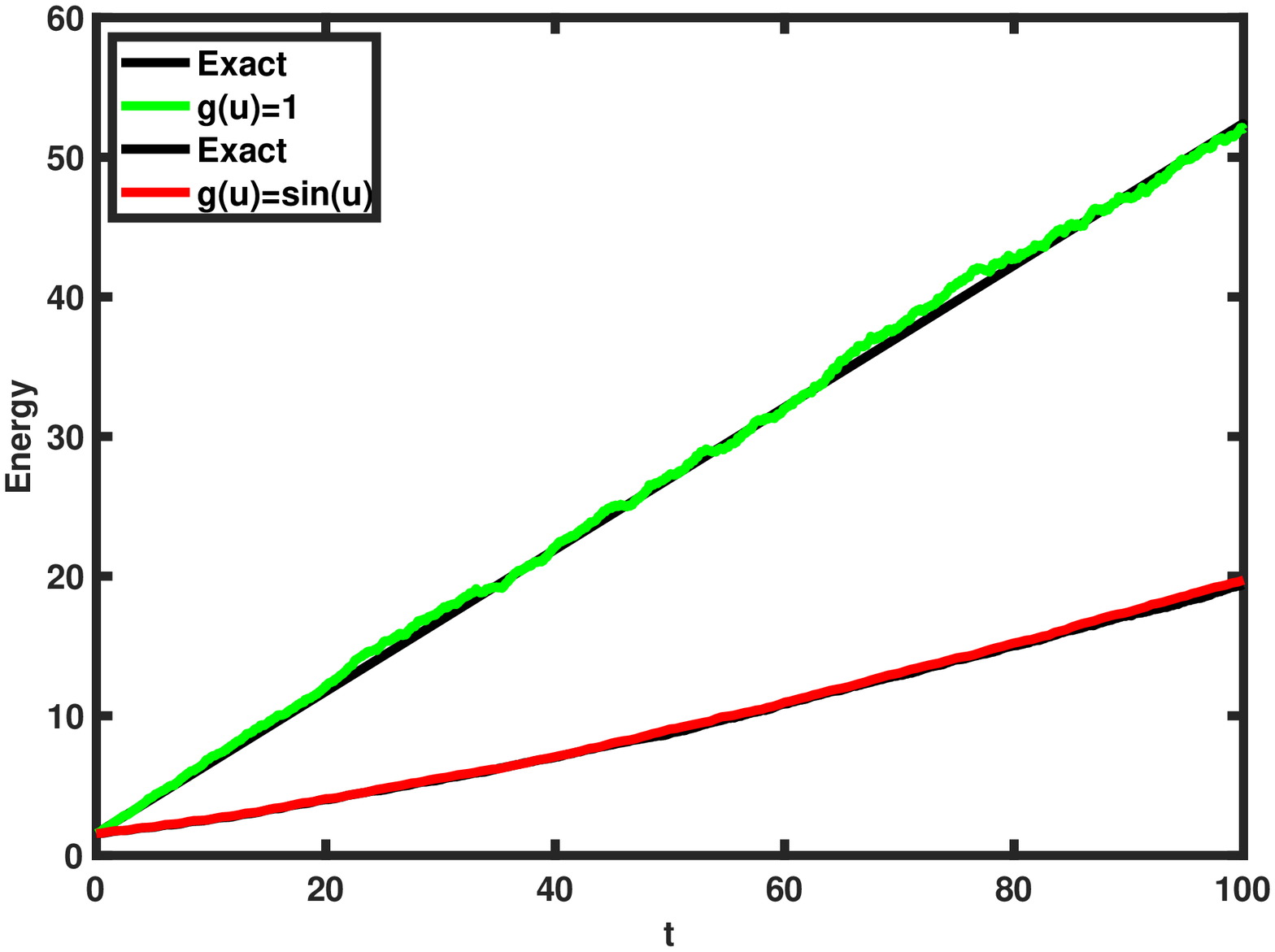}
		\end{minipage}
	}
	\caption{ Averaged energy evolution of \eqref{smspk11} (left: $f(u) = 0$, right: $f(u) = u$) with $\Delta t= 1/20, h=\pi/20.$}
	\label{fig3}
\end{figure}


If $\tilde f(u)$ is at most quadratic, then the fully-discrete method \eqref{eq6.1}-\eqref{eq6.10}  under the symplectic condition \eqref{symplecticcondition} preserves the discrete averaged energy evolution law. 
This property is illustrated by Fig.  \ref{fig3}, from which plots the quantity $\frac{\Delta x}{2}\mathbb E[\sum_{i=1}^{s}(v_i^n)^2+ (w_i^n)^2 + 2\tilde f(u_i^n)],$ $n =1,\ldots,N,$ for \eqref{smspk11}  in two cases: (1) $f(u) = 0 $ and (2) $f(u) = u$, respectively.  
The reference line (black line) in Fig.  \ref{fig3} stands for the averaged energy evolution law of the exact solution. 
It can be observed that \eqref{smspk11} preserves perfectly the averaged energy evolution law.
In detail, when $g(u)=1,$ \eqref{smspk11} reproduces the linear growth of the averaged energy, which coincides with the theoretical results.

Now we turn to the stochastic nonlinear Schr\"odinger equation \eqref{eq2.6}.
Repeating the similar procedures as in the case of the stochastic wave equation, we first split the stochastic nonlinear Schr\"odinger equation into on $[t_0,t_1]$ a deterministic system with random input
\begin{equation}\label{eq7.0}
 \left\{
\begin{aligned}
&\overline{q}_t -\overline{v}_x = \left(\overline{p}^{2}+\overline{q}^{2}\right)\overline{p} , \\
&\overline{p}_t +\overline{w}_x = -\left(\overline{p}^{2}+\overline{q}^{2}\right)\overline{q} , \\
&\overline{p}_x  = \overline{v},~\overline{q}_x  = \overline{w},\\
&\overline{p} (t_0)= p(t_{0}), ~\overline{q} (t_0)= q(t_{0}),
\end{aligned} \right.
\end{equation}
and a stochastic system 
\begin{equation}\label{eq7.01}
\left\{
\begin{aligned}
& p_x =0,~q_x =0,\\
& v_x =0,~ w_x =0,\\
& d q = -p \circ dW(t),\\
& d p = q \circ dW(t),\\
&p(t_0)=\overline{p}(t_{1}), ~q(t_{0})=\overline{q}(t_{1}),
\end{aligned}\right.
\end{equation}
Using $s$-stage and $r$-stage  Runge--Kutta methods \eqref{butchert1} 
with $s, r\geq 1$ to discretize \eqref{eq7.0}, together with the symplectic Euler method applied to \eqref{eq7.01}, yields the fully-discrete method 
\begin{align}
\label{eq7.1}
&P_{i}^{m} = p_{0}^{m}  + {\Delta x}\sum_{j=1}^{s}a_{ij}V_{j}^{m},\quad
Q_{i}^{m} = q_{0}^{m}  + {\Delta x}\sum_{j=1}^{s}a_{ij}\mathcal{W}_{j}^{m},\nonumber\\
&\mathcal{W}_{i}^{m} = w_{0}^{m}  +{\Delta x}\sum_{j=1}^{s}a_{ij}\delta_x\mathcal{W}_{j}^{m},\quad
V_{i}^{m} = v_{0}^{m}+ {\Delta x}\sum_{j=1}^{s}a_{ij}\delta_xV_{j}^{m},\nonumber\\
&\overline{p}_{1}^{m}= p_{0}^{m} +{\Delta x}\sum_{i=1}^{s}b_{i}V_{i}^{m},\quad
\overline{q}_{1}^{m} = q_{0}^{m} + {\Delta x}\sum_{i=1}^{s}b_{i}\mathcal{W}_{i}^{m},\nonumber\\
&\overline{w}_{1}^{m} = w_{0}^{m} +{\Delta x}\sum_{i=1}^{s}b_{i}\delta_x\mathcal{W}_{i}^{m},\quad
\overline{v}_{1}^{m} = v_{0}^{m} + {\Delta x}\sum_{i=1}^{s}b_{i}\delta_xV_{i}^{m},\nonumber\\
 &Q_{i}^{m} = q_{i}^{0} +  \Delta t \sum_{n=1}^{r}\tilde{a}_{nm}\left(\delta_xV_{i}^{n} + (({P_{i}^{n}})^2 + ({Q_{i}^{n}})^2)P_{i}^{n}\right),\\
 &P_{i}^{m} = p_{i}^{0} +  \Delta t \sum_{n=1}^{r}\tilde{a}_{nm}\left(-\delta_x\mathcal{W}_{i}^{n}-(({P_{i}^{n}})^2 + ({Q_{i}^{n}})^2)Q_{i}^{n}\right),\nonumber\\
&\overline{q}_{i}^{1} = q_{i}^{0} +  \Delta t \sum_{m=1}^{r}\tilde{b}_{m}\left(\delta_xV_{i}^{m} + ((P_{i}^{m})^2 + (Q_{i}^{m})^2)P_{i}^{m}\right),\nonumber\\
&\overline{p}_{i}^{1} = p_{i}^{0} +  \Delta t \sum_{m=1}^{r}\tilde{b}_{m}\left(-\delta_x\mathcal{W}_{i}^{m} - ((P_{i}^{m})^2 + (Q_{i}^{m})^2)Q_{i}^{m}\right),\nonumber\\ 
&p_{1}^{m} = \overline{p}_{1}^{m}  ,~q_{1}^{m} = \overline{q}_{1}^{m}  ,~v_{1}^{m} = \overline{v}_{1}^{m}  ,~w_{1}^{m} = \overline{w}_{1}^{m}  ,\nonumber\\
&q_{i}^{1} =  \overline{q}_{i}^{1} - \overline{p}_{i}^{1}\Delta W_{i}^{1}, ~p_{i}^{1} =  \overline{p}_{i}^{1} + q_{i}^{1}\Delta W_{i}^{1},\nonumber 
\end{align}	
where  $P_{i}^{m}\approx p(c_i\Delta x, \tilde c_m\Delta t )$, $p_{i}^{0}\approx p(c_i\Delta x,0)$, $p_{i}^{1}\approx p(c_i\Delta x,\Delta t)$, $\overline{p}_{i}^{1}\approx \overline{p}(c_i\Delta x,\Delta t)$,  $p_{0}^{m}\approx p(0,\tilde c_m\Delta t)$, $p_{1}^{m}\approx p(\Delta x,\tilde c_m\Delta t)$, $\overline{p}_{1}^{m}\approx \overline{p}(\Delta x,\tilde c_m\Delta t)$, etc., with $c_i =\sum_{j=1}^{s}a_{ij}$, $\tilde c_m =\sum_{n=1}^{r}\tilde{a}_{mn}$, $i =1,\dots, s$,  $m=1,\dots, r$.  Similar to Theorem \ref{thm4}, we obtain that the fully-discrete method \eqref{eq7.1} preserves the discrete multi-symplectic conservation law.

\begin{thm} \label{thm5}
Under the symplectic condition \eqref{symplecticcondition},  
the fully-discrete method \eqref{eq7.1} preserves the discrete multi-symplectic conservation law 
\begin{equation*}\begin{aligned}
&\sum_{i=1}^{s}\frac{b_{i}}{\Delta t}\left(\mathrm{d}q_{i}^{1} \wedge \mathrm{d}p_{i}^{1} - \mathrm{d}q_{i}^{0} \wedge \mathrm{d}p_{i}^{0}\right)+ \sum_{m=1}^{r}\frac{\tilde{b}_m}{\Delta x} \left(\mathrm{d} p_{1}^{m} \wedge \mathrm{d}v_{1}^{m}-\mathrm{d} p_{0}^{m} \wedge \mathrm{d} v_{0}^{m} + \mathrm{d} q_{1}^{m} \wedge \mathrm{d}w_{1}^{m}-\mathrm{d} q_{0}^{m} \wedge \mathrm{d} w_{0}^{m}\right)=0.
\end{aligned}\end{equation*}
\end{thm}

Analogously, in the case of the  stochastic KdV equation with additive noise \eqref{eq2.10}, we first decompose it on $t \in [t_0, t_{1}]$
into a deterministic system with random input
\begin{equation}\label{eq9.0}
 \left\{
\begin{aligned}
&\overline{u}_t +2\overline{v}_x= 0 , \\
&\overline{\rho}_t +2\beta\overline{w}_x= 2\overline{v} -\overline{u}^{2}, \\
&\overline{u}_x  = \overline{w},\\
&\overline{\rho}_x  = \overline{u},\\
&\overline{u}(t_{0})= u(t_{0}), ~\overline{\rho}(t_{0})=\rho(t_{0}),
\end{aligned} \right.
\end{equation}
and a  stochastic system 
\begin{equation}\label{eq9.01}
 \left\{
\begin{aligned}
& v_x =0,~w_x =0,\\
& u_x =0,~ \rho_x =0,\\
& d u=  2\lambda \circ dW(t),\\
& \rho_t =  0,\\
&u(t_{0})=\overline{u}(t_{1}), ~\rho(t_{0})=\overline{\rho}(t_{1}).
\end{aligned} \right.
\end{equation}
Next, we take advantage of  $s, r$-stage symplectic Runge--Kutta methods, where $s, r\geq1,$  to numerically solve the deterministic Hamiltonian PDE \eqref{eq9.0} and use symplectic Euler method to  approximate \eqref{eq9.01}, respectively. The resulting numerical method on $t \in [t_0, t_{1}]$ is as follows
 \begin{align}
 \label{eq9.1}
 &V_{i}^{m} = v_{0}^{m}  + {\Delta x}\sum_{j=1}^{s}a_{ij}\delta_xV_{j}^{m},\quad \mathcal{W}_{i}^{m} = w_{0}^{m}  + {\Delta x}\sum_{j=1}^{s}a_{ij}\delta_x\mathcal{W}_{j}^{m},\nonumber\\
 &U_{i}^{m} = u_{0}^{m}  + {\Delta x}\sum_{j=1}^{s}a_{ij}\mathcal{W}_{j}^{m},\quad \mathcal{P}_{i}^{m} = \rho_{0}^{m}  + {\Delta x}\sum_{j=1}^{s}a_{ij}U_{j}^{m},\nonumber\\
 &\overline{v}_{1}^{m} = v_{0}^{m} + {\Delta x}\sum_{i=1}^{s}b_{i}\delta_xV_{i}^{m},\quad
 \overline{w}_{1}^{m} = w_{0}^{m} + {\Delta x}\sum_{i=1}^{s}b_{i}\delta_x\mathcal{W}_{i}^{m},\nonumber\\
 &\overline{u}_{1}^{m} = u_{0}^{m} + {\Delta x}\sum_{i=1}^{s}b_{i}\mathcal{W}_{i}^{m},\quad
 \overline{\rho}_{1}^{m} = \rho_{0}^{m} + {\Delta x}\sum_{i=1}^{s}b_{i}U_{i}^{m},\\
 &U_{i}^{m} = u_{i}^{0} +  \Delta t \sum_{n=1}^{r}\tilde{a}_{nm}\left(-2\delta_xV_{i}^{n} \right),\quad\mathcal{P}_{i}^{m} = \rho_{i}^{0} +  \Delta t \sum_{n=1}^{r}\tilde{a}_{nm}\left(-2\beta \delta_x\mathcal{W}_{i}^{n} + 2V_{i}^{n} - ({U_{i}^{n}})^2 \right),\nonumber\\
&\overline{u}_{i}^{1} = u_{i}^{0} +  \Delta t \sum_{m=1}^{r}\tilde{b}_{m}\left(-2\delta_xV_{i}^{m} \right),\quad 
\overline{ \rho}_{i}^{1} =  \rho_{i}^{0} +  \Delta t \sum_{m=1}^{r}\tilde{b}_{m}\left(-2\beta \delta_x\mathcal{W}_{i}^{m} + 2V_{i}^{m} - (U_{i}^{m})^2\right),\nonumber\\
&v_{1}^{m} = \overline{v}_{1}^{m}  ,~w_{1}^{m} = \overline{w}_{1}^{m},~u_{1}^{m} = \overline{u}_{1}^{m},~ \rho_{1}^{m} = \overline{\rho}_{1}^{m}, ~\rho_{i}^{1} = \overline{\rho}_{i}^{1}  ,~u_{i}^{1} =  \overline{u}_{i}^{1} + 2\lambda \Delta W_{i}^{1}, \nonumber
 \end{align}
 where $i =1,\dots, s$,  $m=1,\dots, r$,  and  $\mathcal{P}_{i}^{m}\approx \rho(c_i\Delta x, \tilde c_m\Delta t )$, $\rho_{i}^{0}\approx \rho(c_i\Delta x,0)$, $\rho_{i}^{1}\approx \rho(c_i\Delta x,\Delta t)$, $\overline{\rho}_{i}^{1}\approx \overline{\rho}(c_i\Delta x,\Delta t)$, $\rho_{0}^{m}\approx \rho(0,\tilde c_m\Delta t)$, $\rho_{1}^{m}\approx \rho(\Delta x,\tilde c_m\Delta t)$, $\overline{\rho}_{1}^{m}\approx \overline{\rho}(\Delta x,\tilde c_m\Delta t)$,  etc., with $c_i =\sum_{j=1}^{s}a_{ij}$, $\tilde c_m =\sum_{n=1}^{r}\tilde{a}_{mn}$. Similar to the proof of Theorem \ref{thm4}, we have the following theorem. 
 
\begin{thm} \label{thm6}
Assume that the symplectic condition \eqref{symplecticcondition} holds. 
Then the fully-discrete method \eqref{eq9.1} preserves the discrete multi-symplectic conservation law 
\begin{equation*}
\sum_{i=1}^{s}\frac{b_{i}}{\Delta t}\left(\mathrm{d}\rho_{i}^{1} \wedge \mathrm{d}u_{i}^{1} - \mathrm{d}\rho_{i}^{0} \wedge \mathrm{d}u_{i}^{0}\right)+ \sum_{m=1}^{r} \frac{2\tilde{b}_m}{{\Delta x}}\left( \mathrm{d} \rho_{1}^{m} \wedge \mathrm{d}v_{1}^{m}-\mathrm{d} \rho_{0}^{m} \wedge \mathrm{d} v_{0}^{m} + \beta \mathrm{d} w_{1}^{m} \wedge \mathrm{d}u_{1}^{m}- \beta\mathrm{d} w_{0}^{m} \wedge \mathrm{d} u_{0}^{m}\right)=0.
\end{equation*}
\end{thm}

Similarly, for the  stochastic Maxwell equation with multiplicative noise \eqref{Maxwell's equation}, we decompose it  on $t \in [t_0, t_{1}]$ into a deterministic PDE with random initial value
\begin{equation}\label{eq10.0}
 \left\{
\begin{aligned}
&(\overline{E}_1)_t + (\overline{H}_2)_z - (\overline{H}_3)_y= 0 , 
 (\overline{E}_2)_t + (\overline{H}_3)_x - (\overline{H}_1)_z= 0, \\
&(\overline{E}_3)_t + (\overline{H}_1)_y - (\overline{H}_2)_x= 0 , 
 (\overline{H}_1)_t + (\overline{E}_3)_y - (\overline{E}_2)_z= 0, \\
&(\overline{H}_2)_t + (\overline{E}_1)_z - (\overline{E}_3)_x= 0 , 
 (\overline{H}_3)_t + (\overline{E}_2)_x - (\overline{E}_1)_y= 0, \\
&\overline{E}_i(t_{0})= E_i(t_{0}), ~\overline{H}_i(t_{0})= H_i(t_{0}),~ i = 1,2,3,
\end{aligned} \right.
\end{equation}
and a  stochastic system
\begin{equation}
\label{eq10.1}
 \left\{
\begin{aligned}
& \mathscr{D}_{1} \mbf H_x=0,~\mathscr{D}_{2} \mbf H_y=0,\mathscr{D}_{3} \mbf H_z=0,\\
& \mathscr{D}_{1} \mbf E_x=0,~\mathscr{D}_{2} \mbf E_y=0,\mathscr{D}_{3} \mbf E_z=0,\\
&  \mbf H_t = \lambda \mbf E \circ dW(t),  \mbf E_t = -\lambda \mbf H \circ dW(t), \\
&\mbf H(t_{0})=\overline{\mbf H}(t_{1}), ~\mbf E(t_{0})=\overline{\mbf E}(t_{1}).
\end{aligned} \right.
\end{equation}
By exploiting $s$-stage and $r$-stage symplectic Runge--Kutta methods to discretize \eqref{eq10.0} and symplectic Euler method to discretize \eqref{eq10.1}, we obtain the numerical method on $t \in [t_0, t_{1}]$ as follows
 \begin{align}
 \label{eq10.1}
 &\mathscr{D}_{1}(\mathbf H)^m_{kln} = \mathscr{D}_{1}(\mathbf H)^m_{0ln} + \Delta x  \sum_{j=1}^{s} a_{kj}^{(1)}\mathscr{D}_{1}(\delta_x \mathbf H)^m_{jln}, 
 \quad \mathscr{D}_{1}(\mathbf E)^m_{kln} = \mathscr{D}_{1}(\mathbf E)^m_{0ln} + \Delta x  \sum_{j=1}^{s} a_{kj}^{(1)}\mathscr{D}_{1}(\delta_x \mathbf E)^m_{jln}, \nonumber\\
 & \mathscr{D}_{2}(\mathbf H)^m_{kln} = \mathscr{D}_{2}(\mathbf H)^m_{k0n} + \Delta y  \sum_{j=1}^{s} a_{lj}^{(2)}\mathscr{D}_{2}(\delta_y \mathbf H)^m_{kjn}, 
\quad \mathscr{D}_{2}(\mathbf E)^m_{kln} = \mathscr{D}_{2}(\mathbf E)^m_{k0n} + \Delta y  \sum_{j=1}^{s} a_{lj}^{(2)}\mathscr{D}_{2}(\delta_y \mathbf E)^m_{kjn}, \nonumber\\
  &\mathscr{D}_{3}(\mathbf H)^m_{kln} = \mathscr{D}_{3}(\mathbf H)^m_{kl0} + \Delta z  \sum_{j=1}^{s} a_{nj}^{(3)}\mathscr{D}_{3}(\delta_z \mathbf H)^m_{klj}, 
\quad \mathscr{D}_{3}(\mathbf E)^m_{kln} = \mathscr{D}_{3}(\mathbf E)^m_{kl0} + \Delta z  \sum_{j=1}^{s} a_{nj}^{(3)}\mathscr{D}_{3}(\delta_z \mathbf E)^m_{klj}, \nonumber\\
& (\mathbf H)^m_{kln} = (\mathbf H)^0_{kln} + \Delta t  \sum_{i=1}^{r} \tilde{a}_{mi}\left( -\mathscr{D}_{1}(\delta_x \mathbf E)^i_{kln}  -\mathscr{D}_{2}(\delta_y \mathbf E)^i_{kln}  -\mathscr{D}_{3}(\delta_z \mathbf E)^i_{kln} \right),   \nonumber\\
& (\mathbf E)^m_{kln} = (\mathbf E)^0_{kln} + \Delta t  \sum_{i=1}^{r} \tilde{a}_{mi}\left( \mathscr{D}_{1}(\delta_x \mathbf H)^i_{kln}  +\mathscr{D}_{2}(\delta_y \mathbf H)^i_{kln}  + \mathscr{D}_{3}(\delta_z \mathbf H)^i_{kln} \right),   \nonumber\\
 &\mathscr{D}_{1}(\overline{\mathbf  H})^m_{1ln} = \mathscr{D}_{1}(\mathbf \overline{H})^m_{0ln} + \Delta x  \sum_{k=1}^{s} b_{k}^{(1)}\mathscr{D}_{1}(\delta_x \mathbf H)^m_{kln}, 
\quad \mathscr{D}_{1}(\overline{\mathbf  E})^m_{1ln} = \mathscr{D}_{1}(\mathbf E)^m_{0ln} + \Delta x  \sum_{k=1}^{s} b_{k}^{(1)}\mathscr{D}_{1}(\delta_x \mathbf E)^m_{kln}, \nonumber\\
& \mathscr{D}_{2}(\overline{\mathbf  H})^m_{k1n} = \mathscr{D}_{2}(\mathbf \overline{H})^m_{k0n} + \Delta y  \sum_{l=1}^{s} b_{l}^{(2)}\mathscr{D}_{2}(\delta_y \mathbf H)^m_{kln}, 
\quad\mathscr{D}_{2}(\overline{\mathbf  E})^m_{k1n} = \mathscr{D}_{2}(\mathbf E)^m_{k0n} + \Delta y  \sum_{l=1}^{s} b_{l}^{(2)}\mathscr{D}_{2}(\delta_y \mathbf E)^m_{kln}, \\
&\mathscr{D}_{3}(\overline{\mathbf  H})^m_{kl1} = \mathscr{D}_{3}(\mathbf \overline{H})^m_{kl0} + \Delta z  \sum_{n=1}^{s} b_{n}^{(3)}\mathscr{D}_{3}(\delta_z \mathbf H)^m_{kln},
\quad \mathscr{D}_{3}(\overline{\mathbf  E})^m_{kl1} = \mathscr{D}_{3}(\mathbf E)^m_{kl0} + \Delta z  \sum_{n=1}^{s} b_{n}^{(3)}\mathscr{D}_{3}(\delta_z \mathbf E)^m_{kln}, \nonumber\\
& (\overline{\mathbf  H})^1_{kln} = (\mathbf H)^0_{kln} + \Delta t  \sum_{m=1}^{r} \tilde{b}_{m}\left( -\mathscr{D}_{1}(\delta_x \mathbf E)^m_{kln}  -\mathscr{D}_{2}(\delta_y \mathbf E)^m_{kln}  -\mathscr{D}_{3}(\delta_z \mathbf E)^m_{kln} \right) ,\nonumber\\
& (\overline{\mathbf  E})^1_{kln} = (\mathbf E)^0_{kln} + \Delta t  \sum_{m=1}^{r} \tilde{b}_{m}\left( \mathscr{D}_{1}(\delta_x \mathbf H)^m_{kln}  +\mathscr{D}_{2}(\delta_y \mathbf H)^m_{kln}  + \mathscr{D}_{3}(\delta_z \mathbf H)^m_{kln} \right) ,\nonumber\\
& \mathscr{D}_{1} (\mathbf  H)^m_{1ln} =  \mathscr{D}_{1} (\overline{\mathbf  H})^m_{1ln}, \quad \mathscr{D}_{1} (\mathbf  E)^m_{1ln} =  \mathscr{D}_{1} (\overline{\mathbf  E})^m_{1ln}, \quad  \mathscr{D}_{2} (\mathbf  H)^m_{k1n} =  \mathscr{D}_{2} (\overline{\mathbf  H})^m_{k1n}, \quad  \mathscr{D}_{2} (\mathbf  E)^m_{k1n} =  \mathscr{D}_{2} (\overline{\mathbf  E})^m_{k1n}, \nonumber\\
&\mathscr{D}_{3} (\mathbf  H)^m_{kl1} =  \mathscr{D}_{3} (\overline{\mathbf  H})^m_{kl1}, \quad  \mathscr{D}_{3} (\mathbf  E)^m_{kl1} =  \mathscr{D}_{3} (\overline{\mathbf  E})^m_{kl1},\nonumber\\
& (\mathbf  E)^1_{kln} = (\overline{\mathbf  E})^1_{kln} - \lambda(\overline{\mathbf  H})^1_{kln} \Delta W_{kln}^{1},\quad (\mathbf  H)^1_{kln} = (\overline{\mathbf  H})^1_{kln} + \lambda(\mathbf  E)^1_{kln} \Delta W_{kln}^{1},\nonumber
\end{align}
where $\delta_x, \delta_y, \delta_z $ are discretizations of partial derivatives $\partial_x, \partial_y, \partial_z$, $ {\mathbf  E}^{0}_{kln} \approx \mathbf  E(c^{(1)}_k\Delta x, c^{(2)}_l\Delta y, c^{(3)}_n\Delta z, 0)$, $ {\mathbf  E}^{m}_{kln} \approx \mathbf  E(c^{(1)}_k\Delta x, c^{(2)}_l\Delta y, c^{(3)}_n\Delta z, \tilde c_m\Delta t )$,  $ {\mathbf  E}^{1}_{kln} \approx \mathbf  E(c^{(1)}_k\Delta x, c^{(2)}_l\Delta y, c^{(3)}_n\Delta z,\Delta t )$,  $ {\overline{\mathbf  E} }^{1}_{kln} \approx \overline{\mathbf  E} (c^{(1)}_k\Delta x, c^{(2)}_l\Delta y, $\\$c^{(3)}_n\Delta z, \Delta t )$, $ {\mathbf  E}^{m}_{0ln} \approx \mathbf  E(0, c^{(2)}_l\Delta y, c^{(3)}_n\Delta z, \tilde c_m\Delta t )$, $ {\mathbf  E}^{m}_{1ln} \approx \mathbf  E(\Delta x, c^{(2)}_l\Delta y, c^{(3)}_n\Delta z, \tilde c_m\Delta t )$, $\overline{\mathbf  E}^{m}_{1ln} \approx \overline{\mathbf  E}(\Delta x, c^{(2)}_l\Delta y,$\\$ c^{(3)}_n\Delta z, \tilde c_m\Delta t )$, $ {\mathbf  E}^{m}_{k0n} \approx \mathbf  E(c^{(1)}_k\Delta x, 0, c^{(3)}_n\Delta z, \tilde c_m\Delta t )$, $ {\mathbf  E}^{m}_{k1n} \approx \mathbf  E(c^{(1)}_k\Delta x, \Delta y, c^{(3)}_n\Delta z, \tilde c_m\Delta t )$,  $ \overline{\mathbf  E}^{m}_{k1n} \approx \overline{\mathbf  E}(c^{(1)}_k\Delta x, $\\$\Delta y, c^{(3)}_n\Delta z, \tilde c_m\Delta t )$,  $ {\mathbf  E}^{m}_{kl0} \approx \mathbf  E(c^{(1)}_k\Delta x, c^{(2)}_l\Delta y, 0, \tilde c_m\Delta t )$,  $ {\mathbf  E}^{m}_{kl1} \approx \mathbf  E(c^{(1)}_k\Delta x, c^{(2)}_l\Delta y, \Delta z, \tilde c_m\Delta t )$,  $ \overline{\mathbf  E}^{m}_{kl1} \approx \overline{\mathbf  E}(c^{(1)}_k\Delta x,$\\$ c^{(2)}_l\Delta y, \Delta z, \tilde c_m\Delta t )$, etc., with $c^{(1)}_k =\sum_{j=1}^{s}a_{k j}^{(1)}$, $c^{(2)}_l =\sum_{j=1}^{s}a_{lj}^{(2)}$, $c^{(3)}_n =\sum_{j=1}^{s}a_{nj}^{(3)}$, $\tilde c_m =\sum_{n=1}^{r}\tilde{a}_{mn}$, $ 1\leq j, k,l, n \leq s, ~ 1\leq i,m \leq r.$ Moreover, the noise increment $\Delta W_{kln}^{1}:=W\left(t_{1}, x_{k},y_{l}, z_{n}\right)-W\left(t_{0}, x_{k},y_{l}, z_{n}\right)$. 
Similar to the proof of Theorem \ref{thm4}, we have the following theorem.

\begin{thm} \label{thm7}
Suppose that the symplectic condition \eqref{symplecticcondition} holds. 
Then the fully-discrete method \eqref{eq10.1} preserves the discrete multi-symplectic conservation law 
\begin{equation*}\begin{aligned}
&\frac{1}{\Delta t }\sum_{k=1}^{s}\sum_{l=1}^{s}\sum_{n=1}^{s}b^{(1)}_{k}b^{(2)}_{l}b^{(3)}_{n}\left(\mathrm{d}(\mathbf E)^1_{kln}  \wedge \mathrm{d}(\mathbf H)^1_{kln}  - \mathrm{d}(\mathbf E)^0_{kln}  \wedge \mathrm{d}(\mathbf H)^0_{kln}  \right)\\
&+\frac{1}{2\Delta x }\sum_{m=1}^{r}\sum_{l=1}^{s}\sum_{n=1}^{s}\tilde{b}_{m}b^{(2)}_{l}b^{(3)}_{n}\left(\mathrm{d}(\mathbf U)^m_{1ln}  \wedge K_{1}\mathrm{d}(\mathbf U)^m_{1ln}  - \mathrm{d}(\mathbf U)^m_{0ln}  \wedge K_{1}\mathrm{d}(\mathbf U)^m_{0ln} \right)\\
&+\frac{1}{2\Delta y }\sum_{m=1}^{r}\sum_{k=1}^{s}\sum_{n=1}^{s}\tilde{b}_{m}b^{(1)}_{k}b^{(3)}_{n}\left(\mathrm{d}(\mathbf U)^m_{k1n}  \wedge K_{2}\mathrm{d}(\mathbf U)^m_{k1n}  - \mathrm{d}(\mathbf U)^m_{k0n}  \wedge K_{2}\mathrm{d}(\mathbf U)^m_{k0n} \right)\\
&+\frac{1}{2\Delta z }\sum_{m=1}^{r}\sum_{k=1}^{s}\sum_{l=1}^{s}\tilde{b}_{m}b^{(1)}_{k}b^{(2)}_{l}\left(\mathrm{d}(\mathbf U)^m_{kl1}  \wedge K_{3}\mathrm{d}(\mathbf U)^m_{kl1}  - \mathrm{d}(\mathbf U)^m_{kl0}  \wedge K_{3}\mathrm{d}(\mathbf U)^m_{kl0} \right) =0
\end{aligned}\end{equation*}
with $(\mathbf U)^{m}_{0ln}=(\mathbf ((\mathbf H)^{m}_{0ln})^\top, ((\mathbf E)^{m}_{0ln})^\top)^\top, (\mathbf U)^{m}_{1ln}=(\mathbf ((\mathbf H)^{m}_{1ln})^\top, ((\mathbf E)^{m}_{1ln})^\top)^\top,(\mathbf U)^{m}_{k0n}=(\mathbf ((\mathbf H)^{m}_{k0n})^\top, ((\mathbf E)^{m}_{k0n})^\top)^\top,$\\  $(\mathbf U)^{m}_{k1n}=(\mathbf ((\mathbf H)^{m}_{k1n})^\top, ((\mathbf E)^{m}_{k1n})^\top)^\top, (\mathbf U)^{m}_{kl0}=(\mathbf ((\mathbf H)^{m}_{kl0})^\top, ((\mathbf E)^{m}_{kl0})^\top)^\top, (\mathbf U)^{m}_{kl1}=(\mathbf ((\mathbf H)^{m}_{kl1})^\top, ((\mathbf E)^{m}_{kl1})^\top)^\top.$
\end{thm}

\begin{rem}
We would like to mention that, in the framework of splitting multi-symplectic Runge--Kutta method,  other multi-symplectic methods can be used to discretize the deterministic Hamiltonian PDE.  By combining with the symplectic Euler method applied to the stochastic subsystem, one can obtain a class of  multi-symplectic methods.
\end{rem}

\section{Multi-symplectic partitioned Runge--Kutta method} 
As we know, symplectic partitioned Runge--Kutta methods, which are the generations of symplectic Runge--Kutta methods, are powerful tools for the construction of symplectic methods for solving stochastic  Hamiltonian ordinary differential equations numerically. 
For separate stochastic  Hamiltonian ordinary differential equations, some symplectic partitioned Runge--Kutta method is explicit, which reduces the computational cost. 
In this section, we construct the third kind of multi-symplectic methods, i.e., multi-symplectic partitioned Runge--Kutta methods,  for nonlinear stochastic wave equation,  stochastic nonlinear Schr\"odinger equation,  stochastic KdV equation and stochastic Maxwell equation by means of the symplectic partitioned Runge--Kutta method in both spatial and temporal directions.  
Further, we present the multi-symplectic conditions.

For the nonlinear stochastic wave equation \eqref{eq2.4}, we proceed to take advantage of $s$-stage partitioned Runge--Kutta method  $(c^{(1)}, A^{(1)}, b^{(1)})$ and $(c^{(2)}, A^{(2)}, b^{(2)})$, i.e., 
\begin{equation}
\label{butchert2}
\begin{array}{c|ccc}
c_{1}^{(1)} &  a_{11}^{(1)} &\dots &a_{1s}^{(1)}\\
\vdots &   \vdots &&\vdots \\
c_{s}^{(1)} &a_{s1}^{(1)} &\dots &a_{ss}^{(1)}\\
\hline
& b_{1}^{(1)}& \dots&  b_{s}^{(1)}
\end{array},\qquad
\begin{array}{c|ccc}
c_{1}^{(2)} &  a_{11}^{(2)} &\dots &a_{1s}^{(2)}\\
\vdots &   \vdots &&\vdots \\
c_{s}^{(2)} &a_{s1}^{(2)} &\dots &a_{ss}^{(2)}\\
\hline
& b_{1}^{(2)}& \dots&  b_{s}^{(2)}
\end{array}, 
\end{equation}
in the spatial direction, and $r$-stage partitioned Runge--Kutta method $( \tilde{c}^{(1)}, \tilde{A}^{(1)}, \tilde{b}^{(1)})$,  $(\tilde{c}^{(2)},  \tilde{A}^{(2)}, \tilde{b}^{(2)}),$ i.e., 
\begin{equation}
\label{butchert3}
\begin{array}{c|ccc}
\tilde c_{1}^{(1)} &  \tilde a_{11}^{(1)} &\dots & \tilde a_{1r}^{(1)}\\
\vdots &   \vdots &&\vdots \\
\tilde c_{r}^{(1)} &\tilde a_{r1}^{(1)} &\dots &\tilde a_{rr}^{(1)}\\
\hline
& \tilde b_{1}^{(1)}& \dots&  \tilde b_{r}^{(1)}
\end{array}, \quad
\begin{array}{c|ccc}
\tilde c_{1}^{(2)} &  \tilde a_{11}^{(2)} &\dots & \tilde a_{1r}^{(2)}\\
\vdots &   \vdots &&\vdots \\
\tilde c_{r}^{(2)} &\tilde a_{r1}^{(2)} &\dots &\tilde a_{rr}^{(2)}\\
\hline
& \tilde b_{1}^{(2)}& \dots&  \tilde b_{r}^{(2)}
\end{array}, 
\end{equation}
together with an $r$-stage Runge--Kutta method 
$( \bar{c}, \bar{A}, \bar{b})$
in the temporal direction, respectively, where $s,r\in\mathbb N_+.$
The resulting fully-discrete method is as follows
\begin{subequations}
\begin{align}
\label{eq5.1}
&U_{i}^{m} = u_{i}^{0}  + \Delta t\sum_{n=1}^{r}\tilde{a}^{(1)}_{nm}V_{i}^{n},\quad
V_{i}^{m} = v_{i}^{0}  + \Delta t\sum_{n=1}^{r}\tilde{a}^{(2)}_{nm}\left(\delta_x\mathcal{W}_{i}^{n} - f(U_{i}^{n})\right) + \Delta W_{i}^{1}\sum_{n=1}^{r}\bar{a}_{nm}g(U_{i}^{n}),\\
 \label{eq5.3}
 &u_{i}^{1} = u_{i}^{0}  + \Delta t\sum_{m=1}^{r}\tilde{b}^{(1)}_{m}V_{i}^{m},\quad 
 v_{i}^{1} = v_{i}^{0}  + \Delta t\sum_{m=1}^{r}\tilde{b}^{(2)}_{m} \left(\delta_x\mathcal{W}_{i}^{m} - f(U_{i}^{m})\right)+ \Delta W_{i}^{1}\sum_{m=1}^{r}\bar{b}_{m}g(U_{i}^{m}),\\
  \label{eq5.5}
 &U_{i}^{m} = u_{0}^{m}  + {\Delta x}\sum_{j=1}^{s} a^{(1)}_{ij}\mathcal{W}_{j}^{m},\quad
 \mathcal{W}_{i}^{m} = w_{0}^{m}  + {\Delta x}\sum_{j=1}^{s} a^{(2)}_{ij}\delta_x\mathcal{W}_{j}^{m},\\
  \label{eq5.7}
 &u_{1}^{m} = u_{0}^{m}  + {\Delta x}\sum_{i=1}^{s} b^{(1)}_{i}\mathcal{W}_{i}^{m},\quad 
 w_{1}^{m} = w_{0}^{m}  + {\Delta x}\sum_{i=1}^{s} b^{(2)}_{i}\delta_x\mathcal{W}_{i}^{m},
\end{align}
\end{subequations}
where   $U_{i}^{m}\approx u(c_i^{(1)}\Delta x, \tilde c_m^{(1)}\Delta t )$, $u_{i}^{0}\approx u(c_i^{(1)}\Delta x,0)$, $u_{i}^{1}\approx u(c_i^{(1)}\Delta x,\Delta t)$, $u_{0}^{m}\approx u(0,\tilde c_m^{(1)}\Delta t)$, $u_{1}^{m}\approx u(\Delta x,\tilde c_m^{(1)}\Delta t)$,  etc., with $c_i^{(1)} =\sum_{j=1}^{s}a_{ij}^{(1)}$, $\tilde c_m ^{(1)}=\sum_{n=1}^{r}\tilde{a}_{mn}^{(1)}$ for $i =1,\dots, s, m=1,\dots, r$.

\begin{thm} \label{thm7}
Suppose that
\begin{subequations}
\begin{align}
\label{eq5.9}
&\bar{a}_{nm}\tilde{b}^{(1)}_{m} + \tilde{a}^{(1)}_{m,n} \bar{b}_{n}- \tilde{b}^{(1)}_{m}\bar{b}_{n} = 0,\\
 \label{eq5.10}
&\tilde{a}^{(2)}_{nm}\tilde{b}^{(1)}_{m} + \tilde{a}^{(1)}_{m,n}\tilde{b}^{(2)}_{n} - \tilde{b}^{(1)}_{m}\tilde{b}^{(2)}_{n}  = 0,\\
 \label{eq5.11}
& a^{(2)}_{ij}b^{(1)}_{i} + a^{(1)}_{j, i}b^{(2)}_{j} - b^{(1)}_{i}b^{(2)}_{j}  = 0,
\end{align}
\end{subequations} 
for $ 1\leq i, j \leq s, ~ 1\leq m, n \leq r.$  Then the fully-discrete method \eqref{eq5.1}-\eqref{eq5.7} admits the discrete multi-symplectic conservation law 
\begin{equation} \label{eq5.12}
\sum_{i=1}^{s}b^{(2)}_{i}\frac{1}{\Delta t }\left( \mathrm{d}u_{i}^{1} \wedge \mathrm{d}v_{i}^{1} - \mathrm{d}u_{i}^{0} \wedge \mathrm{d}v_{i}^{0} \right) - \sum_{m=1}^{r}\tilde{b}^{(2)}_{m} \frac{1}{{\Delta x} } \left( \mathrm{d}u_{1}^{m} \wedge \mathrm{d}w_{1}^{m} - \mathrm{d}u_{0}^{m} \wedge \mathrm{d}w_{0}^{m} \right)=0.
\end{equation} 
\end{thm}
\noindent{\textbf{Proof. }}
From  \eqref{eq5.3} it follows that
\begin{equation*}\begin{aligned}
&\left( \mathrm{d}u_{i}^{1} \wedge \mathrm{d}v_{i}^{1} - \mathrm{d}u_{i}^{0} \wedge \mathrm{d}v_{i}^{0} \right)\\
= &\Big( \mathrm{d}u_{i}^{0}  + \Delta t\sum_{m=1}^{r}\tilde{b}^{(1)}_{m}\mathrm{d}V_{i}^{m}  \Big)\wedge  \Big(\mathrm{d}  v_{i}^{0}  + \Delta t\sum_{m=1}^{r}\tilde{b}^{(2)}_{m}\mathrm{d} \left(\delta_x\mathcal{W}_{i}^{m} - f(U_{i}^{m})\right) + \Delta W_{i}^{1}\sum_{m=1}^{r}\bar{b}_{m}\mathrm{d}  g(U_{i}^{m})\Big)-\mathrm{d}u_{i}^{0} \wedge \mathrm{d}v_{i}^{0}.
\end{aligned}\end{equation*}
Based on \eqref{eq5.1} we derive
\begin{equation*}\begin{aligned}
\mathrm{d}u_{i}^{0} =  \mathrm{d} U_{i}^{m} - \Delta t\sum_{n=1}^{r}\tilde{a}^{(1)}_{nm}\mathrm{d} V_{i}^{n},\quad \mathrm{d} v_{i}^{0} = \mathrm{d} V_{i}^{m}  - \Delta t\sum_{n=1}^{r}\tilde{a}^{(2)}_{nm}\mathrm{d} \left(\delta_x\mathcal{W}_{i}^{n} - f(U_{i}^{n})\right) - \Delta W_{i}^{1}\sum_{n=1}^{r}\bar{a}_{nm}\mathrm{d} g(U_{i}^{n}),
\end{aligned}\end{equation*}
which implies
\begin{equation*}\begin{aligned}
&\frac{1}{\Delta t}\left( \mathrm{d}u_{i}^{1} \wedge \mathrm{d}v_{i}^{1} - \mathrm{d}u_{i}^{0} \wedge \mathrm{d}v_{i}^{0} \right)\\
=& -\Delta W_{i}^{1}\sum_{m,n=1}^{r}\left( \bar{a}_{nm}\tilde{b}^{(1)}_{m} + \tilde{a}^{(1)}_{m,n}\bar{b}_{n} - \tilde{b}^{(1)}_{m} \bar{b}_{n}\right)\mathrm{d}V_{i}^{m}\wedge \mathrm{d}g(U_{i}^{n})\\
& -\Delta t\sum_{m,n=1}^{r}\left( \tilde{a}^{(2)}_{nm}\tilde{b}^{(1)}_{m} + \tilde{a}^{(1)}_{m,n}\tilde{b}^{(2)}_{n} - \tilde{b}^{(1)}_{m}\tilde{b}^{(2)}_{n} \right) \mathrm{d} V_{i}^{m}\wedge
\mathrm{d}(\delta_x\mathcal{W}_{i}^{n} - f(U_{i}^{n}))  +\sum_{m=1}^{r}\tilde{b}^{(2)}_{m}\mathrm{d}U_{i}^{m}\wedge \mathrm{d} (\delta_x\mathcal{W}_{i}^{m}).
\end{aligned}\end{equation*}
Making use of \eqref{eq5.9} and \eqref{eq5.10} leads to
\begin{equation}\label{eq5.14}
\frac{1}{\Delta t}\left( \mathrm{d}u_{i}^{1} \wedge \mathrm{d}v_{i}^{1} - \mathrm{d}u_{i}^{0} \wedge \mathrm{d}v_{i}^{0} \right)
= \sum_{m=1}^{r}\tilde{b}^{(2)}_{m}\mathrm{d}U_{i}^{m}\wedge \mathrm{d} (\delta_x\mathcal{W}_{i}^{m}).
\end{equation}
Similarly, by means of  \eqref{eq5.5}, we derive 
\begin{equation*}\begin{aligned}
&\frac{1}{{\Delta x}}\left( \mathrm{d}u_{1}^{m} \wedge \mathrm{d}w_{1}^{m} - \mathrm{d}u_{0}^{m} \wedge \mathrm{d}w_{0}^{m} \right)\\
=& -{\Delta x}\sum_{i, j =1}^{s}\left( a^{(2)}_{ij}b^{(1)}_{i} + a^{(1)}_{j,i}b^{(2)}_{j} - b^{(1)}_{i}b^{(2)}_{j} \right)
\mathrm{d}\mathcal{W}_{i}^{m} \wedge \mathrm{d} (\delta_x\mathcal{W}_{j}^{m} )  
+  \sum_{i=1}^{s}b^{(2)}_{i}\mathrm{d}U_{i}^{m} \wedge \mathrm{d}(\delta_x\mathcal{W}_{i}^{m}) .
\end{aligned}\end{equation*}
By utilizing \eqref{eq5.11}, we obtain
\begin{equation}\label{eq5.15}
\frac{1}{{\Delta x}}\left( \mathrm{d}u_{1}^{m} \wedge \mathrm{d}w_{1}^{m} - \mathrm{d}u_{0}^{m} \wedge \mathrm{d}w_{0}^{m} \right)
= \sum_{i=1}^{s}b^{(2)}_{i}\mathrm{d}U_{i}^{m} \wedge \mathrm{d}(\delta_x\mathcal{W}_{i}^{m}) .
\end{equation}
Combining  \eqref{eq5.14} and  \eqref{eq5.15}, we have
\begin{equation*}
\sum_{i=1}^{s}b^{(2)}_{i}\frac{1}{\Delta t }\left( \mathrm{d}u_{i}^{1} \wedge \mathrm{d}v_{i}^{1} - \mathrm{d}u_{i}^{0} \wedge \mathrm{d}v_{i}^{0} \right) - \sum_{m=1}^{r}\tilde{b}^{(2)}_{m} \frac{1}{{\Delta x} } \left( \mathrm{d}u_{1}^{m} \wedge \mathrm{d}w_{1}^{m} - \mathrm{d}u_{0}^{m} \wedge \mathrm{d}w_{0}^{m} \right)=0,
\end{equation*}
which completes the proof.\hfill$\qed$

\begin{figure}[h]
	\centering
	\subfigure{
		\begin{minipage}{15cm}
			\centering
\includegraphics[height=4cm,width=4.5cm]{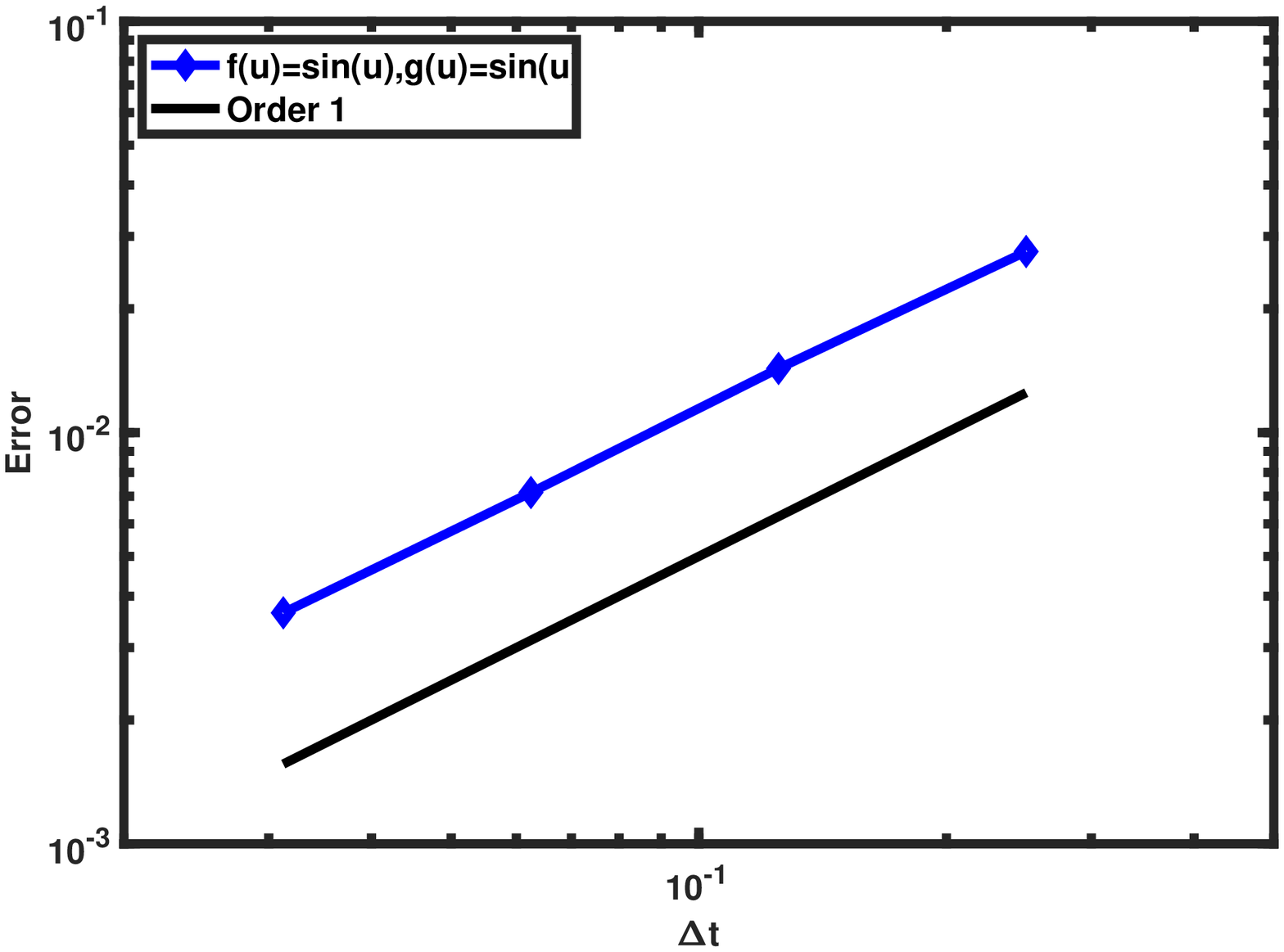}
\includegraphics[height=4cm,width=4.5cm]{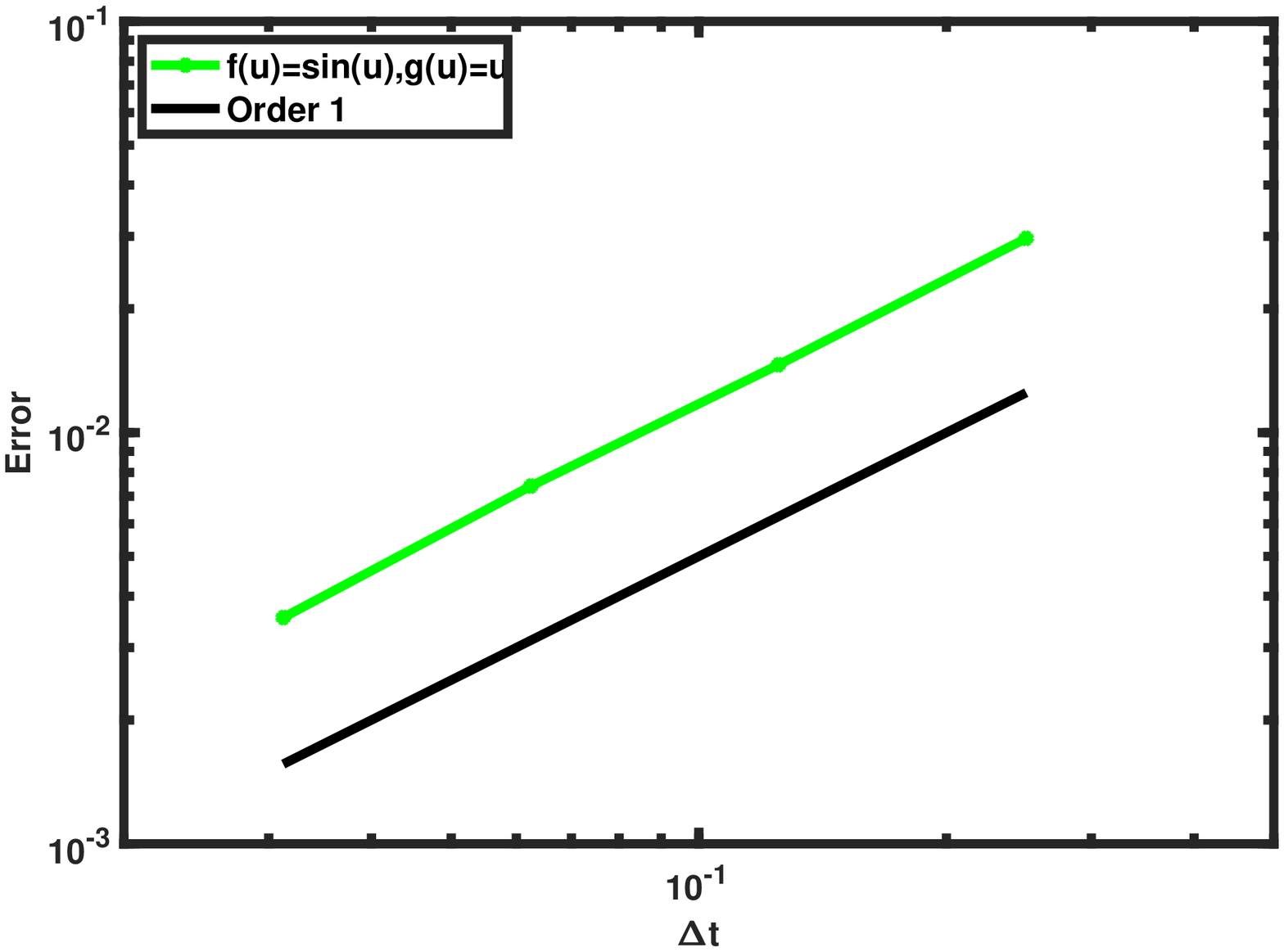}
\includegraphics[height=4cm,width=4.5cm]{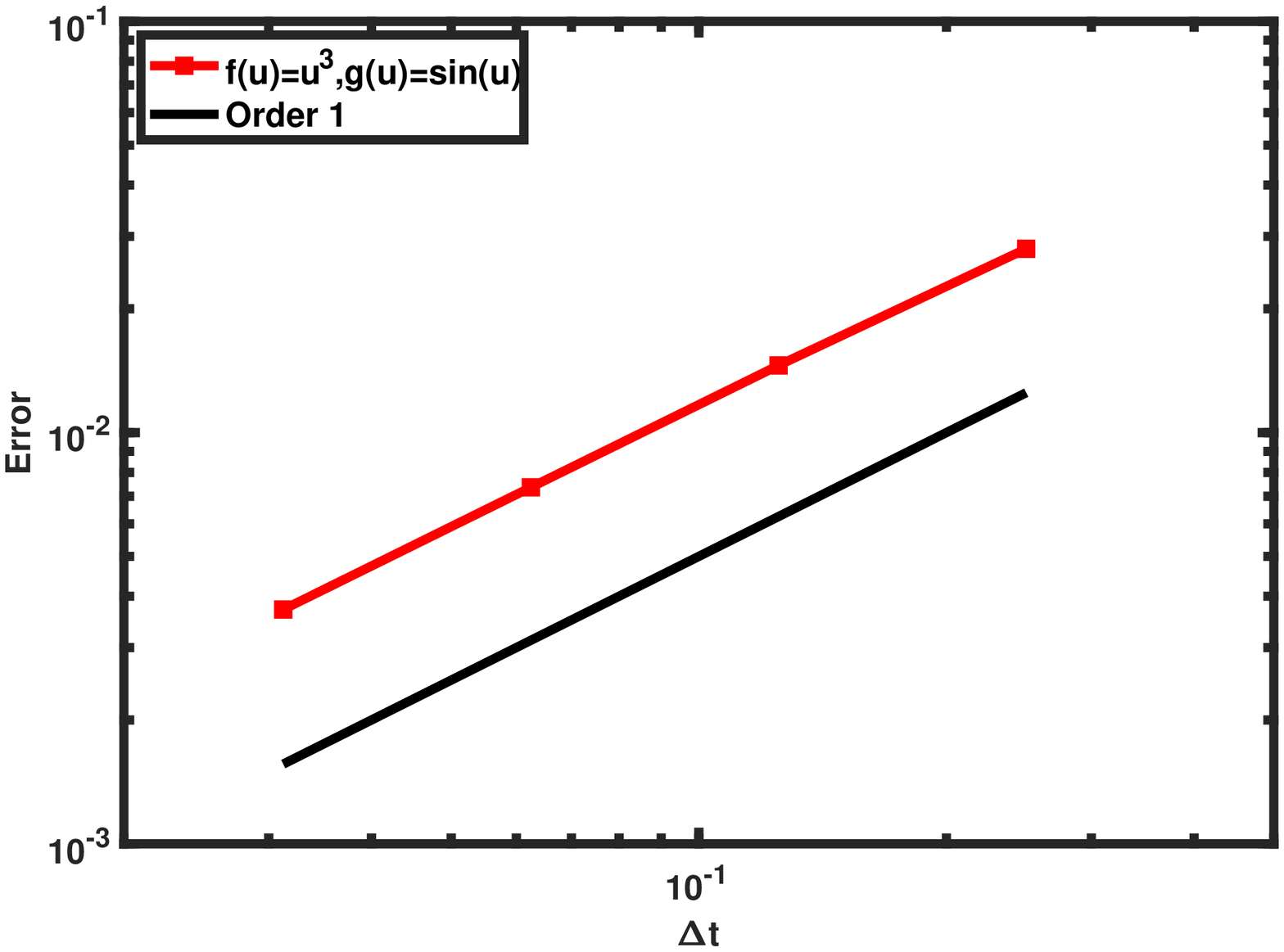}
		\end{minipage}
	}
\caption{Mean-square convergence order of \eqref{eq5.1ex} in temporal direction in the cases of (1) $f(u)= \sin(u), g(u) = \sin(u)$ (2) $f(u)= \sin(u), g(u) = u$ and (3) $f(u)= u^3, g(u) = \sin(u).$} 
	\label{fig4}
\end{figure}

\begin{ex}
Let $s=r=1$ and the Butcher tableaux in both \eqref{butchert2} and \eqref{butchert3} be
$
\begin{array}{c|c}
\frac 12 &  \frac 12\\
\hline
&1
\end{array},
$
we obtain an explicit numerical method for the nonlinear stochastic wave equation \eqref{eq2.4} as follows
\begin{equation}
\label{eq5.1ex}
\begin{aligned}
&U_{1}^{1} = u_{1}^{0}  + \frac{\Delta t} 2V_{1}^{1},\quad
V_{1}^{1} = v_{1}^{0}  + \frac{\Delta t}2\left(\delta_x\mathcal{W}_{1}^{1} - f(U_{1}^{1})\right) + \frac {\Delta W_{1}^{1}} 2g(U_{1}^{1}),\\
&u_{1}^{1} = u_{1}^{0}  + \Delta tV_{1}^{1},\quad 
v_{1}^{1} = v_{1}^{0}  + \Delta t \left(\delta_x\mathcal{W}_{1}^{1} - f(U_{1}^{1})\right)+ \Delta W_{1}^{1}g(U_{1}^{1}),\\
&U_{1}^{1} = u_{0}^{1}  + \frac {\Delta x}2\mathcal{W}_{1}^{1},\quad
\mathcal{W}_{1}^{1} = w_{0}^{1}  +  \frac {\Delta x}2\delta_x\mathcal{W}_{1}^{1},\\
&u_{1}^{1} = u_{0}^{1}  + {\Delta x}\mathcal{W}_{1}^{1},\quad 	w_{1}^{1} = w_{0}^{1}  + {\Delta x}\delta_x\mathcal{W}_{1}^{1}.
\end{aligned}
\end{equation}
Now we perform experiments by applying \eqref{eq5.1ex} to the 1-dimensional nonlinear stochastic wave equation,
and consider the same problem as in Example \ref{ex2.1}.  Table 3 shows the mean-square convergence error  against $\Delta t = 2^{-s}, s = 2, 3, 4, 5$ on log-log scale at time $T = 1.$ 
The exact solution is regarded as the numerical approximation obtained by a fine mesh with $\Delta t=2^{-8}, \Delta x = 2^{-7}\pi$. 
Fig. \ref{fig4} shows  that the mean-square convergence order of the proposed numerical method is $1$ in time. 		
\end{ex}

\begin{table}[htbp]
	\setlength{\abovecaptionskip}{0pt}
	\setlength{\belowcaptionskip}{3pt}
	\centering
       \caption{\label{tab3}Mean-square errors of \eqref{eq5.1ex} in time.}
	\begin{tabularx}{0.84\textwidth}{|c|c|c|c|}
		\hline
		& \multicolumn{1}{c|}{$f(u)= \sin(u), g(u) = \sin(u)$}  & \multicolumn{1}{c|}{$f(u)= \sin(u), g(u) = u$}& \multicolumn{1}{c|}{$f(u)= u^3, g(u) = \sin(u)$}\\
		\hline
		$\Delta t$    &$L^2$~error          &$L^2$~error         &$L^2$~error          \\
		\hline
		$2^{-2}$        & 2.7559e-02         & 2.9672e-02          & 2.7987e-02            \\
		$2^{-3}$        & 1.4317e-02         &1.4617e-02           &1.4566e-02             \\
		$2^{-4}$        & 7.1487e-03         &7.4123e-03           & 7.3600e-03            \\
		$2^{-5}$        & 3.6446e-03         & 3.5495e-03          & 3.7130e-03             \\
		\hline
	\end{tabularx}
\end{table}

Analogous to the nonlinear stochastic wave equation, for the  stochastic nonlinear Schr\"odinger equation \eqref{eq2.6}, based on  $s$-stage partitioned Runge--Kutta methods 
	\begin{equation}
	\label{butchert4}
	\begin{array}{c|ccc}
	c_{1}^{(1)} &  a_{11}^{(1)} &\dots &a_{1s}^{(1)}\\
	\vdots &   \vdots &&\vdots \\
	c_{s}^{(1)} &a_{s1}^{(1)} &\dots &a_{ss}^{(1)}\\
	\hline
	& b_{1}^{(1)}& \dots&  b_{s}^{(1)}
	\end{array},\quad
	\begin{array}{c|ccc}
	c_{1}^{(2)} &  a_{11}^{(2)} &\dots &a_{1s}^{(2)}\\
	\vdots &   \vdots &&\vdots \\
	c_{s}^{(2)} &a_{s1}^{(2)} &\dots &a_{ss}^{(2)}\\
	\hline
	& b_{1}^{(2)}& \dots&  b_{s}^{(2)}
	\end{array}, \quad
		\begin{array}{c|ccc}
	c_{1}^{(2)} &  a_{11}^{(3)} &\dots &a_{1s}^{(3)}\\
	\vdots &   \vdots &&\vdots \\
	c_{s}^{(3)} &a_{s1}^{(3)} &\dots &a_{ss}^{(3)}\\
	\hline
	& b_{1}^{(3)}& \dots&  b_{s}^{(3)}
	\end{array}, \quad
		\begin{array}{c|ccc}
	c_{1}^{(4)} &  a_{11}^{(4)} &\dots &a_{1s}^{(4)}\\
	\vdots &   \vdots &&\vdots \\
	c_{s}^{(4)} &a_{s1}^{(4)} &\dots &a_{ss}^{(4)}\\
	\hline
	& b_{1}^{(4)}& \dots&  b_{s}^{(4)}
	\end{array}, 
	\end{equation}
and $r$-stage partitioned Runge--Kutta methods 
\begin{equation}\label{butchert5}
\begin{array}{c|ccc}
\tilde c_{1}^{(1)} &  \tilde a_{11}^{(1)} &\dots & \tilde a_{1r}^{(1)}\\
\vdots &   \vdots &&\vdots \\
\tilde c_{r}^{(1)} &\tilde a_{r1}^{(1)} &\dots &\tilde a_{rr}^{(1)}\\
\hline
& \tilde b_{1}^{(1)}& \dots&  \tilde b_{r}^{(1)}
\end{array}, \quad
\begin{array}{c|ccc}
\tilde c_{1}^{(2)} &  \tilde a_{11}^{(2)} &\dots & \tilde a_{1r}^{(2)}\\
\vdots &   \vdots &&\vdots \\
\tilde c_{r}^{(2)} &\tilde a_{r1}^{(2)} &\dots &\tilde a_{rr}^{(2)}\\
\hline
& \tilde b_{1}^{(2)}& \dots&  \tilde b_{r}^{(2)}
\end{array},\quad
\begin{array}{c|ccc}
\bar c_{1}^{(1)} &  \bar a_{11}^{(1)} &\dots & \bar a_{1r}^{(1)}\\
\vdots &   \vdots &&\vdots \\
\bar c_{r}^{(1)} &\bar a_{r1}^{(1)} &\dots &\bar a_{rr}^{(1)}\\
\hline
& \bar b_{1}^{(1)}& \dots&  \bar b_{r}^{(1)}
\end{array}, \quad
\begin{array}{c|ccc}
\bar c_{1}^{(2)} &  \bar a_{11}^{(2)} &\dots & \bar a_{1r}^{(2)}\\
\vdots &   \vdots &&\vdots \\
\bar c_{r}^{(2)} &\bar a_{r1}^{(2)} &\dots &\bar a_{rr}^{(2)}\\
\hline
& \bar b_{1}^{(2)}& \dots&  \bar b_{r}^{(2)}
\end{array},
\end{equation}
where $s, r\in\mathbb N_+$, 
we deduce the following 
fully-discrete method
\begin{align}
\label{eq5.30}
 &Q_{i}^{m} = q_{i}^{0}  + \Delta t\sum_{n=1}^{r}\tilde{a}^{(1)}_{nm}\left(\delta_xV_{i}^{n} + ((P_{i}^{n})^2+(Q_{i}^{n})^2)P_{i}^{n}\right) - \Delta W_{i}^{1}\sum_{n=1}^{r}\bar{a}^{(1)}_{nm}P_{i}^{n},\nonumber\\
 &P_{i}^{m} = p_{i}^{0}  + \Delta t\sum_{n=1}^{r}\tilde{a}^{(2)}_{nm}\left(-\delta_x\mathcal{W}_{i}^{n} - ((P_{i}^{n})^2+(Q_{i}^{n})^2)Q_{i}^{n}\right) + \Delta W_{i}^{1}\sum_{n=1}^{r}\bar{a}^{(2)}_{nm}Q_{i}^{n},\nonumber\\
 & q_{i}^{1} = q_{i}^{0}  + \Delta t\sum_{m=1}^{r}\tilde{b}^{(1)}_{m}\left(\delta_xV_{i}^{m} + ((P_{i}^{m})^2+(Q_{i}^{m})^2)P_{i}^{m}\right) - \Delta W_{i}^{1}\sum_{m=1}^{r}\bar{b}^{(1)}_{m}P_{i}^{m},\nonumber\\
 &p_{i}^{1} = p_{i}^{0}  + \Delta t\sum_{m=1}^{r}\tilde{b}^{(2)}_{m}\left(-\delta_x\mathcal{W}_{i}^{m} - ((P_{i}^{m})^2+(Q_{i}^{m})^2)Q_{i}^{m}\right) + \Delta W_{i}^{1}\sum_{m=1}^{r}\bar{b}^{(2)}_{m}Q_{i}^{m},\nonumber\\
 &V_{i}^{m} = v_{0}^{m}  + {\Delta x}\sum_{j=1}^{s} a^{(1)}_{ij}\delta_xV_{j}^{m},\quad \mathcal{W}_{i}^{m} = w_{0}^{m}  + {\Delta x}\sum_{j=1}^{s} a^{(2)}_{ij}\delta_x\mathcal{W}_{j}^{m},\\
 &P_{i}^{m} = p_{0}^{m}  + {\Delta x}\sum_{j=1}^{s} a^{(3)}_{ij}V_{j}^{m},\quad Q_{i}^{m} = q_{0}^{m}  + {\Delta x}\sum_{j=1}^{s} a^{(4)}_{ij}\mathcal{W}_{j}^{m},\nonumber\\
 &v_{1}^{m} = v_{0}^{m}  + {\Delta x}\sum_{i=1}^{s} b^{(1)}_{i}\delta_xV_{i}^{m},\quad w_{1}^{m} = w_{0}^{m}  + {\Delta x}\sum_{i=1}^{s} b^{(2)}_{i}\delta_x\mathcal{W}_{i}^{m},\nonumber\\
 &p_{1}^{m} = p_{0}^{m}  + {\Delta x}\sum_{i=1}^{s} b^{(3)}_{i}V_{i}^{m},\quad q_{1}^{m} = q_{0}^{m}  + {\Delta x}\sum_{i=1}^{s} b^{(4)}_{i}\mathcal{W}_{i}^{m},\nonumber
\end{align}
where   $Q_{i}^{m}\approx q(c_i^{(4)}\Delta x, \tilde c_m^{(1)}\Delta t )$, $q_{i}^{0}\approx q(c_i^{(4)}\Delta x,0)$, $q_{i}^{1}\approx q(c_i^{(4)}\Delta x,\Delta t)$,  $q_{0}^{m}\approx q(0,\tilde c_m^{(1)}\Delta t)$, $q_{1}^{m}\approx q(\Delta x,\tilde c_m^{(1)}\Delta t)$, etc., with $c_i^{(4)} =\sum_{j=1}^{s}a_{ij}^{(4)}$, $\tilde c_m^{(1)} =\sum_{n=1}^{r}\tilde{a}_{mn}^{(1)}$ for $i =1,\dots, s$,  $m=1,\dots, r$.  Similar to Theorem \ref{thm7}, we obtain the following result.
\begin{thm} \label{thm8}
If the following conditions
\begin{align*}
&\tilde{a}^{(1)}_{m,n}\tilde{b}^{(2)}_{n} + \tilde{b}^{(1)}_{m}\tilde{a}^{(2)}_{nm} - \tilde{b}^{(1)}_{m}\tilde{b}^{(2)}_{n}=0,\\
&\bar{a}^{(1)}_{m,n}\tilde{b}^{(2)}_{n} + \bar{b}^{(1)}_{m}\tilde{a}^{(2)}_{nm} - \bar{b}^{(1)}_{m}\tilde{b}^{(2)}_{n}=0,\\
 &\tilde{a}^{(1)}_{m,n}\bar{b}^{(2)}_{n} + \tilde{b}^{(1)}_{m}\bar{a}^{(2)}_{nm} - \tilde{b}^{(1)}_{m}\bar{b}^{(2)}_{n}=0,\\
&\bar{a}^{(1)}_{m,n}\bar{b}^{(2)}_{n} + \bar{b}^{(1)}_{m}\bar{a}^{(2)}_{nm} - \bar{b}^{(1)}_{m}\bar{b}^{(2)}_{n}=0,\\
&a^{(3)}_{ij}b^{(1)}_{j} + b^{(3)}_{i}a^{(1)}_{j, i} - b^{(3)}_{i}b^{(1)}_{j}=0,\\
&a^{(2)}_{ij}b^{(4)}_{i} + b^{(2)}_{j}a^{(4)}_{j, i} - b^{(4)}_{i}b^{(2)}_{j}=0,\\
&b^{(1)}_{i} =b^{(2)}_{i},~~ \tilde{b}^{(1)}_{m} =\tilde{b}^{(2)}_{m},
\end{align*}
for $ 1\leq i, j \leq s, ~ 1\leq m, n \leq r,$  hold, the fully-discrete method \eqref{eq5.30} possesses the discrete multi-symplectic conservation law
\begin{equation}\begin{aligned}
&\sum_{i=1}^{s}b^{(1)}_{i}\frac{1}{\Delta t }\left( \mathrm{d}q_{i}^{1} \wedge \mathrm{d}p_{i}^{1} - \mathrm{d}q_{i}^{0} \wedge \mathrm{d}p_{i}^{0} \right)\\
&+\sum_{m=1}^{r}\tilde{b}^{(1)}_{m} \frac{1}{{\Delta x}}\left( \mathrm{d}p_{1}^{m} \wedge \mathrm{d}v_{1}^{m} - \mathrm{d}p_{0}^{m} \wedge \mathrm{d}v_{0}^{m} +\mathrm{d}q_{1}^{m} \wedge \mathrm{d}w_{1}^{m} - \mathrm{d}q_{0}^{m} \wedge \mathrm{d}w_{0}^{m} 
\right)=0.
\end{aligned}\end{equation} 
\end{thm}

In the case of the stochastic KdV equation \eqref{eq2.10}, exploiting similar procedures that applying   
$s$-stage partitioned Runge--Kutta methods \eqref{butchert4}, and $r$-stage Runge--Kutta methods \eqref{butchert5} with $s, r\in\mathbb N_+$ yields the following multi-symplectic method 
\begin{align}
\label{eq5.50}
 &U_{i}^{m} = u_{i}^{0}  + \Delta t\sum_{n=1}^{r}\left(-2\delta_xV_{i}^{n}\right)\tilde{a}^{(1)}_{nm} + \sum_{n=1}^{r}2\lambda \Delta W_{i}^{1}\bar{a}_{nm},\nonumber\\
 &\mathcal{P}_{i}^{m} = \rho_{i}^{0}  + \Delta t\sum_{n=1}^{r}\left(-2 \beta\delta_x\mathcal{W}_{i}^{n}+2V_{i}^{n}-(U_{i}^{n})^2 \right)\tilde{a}^{(2)}_{nm} ,\nonumber\\ 
 & u_{i}^{1} = u_{i}^{0}  + \Delta t\sum_{m=1}^{r}\left(-2\delta_xV_{i}^{m}\right)\tilde{b}^{(1)}_{m} +\sum_{m=1}^{r}2\lambda \Delta W_{i}^{1}\bar{b}_{m},\nonumber\\
 &\rho_{i}^{1} = \rho_{i}^{0}  + \Delta t\sum_{m=1}^{r}\left(-2 \beta\delta_x\mathcal{W}_{i}^{m}+2V_{i}^{m}-(U_{i}^{m})^2\right)\tilde{b}^{(2)}_{m},\nonumber\\
 &U_{i}^{m} = u_{0}^{m}  + {\Delta x}\sum_{j=1}^{s} a^{(1)}_{ij}\mathcal{W}_{j}^{m},\quad \mathcal{P}_{i}^{m} = \rho_{0}^{m}  + {\Delta x}\sum_{j=1}^{s} a^{(2)}_{ij}U_{j}^{m},\\
 &V_{i}^{m} = v_{0}^{m}  + {\Delta x}\sum_{j=1}^{s} a^{(3)}_{ij}\delta_xV_{j}^{m},
 \quad
 \mathcal{W}_{i}^{m} = w_{0}^{m}  + {\Delta x}\sum_{j=1}^{s} a^{(4)}_{ij}\delta_x\mathcal{W}_{j}^{m},\nonumber\\ 
 &u_{1}^{m} = u_{0}^{m}  + {\Delta x}\sum_{i=1}^{s} b^{(1)}_{i}\mathcal{W}_{i}^{m},\quad 
 \rho_{1}^{m} = \rho_{0}^{m}  + {\Delta x}\sum_{i=1}^{s} b^{(2)}_{i}U_{i}^{m},\nonumber\\
 &v_{1}^{m} = v_{0}^{m}  + {\Delta x}\sum_{i=1}^{s} b^{(3)}_{i}\delta_xV_{i}^{m},\quad 
 w_{1}^{m} = w_{0}^{m}  + {\Delta x}\sum_{i=1}^{s} b^{(4)}_{i}\delta_x\mathcal{W}_{i}^{m},\nonumber
\end{align}
where
\begin{align*}
&\tilde{a}^{(1)}_{m,n}\tilde{b}^{(2)}_{n} + \tilde{b}^{(1)}_{m}\tilde{a}^{(2)}_{nm} - \tilde{b}^{(1)}_{m}\tilde{b}^{(2)}_{n}=0,\\
&\bar{a}_{m,n}\tilde{b}^{(2)}_{n} + \bar{b}_{m}\tilde{a}^{(2)}_{nm} - \bar{b}_{m}\tilde{b}^{(2)}_{n}=0,\\
&a^{(3)}_{ij}b^{(2)}_{i} + b^{(3)}_{j}a^{(2)}_{j, i} - b^{(2)}_{i}b^{(3)}_{j}=0,\\
&a^{(1)}_{ij}b^{(4)}_{i} + b^{(1)}_{j}a^{(4)}_{j, i} - b^{(4)}_{i}b^{(1)}_{j}=0,\\
&b^{(2)}_{i} = b^{(3)}_{i},~b^{(2)}_{i} = b^{(4)}_{i},~ \tilde{b}^{(1)}_{m} = \tilde{b}^{(2)}_{m},
\end{align*} and  $U_{i}^{m}\approx u(c^{(1)}_i\Delta x, \tilde c^{(1)}_m\Delta t )$, $u_{i}^{0}\approx u(c^{(1)}_i\Delta x,0)$, $u_{i}^{1}\approx u(c^{(1)}_i\Delta x,\Delta t)$, $u_{0}^{m}\approx u(0,\tilde c^{(1)}_m\Delta t)$, $u_{1}^{m}\approx u(\Delta x,\tilde c^{(1)}_m\Delta t)$, etc., with $c^{(1)}_i =\sum_{j=1}^{s}a_{ij}^{(1)}$, $\tilde c^{(1)}_m =\sum_{n=1}^{r}\tilde{a}_{mn}^{(1)}$ for $i =1,\dots, s$,  $m=1,\dots, r.$ 
Making use of the same arguments as in the proof of Theorem \ref{thm7}, the associated discrete multi-symplectic conservation law reads
\begin{equation}\begin{aligned}
&\sum_{i=1}^{s}\frac{b^{(2)}_{i}}{\Delta t }\left(\mathrm{d}\rho_{i}^{1} \wedge \mathrm{d}u_{i}^{1} - \mathrm{d}\rho_{i}^{0} \wedge \mathrm{d}u_{i}^{0} \right)\\
&+\sum_{m=1}^{r} \frac{\tilde{b}^{(2)}_{m}}{{\Delta x}}\left( 2\mathrm{d}\rho_{1}^{m} \wedge \mathrm{d}v_{1}^{m} - 2\mathrm{d}\rho_{0}^{m} \wedge \mathrm{d}v_{0}^{m} +2\beta\mathrm{d}w_{1}^{m} \wedge \mathrm{d}u_{1}^{m} -2\beta \mathrm{d}w_{0}^{m} \wedge \mathrm{d}u_{0}^{m} 
\right)=0.
\end{aligned}\end{equation}


For the stochastic Maxwell equation \eqref{Maxwell's equation}, adopting  $s$-stage partitioned Runge--Kutta methods with Butcher tableaux $(c^{(1)}, A^{(1)}, b^{(1)})$, $(c^{(2)}, A^{(2)}, b^{(2)})$ in $x$ direction,  $(c^{(3)}, A^{(3)}, b^{(3)})$ and $(c^{(4)}, A^{(4)}, b^{(4)})$ in $y$ direction, which are presented in \eqref{butchert4},  $(c^{(5)}, A^{(5)}, b^{(5)})$ and $(c^{(6)}, A^{(6)}, b^{(6)})$ in $z$ direction as follows
\begin{equation} 
	\begin{array}{c|ccc}
	c_{1}^{(5)} &  a_{11}^{(5)} &\dots &a_{1s}^{(5)}\\
	\vdots &   \vdots &&\vdots \\
	c_{s}^{(5)} &a_{s1}^{(5)} &\dots &a_{ss}^{(5)}\\
	\hline
	& b_{1}^{(5)}& \dots&  b_{s}^{(5)}
	\end{array}, \quad
		\begin{array}{c|ccc}
	c_{1}^{(6)} &  a_{11}^{(6)} &\dots &a_{1s}^{(6)}\\
	\vdots &   \vdots &&\vdots \\
	c_{s}^{(6)} &a_{s1}^{(6)} &\dots &a_{ss}^{(6)}\\
	\hline
	& b_{1}^{(6)}& \dots&  b_{s}^{(6)}
	\end{array}, 
	\end{equation} 	
and $r$-stage partitioned Runge--Kutta methods \eqref{butchert5}
in the temporal direction, respectively, where $s, r\in\mathbb N_+.$ The resulting numerical method is as follows
\begin{align}\label{MMM}
 &\mathscr{D}_{1}(\mathbf H)^m_{kln} = \mathscr{D}_{1}(\mathbf H)^m_{0ln} + \Delta x  \sum_{j=1}^{s} a_{kj}^{(1)}\mathscr{D}_{1}(\delta_x \mathbf H)^m_{jln}, 
\quad \mathscr{D}_{2}(\mathbf H)^m_{kln} = \mathscr{D}_{2}(\mathbf H)^m_{k0n} + \Delta y  \sum_{j=1}^{s} a_{lj}^{(2)}\mathscr{D}_{2}(\delta_y \mathbf H)^m_{kjn}, \nonumber\\
  &\mathscr{D}_{3}(\mathbf H)^m_{kln} = \mathscr{D}_{3}(\mathbf H)^m_{kl0} + \Delta z  \sum_{j=1}^{s} a_{nj}^{(3)}\mathscr{D}_{3}(\delta_z \mathbf H)^m_{klj}, 
\quad \mathscr{D}_{1}(\mathbf E)^m_{kln} = \mathscr{D}_{1}(\mathbf E)^m_{0ln} + \Delta x  \sum_{j=1}^{s} a_{kj}^{(4)}\mathscr{D}_{1}(\delta_x \mathbf E)^m_{jln}, \nonumber\\
  &\mathscr{D}_{2}(\mathbf E)^m_{kln} = \mathscr{D}_{2}(\mathbf E)^m_{k0n} + \Delta y  \sum_{j=1}^{s} a_{lj}^{(5)}\mathscr{D}_{2}(\delta_y \mathbf E)^m_{kjn}, 
\quad \mathscr{D}_{3}(\mathbf E)^m_{kln} = \mathscr{D}_{3}(\mathbf E)^m_{kl0} + \Delta z  \sum_{j=1}^{s} a_{nj}^{(6)}\mathscr{D}_{3}(\delta_z \mathbf E)^m_{klj}, \nonumber\\
& (\mathbf H)^m_{kln} = (\mathbf H)^0_{kln} + \Delta t  \sum_{i=1}^{r} \tilde{a}_{mi}^{(1)}\left( -\mathscr{D}_{1}(\delta_x \mathbf E)^i_{kln}  -\mathscr{D}_{2}(\delta_y \mathbf E)^i_{kln}  -\mathscr{D}_{3}(\delta_z \mathbf E)^i_{kln} \right)  + \lambda \Delta W_{kln}^{1}\sum_{i=1}^{r}\bar{a}_{mi}^{(1)}(\mathbf E)^i_{kln}, \nonumber\\
& (\mathbf E)^m_{kln} = (\mathbf E)^0_{kln} + \Delta t  \sum_{i=1}^{r} \tilde{a}_{mi}^{(2)}\left( \mathscr{D}_{1}(\delta_x \mathbf H)^i_{kln}  +\mathscr{D}_{2}(\delta_y \mathbf H)^i_{kln}  + \mathscr{D}_{3}(\delta_z \mathbf H)^i_{kln} \right)  - \lambda \Delta W_{kln}^{1}\sum_{i=1}^{r}\bar{a}_{mi}^{(2)}(\mathbf H)^i_{kln}, \\
 &\mathscr{D}_{1}(\mathbf H)^m_{1ln} = \mathscr{D}_{1}(\mathbf H)^m_{0ln} + \Delta x  \sum_{k=1}^{s} b_{k}^{(1)}\mathscr{D}_{1}(\delta_x \mathbf H)^m_{kln}, 
\quad \mathscr{D}_{2}(\mathbf H)^m_{k1n} = \mathscr{D}_{2}(\mathbf H)^m_{k0n} + \Delta y  \sum_{l=1}^{s} b_{l}^{(2)}\mathscr{D}_{2}(\delta_y \mathbf H)^m_{kln},\nonumber \\
  &\mathscr{D}_{3}(\mathbf H)^m_{kl1} = \mathscr{D}_{3}(\mathbf H)^m_{kl0} + \Delta z  \sum_{n=1}^{s} b_{n}^{(3)}\mathscr{D}_{3}(\delta_z \mathbf H)^m_{kln}, 
\quad \mathscr{D}_{1}(\mathbf E)^m_{1ln} = \mathscr{D}_{1}(\mathbf E)^m_{0ln} + \Delta x  \sum_{k=1}^{s} b_{k}^{(4)}\mathscr{D}_{1}(\delta_x \mathbf E)^m_{kln}, \nonumber\\
  &\mathscr{D}_{2}(\mathbf E)^m_{k1n} = \mathscr{D}_{2}(\mathbf E)^m_{k0n} + \Delta y  \sum_{l=1}^{s} b_{l}^{(5)}\mathscr{D}_{2}(\delta_y \mathbf E)^m_{kln}, 
\quad \mathscr{D}_{3}(\mathbf E)^m_{kl1} = \mathscr{D}_{3}(\mathbf E)^m_{kl0} + \Delta z  \sum_{n=1}^{s} b_{n}^{(6)}\mathscr{D}_{3}(\delta_z \mathbf E)^m_{kln}, \nonumber\\
& (\mathbf H)^1_{kln} = (\mathbf H)^0_{kln} + \Delta t  \sum_{m=1}^{r} \tilde{b}_{m}^{(1)}\left( -\mathscr{D}_{1}(\delta_x \mathbf E)^m_{kln}  -\mathscr{D}_{2}(\delta_y \mathbf E)^m_{kln}  -\mathscr{D}_{3}(\delta_z \mathbf E)^m_{kln} \right)  + \lambda \Delta W_{kln}^{1}\sum_{m=1}^{r}\bar{b}_{m}^{(1)}(\mathbf E)^m_{kln}, \nonumber\\
& (\mathbf E)^1_{kln} = (\mathbf E)^0_{kln} + \Delta t  \sum_{m=1}^{r} \tilde{b}_{m}^{(2)}\left( \mathscr{D}_{1}(\delta_x \mathbf H)^m_{kln}  +\mathscr{D}_{2}(\delta_y \mathbf H)^m_{kln}  + \mathscr{D}_{3}(\delta_z \mathbf H)^m_{kln} \right)  - \lambda \Delta W_{kln}^{1}\sum_{m=1}^{r}\bar{b}_{m}^{(2)}(\mathbf H)^m_{kln},\nonumber
\end{align}
where
\begin{align*}
&\tilde{a}^{(2)}_{mi}\tilde{b}^{(1)}_{m} + \tilde{b}^{(2)}_{i}\tilde{a}^{(1)}_{im} - \tilde{b}^{(2)}_{i}\tilde{b}^{(1)}_{m} = 0,
\quad \bar{a}_{mi}^{(2)}\tilde{b}^{(1)}_{m} + \bar{b}_{i}^{(2)}\tilde{a}^{(1)}_{im} - \bar{b}_{i}^{(2)}\tilde{b}^{(1)}_{m} = 0,\\
&\tilde{a}^{(2)}_{mi}\bar{b}_{m}^{(1)} + \tilde{b}^{(2)}_{i}\bar{a}_{im}^{(1)}  - \tilde{b}^{(2)}_{i}\bar{b}_{m}^{(1)}  = 0,
~\quad \bar{a}_{mi}^{(2)}\bar{b}_{m}^{(1)} + \bar{b}_{i}^{(2)}\bar{a}_{im}^{(1)} - \bar{b}_{i}^{(2)}\bar{b}_{m}^{(1)} = 0,\\
& a^{(1)}_{kj}b^{(1)}_{k} + b^{(1)}_{j}a^{(1)}_{jk} - b^{(1)}_{j}b^{(1)}_{k} = 0,
\quad a^{(2)}_{lj}b^{(2)}_{l} + b^{(2)}_{j}a^{(2)}_{jl} - b^{(2)}_{j}b^{(2)}_{l} = 0,\\
& a^{(3)}_{nj}b^{(3)}_{n} + b^{(3)}_{j}a^{(3)}_{jn} - b^{(3)}_{j}b^{(3)}_{n} = 0,
\quad a^{(4)}_{kj}b^{(4)}_{k} + b^{(4)}_{j}a^{(4)}_{jk} - b^{(4)}_{j}b^{(4)}_{k} = 0,\\
& a^{(5)}_{lj}b^{(5)}_{l} + b^{(5)}_{j}a^{(5)}_{jl} - b^{(5)}_{j}b^{(5)}_{l} = 0,
\quad a^{(6)}_{nj}b^{(6)}_{n} + b^{(6)}_{j}a^{(6)}_{jn} - b^{(6)}_{j}b^{(6)}_{n} = 0,\\
&b^{(1)}_{k} =b^{(4)}_{k},~~b^{(2)}_{l} =b^{(5)}_{l},\qquad\qquad\quad  b^{(3)}_{n} =b^{(6)}_{n},~~ \tilde{b}^{(1)}_{m} =\tilde{b}^{(2)}_{m},
\end{align*}
${\mathbf  H}^{0}_{kln} \approx \mathbf  H(c^{(1)}_k\Delta x, c^{(2)}_l\Delta y, c^{(3)}_n\Delta z, 0)$, $ {\mathbf  H}^{m}_{kln} \approx \mathbf  H(c^{(1)}_k\Delta x, c^{(2)}_l\Delta y, c^{(3)}_n\Delta z, \tilde c^{(1)}_m\Delta t )$,  $ {\mathbf  H}^{1}_{kln} \approx \mathbf  H(c^{(1)}_k\Delta x, c^{(2)}_l\Delta y, c^{(3)}_n\Delta z,\Delta t )$,  $ {\mathbf  H}^{m}_{0ln} \approx \mathbf  H(0, c^{(2)}_l\Delta y, c^{(3)}_n\Delta z, \tilde c^{(1)}_m\Delta t )$, $ {\mathbf  H}^{m}_{1ln} \approx \mathbf  H(\Delta x, c^{(2)}_l\Delta y, c^{(3)}_n\Delta z, \tilde c^{(1)}_m\Delta t )$,  $ {\mathbf  H}^{m}_{k0n} \approx \mathbf  H(c^{(1)}_k\Delta x, 0, c^{(3)}_n\Delta z, \tilde c^{(1)}_m\Delta t )$, $ {\mathbf  H}^{m}_{k1n} \approx \mathbf  H(c^{(1)}_k\Delta x, \Delta y, c^{(3)}_n\Delta z,\tilde c^{(1)}_m\Delta t )$,  $ {\mathbf  H}^{m}_{kl0} \approx \mathbf  H(c^{(1)}_k\Delta x, c^{(2)}_l\Delta y, 0, \tilde c^{(1)}_m\Delta t )$,  ${\mathbf  H}^{m}_{kl1} \approx \mathbf  H(c^{(1)}_k\Delta x, c^{(2)}_l\Delta y, \Delta z, \tilde c^{(1)}_m\Delta t )$, etc., with $c^{(1)}_k =\sum_{j=1}^{s}a_{k j}^{(1)}$, $c^{(2)}_l =\sum_{j=1}^{s}a_{lj}^{(2)}$, $c^{(3)}_n =\sum_{j=1}^{s}a_{nj}^{(3)}$, $\tilde c^{(1)}_m =\sum_{n=1}^{r}\tilde{a}^{(1)}_{mn},$ and $ 1\leq j, k,l, n \leq s, ~ 1\leq i,m \leq r.$ Similar to the proof of Theorem \ref{thm7}, this fully-discrete method \eqref{MMM} satisfies the following discrete multi-symplectic conservation law
\begin{align*}
&\frac{1}{\Delta t }\sum_{k=1}^{s}\sum_{l=1}^{s}\sum_{n=1}^{s}b^{(1)}_{k}b^{(2)}_{l}b^{(3)}_{n}\left(\mathrm{d}(\mathbf E)^1_{kln}  \wedge \mathrm{d}(\mathbf H)^1_{kln}  - \mathrm{d}(\mathbf E)^0_{kln}  \wedge \mathrm{d}(\mathbf H)^0_{kln}  \right)\\
&+\frac{1}{2\Delta x }\sum_{m=1}^{r}\sum_{l=1}^{s}\sum_{n=1}^{s}\tilde{b}^{(1)}_{m}b^{(2)}_{l}b^{(3)}_{n}\left(\mathrm{d}(\mathbf U)^m_{1ln}  \wedge K_{1}\mathrm{d}(\mathbf U)^m_{1ln}  - \mathrm{d}(\mathbf U)^m_{0ln}  \wedge K_{1}\mathrm{d}(\mathbf U)^m_{0ln} \right)\\
&+\frac{1}{2\Delta y }\sum_{m=1}^{r}\sum_{k=1}^{s}\sum_{n=1}^{s}\tilde{b}^{(1)}_{m}b^{(1)}_{k}b^{(3)}_{n}\left(\mathrm{d}(\mathbf U)^m_{k1n}  \wedge K_{2}\mathrm{d}(\mathbf U)^m_{k1n}  - \mathrm{d}(\mathbf U)^m_{k0n}  \wedge K_{2}\mathrm{d}(\mathbf U)^m_{k0n} \right)\\
&+\frac{1}{2\Delta z }\sum_{m=1}^{r}\sum_{k=1}^{s}\sum_{l=1}^{s}\tilde{b}^{(1)}_{m}b^{(1)}_{k}b^{(2)}_{l}\left(\mathrm{d}(\mathbf U)^m_{kl1}  \wedge K_{3}\mathrm{d}(\mathbf U)^m_{kl1}  - \mathrm{d}(\mathbf U)^m_{kl0}  \wedge K_{3}\mathrm{d}(\mathbf U)^m_{kl0} \right) =0.
\end{align*} 

\section{Conclusions}
In this paper, three novel multi-symplectic methods are proposed to numerically solve stochastic Hamiltonian PDEs. 
We prove that the meshless LRBF collocation midpoint method, the splitting multi-symplectic Runge--Kutta method and the multi-symplectic partitioned Runge–Kutta method  preserve the discrete multi-symplectic conservation law almost surely. 
In general, these proposed multi-symplectic methods are always implicit, and have better numerical stability in the numerical implementation.  
Unlike the splitting multi-symplectic Runge--Kutta method and the multi-symplectic partitioned Runge–Kutta method, the meshless LRBF collocation midpoint method has high-order accuracy and does not require connection between nodes of the simulation domain, which leads to the liberty in selecting space nodes. 
Due to the geometric structure preserved property of the numerical method for subsystems, the splitting multi-symplectic Runge--Kutta method also has the superiority in preserving the averaged energy evolution law of some stochastic wave equations, as shown in Section 4.
The multi-symplectic partitioned Runge--Kutta method based on symplectic Euler method for some separate stochastic Hamiltonian PDEs, such as stochastic wave equations, is  always explicit, which reduces computational cost. 
 We take the stochastic wave equation as an example to perform numerical experiments, which indicates the validity of the proposed methods.
In fact, there are still many problems of interest which remain to be solved, such as 1)  to systematically construct explicit multi-symplectic methods for nonlinear stochastic Hamiltonian PDEs; 2) to propose numerical methods preserving both the multi-symplecticity and physical properties of stochastic Hamiltonian PDEs; 3) to prove theoretically the strong convergence order of accuracy for the proposed three numerical methods applied to stochastic Hamiltonian PDEs. 
We attempt to study these problems in our future work.

\section*{Acknowledgements}
This work is supported  by National key R\&D Program of China (No. 2020YFA0713701), and by National Natural Science Foundation of China (No. 11971470, No. 11871068, No. 12031020, No. 12022118, No. 12101596, No. 12171047), and by the China Postdoctoral Science Foundation (No. 2021M693339, No. 2021M690163, No. BX2021345).



\section*{References}
\begin{thebibliography}{}

\bibitem{Bouard}
A. de Bouard, A. Debussche.  On the stochastic Kortewegde Vries equation.  J. Funct. Anal., 154 (1998), 215-251.

\bibitem{Cao}
Y. Cao,  L. Yin,  Spectral Galerkin method for stochastic wave equations driven by space-time white noise. Commun. Pure  Appl. Anal.,  6 (2007),  607-617. 


\bibitem{Chen}
C. Chen, J. Hong,  L. Zhang. Preservation of physical properties of stochastic Maxwell equations with additive noise via multi-symplectic methods. J. Comput. Phys., 306 (2016), 500-519.

\bibitem{Cohen}
D.Cohen,  O.Verdier. Multisymplectic discretization of wave map equations. SIAM J. Sci. Comput.,  38 (2016), A953-A972.

\bibitem{Cui}
J. Cui,  J. Hong,  Z. Liu,  W. Zhou. Stochastic symplectic and multi-symplectic methods for nonlinear Schr\"odinger equation with white noise dispersion. J. Comput. Phys.,  342 (2017), 267-285.


\bibitem{Dehghan}
M. Dehghan, A. Shokri.  A numerical method for two-dimensional Schr\"odinger equation using collocation and radial basis functions. Comput. Math. Appl., 54 (2007), 136-146.

\bibitem{Hong1}
J. Hong, L. Ji,  L. Zhang. A multi-symplectic scheme for stochastic Maxwell equations 
with additive noise. J. Comput. Phys.,  268 (2014), 255-268. 

\bibitem{Hong2}
J. Hong,  L. Ji,  L. Zhang,  J. Cai.  An energy-conserving method for stochastic Maxwell equations 
with multiplicative noise. J. Comput. Phys., 351 (2017), 216-229.

\bibitem{Hong4}
J. Hong, H. Liu, G. Sun. The multi-symplecticity of partitioned Runge-Kutta methods for Hamiltonian PDEs. Math. Comp., 75 (2006), 167-181.

\bibitem{Hong3}
J. Hong,  X. Wang,  L. Zhang. Numerical analysis on ergodic limit of approximations for stochastic NLS equation via multi-symplectic scheme. SIAM J. Numer. Anal., 55 (2017), 305-327.

\bibitem{Jiang}
S. Jiang, L. Wang, J. Hong.  Stochastic multi-symplectic integrator for stochastic nonlinear Schr\"odinger equation. 
Commun. Comput. Phys., 14 (2013), 393-411.

\bibitem{Kansa1}
E. Kansa. Multiquadrics-A scattered data approximation scheme with applications to computational fluid 
dynamics- I Surface approximations and partial derivative estimates. Comput. Math. Appl., 19 (1990), 127-145.

\bibitem{Kansa2}
E. Kansa. Multiquadrics-A scattered data approximation scheme with applications to computational fluid 
dynamics- II. Solutions to parabolic, hyperbolic and elliptic partial differential equations. Comput. Math. Appl., 19 (1990), 147-161.

\bibitem{Lee}
C. K. Lee,  X. Liu, S. C. Fan.  Local multiquadric approximation for solving boundary value problems. Comput. Math.,  30 (2003), 396-409.


\bibitem{McLachlan}
R. I. McLachlan, B. N. Ryland,  Y. Sun.  High order multisymplectic Runge-Kutta methods. SIAM J. Sci. Comput. 36 (2014), A2199-A2226.

\bibitem{Roach}
G. Roach, I. Stratis, A. Yannacopoulos.  Mathematical analysis of deterministic and stochastic problems in complex media electromagnetics, Princeton University Press, 2012.

\bibitem{Sarler}
B.  $\rm \check{S}$arler,  R. Vertnik.  Meshless explicit local radial basis function collocation method for diffusion
problems. Comput. Math. Appl.,  51 (2006), 1269-1282.

\bibitem{Song}
M. Song, X. Qian, T. Shen, S. Song. Stochastic conformal schemes for damped stochastic Klein-Gordon equation with additive noise.  J. Comput. Phys., 411 (2020), 109300.

\bibitem{Wu}
Z. Wu, S. Zhang.  A meshless symplectic algorithm for multi-variate Hamiltonian PDEs with radial basis approximation. Eng. Anal. Bound. Elem.,  50 (2015), 258-264. 

\bibitem{Zhang}
S. Zhang. Meshless symplectic and multi-symplectic local RBF
collocation methods for nonlinear Schr\"odinger equation. J. Comput. Phys., 450 (2022), 110820.

\bibitem{ZhangJ}
L. Zhang, L. Ji. Stochastic multi-symplectic Runge-Kutta methods for stochastic Hamiltonian PDEs. Appl. Numer. Math., 135 (2019), 396-406.


\end{thebibliography}
\end{document}